\documentclass[12pt,reqno]{amsart}

\usepackage[a4paper,left=24mm,right=24mm,top=32mm,bottom=30mm]{geometry}
\usepackage{graphicx}
\usepackage{pgf}
\usepackage{import}
\usepackage{amssymb}
\usepackage{amsthm,mathtools}
\usepackage{amsmath,amsfonts,amssymb}
\usepackage{dsfont}
\usepackage{todonotes}
\usepackage{hyperref}
\usepackage{cleveref}
\usepackage[english]{babel}
\usepackage{pdfpages}
\usepackage{bm}
\usepackage[utf8]{inputenc}
\usepackage[mathscr]{euscript}
\usepackage{caption}
\usepackage{subcaption}
\usepackage{cancel}
\usepackage{stackrel}
\usepackage{textcomp}
\usepackage[T1]{fontenc}

\usepackage{commath} 

\usepackage{tikz}
\usetikzlibrary{arrows,topaths}
\usetikzlibrary{positioning}

\DeclareMathOperator{\Id}{Id}
\DeclareMathOperator{\proj}{proj}
\DeclareMathOperator*{\esssup}{ess\,sup}
\newcommand{\Z}{\mathbb{Z}}
\newcommand{\N}{\mathbb{N}}
\newcommand{\R}{\mathbb{R}}
\newcommand{\T}{\mathbb{T}}

\newcommand{\half}{\frac{1}{2}}
\newcommand{\lra}{\longrightarrow}

\newcommand{\ssubset}{\subset\mathrel{\mkern-3mu}\subset}
\newcommand{\FF}{\mathrm{I}\mathrel{\mathchoice{\mkern-12mu}{\mkern-15mu}{}{}}\mathrm{I}}

\newcommand{\WCH}{W_{\mathrel{\mathchoice{}{}{\mkern-2mu}{}}\mathrm{o}}}
\newcommand{\uCH}{u^\mathrm{o}}
\newcommand{\lambdaCH}{\lambda^\mathrm{o}}

\newcommand{\calA}{\mathscr{A}}
\newcommand{\scrA}{\mathscr{A}}
\newcommand{\calP}{\mathcal{P}}
\newcommand{\Chi}{\mathcal{X}}
\newcommand{\calE}{\mathcal{E}}
\newcommand{\Per}{\mathcal{P}}
\newcommand{\calW}{\mathcal{W}}
\newcommand{\calG}{\mathcal{G}}
\newcommand{\calL}{\mathcal{L}}
\newcommand{\calM}{\mathcal{M}}

\newcommand{\calH}{\mathcal{H}}
\newcommand{\calV}{\mathcal{V}}
\newcommand{\calO}{O}
\newcommand{\bigO}{\calO}
\newcommand{\eps}{\varepsilon}
\newcommand{\LL}{\mathcal{L}}
\newcommand{\Ha}{\mathcal{H}}
\newcommand{\mres}{\mathbin{\vrule height 1.6ex depth 0pt width
0.13ex\vrule height 0.13ex depth 0pt width 1.3ex}}

\newcommand{\fraku}{\bm{u}}
\newcommand{\frakv}{\bm{v}}
\newcommand{\frakL}{\bm{L}}
\newcommand{\frakA}{\bm{\mathscr{A}}}
\renewcommand{\d}{\,{\operatorname{d}}}
\newcommand{\dx}{\:\!\d\calL}
\newcommand{\ds}{\:\!\d\calH}

\newcommand\blfootnote[1]{%
  \begingroup
  \renewcommand\thefootnote{}\footnote{#1}%
  \addtocounter{footnote}{-1}%
  \endgroup
}

\DeclareMathOperator{\sdist}{sdist}
\DeclareMathOperator{\loc}{loc}
\DeclareMathOperator{\spann}{span}
\DeclareMathOperator{\dist}{dist}
\DeclareMathOperator{\sgn}{sgn}

\DeclareMathOperator{\range}{range}
\DeclareMathOperator{\inn}{in}
\newcommand{\uei}{u_\eps^{\inn}}
\newcommand{\vei}{v_\eps^{\inn}}
\newcommand{\wei}{w_\eps^{\inn}}

\theoremstyle{plain}
\newtheorem{theorem}{Theorem}[section]
\newtheorem{corollary}[theorem]{Corollary}
\newtheorem{lemma}[theorem]{Lemma}
\newtheorem{proposition}[theorem]{Proposition}
\theoremstyle{definition}
\newtheorem{remark}[theorem]{Remark}

\newtheorem{assumptions}[theorem]{Assumption}
\newtheorem{definition}[theorem]{Definition}
\newtheorem{prelim}[theorem]{Preliminaries}

\begin{document}

\title[Gradient-free diffuse approximation of the Willmore functional]{\textquotedbl{}Gradient-free\textquotedbl{} diffuse approximations of the Willmore functional and Willmore flow}

\author{Nils Dabrock,
Sascha Knüttel*
\and
Matthias R\"oger
}

\begin{abstract}
We introduce new diffuse approximations of the Willmore functional and the Willmore flow. They are based on a corresponding approximation of the perimeter that has been studied by Amstutz-van Goethem [{\em Interfaces Free Bound. 14 (2012)}]. We identify the candidate for the $\Gamma$--convergence, prove the $\Gamma$--limsup statement and justify the convergence to the Willmore flow by an asymptotic expansion. Furthermore, we present numerical simulations that are based on the new approximation.
\end{abstract}
\subjclass[2010]{35R35, 35K65, 65N30}
\keywords{Willmore flow, phase-field model, diffuse
interface, sharp interface limit}

\maketitle
\blfootnote{All authors are affiliated with Technische Universit\"at Dortmund, Fakult\"at f\"ur Mathematik,\\
			\hspace*{0.43cm}Vogelpothsweg 87, D-44227 Dortmund, Germany}
\blfootnote{*Corresponding author. Tel.: +49 231 755 5163, Fax:+49 231 755 5942,\\
			\hspace*{0.43cm}E-Mail: Sascha.Knuettel@math.tu-dortmund.de}
\section{Introduction}
\setlength{\parindent}{0cm}
The Willmore functional
\begin{equation*}
  \mathcal{W}(\Gamma) \coloneqq \int_{\Gamma} H^2(y)\d\mathcal{H}^{n-1}(y),
\end{equation*}
of a $C^2$-regular hypersurface $\Gamma\subseteq\R^n$ with mean curvature $H$ is one of the most prominent examples of a curvature energy.
Such energies have already been considered by Poisson \cite{Poisson1814} and Germain \cite{Germain1821} in the $19^{\text{th}}$ century and appear in a variety of applications, for example as a shape energy of bio membranes as proposed by Canham \cite{Canham1970} and Helfrich \cite{Helfrich1973}.
The Willmore functional in particular has been studied intensively in differential geometry and geometric measure theory by Thomsen and Blaschke \cite{Thomsen1924,BlaschkeThomsem1929} at the beginning of the last century as well as by Willmore \cite{Willmore1965} and more recently by Simon \cite{Simon1993}, Kuwert and Schätzle \cite{KuwertSchaetzle2012}, Riviére \cite{Riviere2008}.
The most spectacular recent contribution is the proof of the Willmore conjecture on the minimal Willmore energy of immersed tori by Marques and Neves \cite{MarquesNeves2014}.

In the case of planar curves, the Willmore energy reduces to Euler's Elastica energy for the bending of a rod.
This energy has an even longer history (see e.g. \cite{Love2013}), has been thoroughly investigated in many contributions \cite{LangerSinger1984, Mumford1994} and still is an active field of research.

Gradient flows for curvature energies as steepest decent dynamics have also attracted a lot of attention.
The Willmore flow, in particular, has been considered in many contributions over the past decades, see for example Simonett \cite{Simonett2001} and Kuwert and Schätzle \cite{KuwertSchaetzle2001,KuwertSchaetzle2002,KuwertSchaetzle2004}.

\medskip
Motivated by phase separation problems and as a tool for numerical simulations, diffuse approximations of curvature energies and in particular the Willmore functional and Willmore flow are widely used.
The most famous example is the phase field approximation going back to De Giorgi \cite{DeGiorgi1991}.
This approximation is based on the Van der Waals--Cahn-Hilliard energy,
given by
\begin{equation}
  \calP^{\text{CH}}_\eps(u) \coloneqq \int_{\Omega} \Big(\frac{\eps}{2} |\nabla u|^2 +
  \frac{1}{\eps}W(u)\Big)\dif\LL^n,
  \label{eq:Pd}
\end{equation}
where $W$ is a suitable double well potential and $u$ is a smooth function on a domain $\Omega\subseteq\R^n$.
To achieve low energy values, the function $u$ has to be close to the wells of the potential except
for thin transition layers with thickness of order $\eps$.
The celebrated result by Modica and Mortola \cite{ModicaMortola1977,Modica1987}
states that this functionals $\Gamma$--converge in the sharp interface
limit $\eps\to 0$ to the perimeter functional $\calP$,
\begin{equation}
  \calP^{\text{CH}}_\eps \to c_0^{\text{CH}} \calP,\quad
  c_0^{\text{CH}} = \int_{-1}^1 \sqrt{2W}\dx^1.
  \label{eq:c0}
\end{equation}
Since the mean curvature is the $L^2$-gradient of the perimeter it is a natural approach to take the $L^2$-gradient
of the diffuse perimeter \eqref{eq:Pd} as a starting point for a diffuse Willmore energy. This
motivates a formal approximation of the Willmore functional, given by
\begin{align}
  \calW^{\text{CH}}_\eps(u) \,&\coloneqq\, \int_{\Omega} \frac{1}{2\eps} \Big(-\eps\Delta u +
  \frac{1}{\eps}W'(u)\Big)^2\dx^n, \label{eq:wd}
\end{align}
which is a modified version of De Giorgi's proposal \cite{DeGiorgi1991}, introduced by Bellettini and Paolini in \cite{BellettiniPaolini1993}.
The $\eps^{-1}$-factor compensates for the volume of the transition layer as we will see in the calculations in chapter \ref{sec:3}.\\
The $\Gamma$--$\limsup$ property was proved in \cite{BellettiniPaolini1993}, the $\Gamma$--convergence
was shown in dimensions 2 and 3 for smooth limit configurations in \cite{RoegerSchaetzle2006}.
This gives a solid justification for this approximation, though many issues concerning the
$\Gamma$--convergence are still to be resolved, in particular concerning non-smooth limit configurations.

Corresponding diffuse approximations of the Willmore flow have been introduced by \cite{MarchDozio1997}
and have been justified by formal asymptotic expansions in \cite{LoretiMarch2000}, see also \cite{Wang2008}.
In a recent article \cite{FeiLiu2021} Fei and Liu prove the convergence of diffuse approximations to the Willmore
flow for well-prepared initial data, as long as the smooth limit flow exists.

Quite a number of numerical schemes for the simulation of the Willmore flow have been proposed. For the treatment
in a sharp-interface approach, we refer to
\cite{BarrettGarckeNuernberg2010,BarrettGarckeNuernberg2008,BarrettGarckeNuernberg2007,
Rusu2005,DeckelnickDziukElliott2005,BonitoNochettoPauletti2010,ElliottStinner2010}. Level set techniques have been used in
\cite{DroskeRumpf2004,BretinMasnouOudet2015}.
Diffuse approximation have been employed in a huge number
of applications, see for example
\cite{DuLiuWang2004,CampeloHernandez-Machado2006,WangDu2008,LowengrubRaetzVoigt2009,DondlMugnaiRoeger2011,
DondlLemenantWojtowytsch2017,EsedogluRaetzRoeger2014,WangJuDu2016,
FrankenRumpfWirth2013,BretinMasnouOudet2015}.

\medskip
A number of alternative diffuse approximations of the Willmore energy have been proposed, in particular to enforce
the $\Gamma$--convergence of approximations for non-smooth limit configurations to the $L^1$ lower semi-continuous
envelope of the Willmore functional \cite{Bellettini1997,Mugnai2013,EsedogluRaetzRoeger2014,RaetzRoeger2021}.
These approximations, however, often lack the simplicity of the standard approximation and its direct relation to applications.

\medskip
A class of nonlocal perimeter approximations can be derived from classical Ising-type models \cite{AlbertiBellettini1998}.
They involve a discrete gradient and a double integral
\begin{equation}
  \calP^{\text{AB}}_\eps(u) = \frac{1}{\eps}\int_\Omega W(u)\dx^n + \frac{\eps}{4}\int_\Omega
  \int_\Omega J_\eps(x-y)\Big(\frac{u(x)-u(y)}{\eps}\Big)^2\d y\d x,
  \label{eq:P_AB}
\end{equation}
where $J_\eps=\eps^{-n}J(\cdot/\eps)$ for a suitable kernel $J$. Alberti and Bellettini \cite{AlbertiBellettini1998} proved the
$\Gamma$--convergence with $\eps\to 0$ to an (in general anisotropic) perimeter functional.
At least on a formal level, one might construct diffuse approximations of the Willmore energy starting from the $L^2$-gradient of $\calP_\eps^{\text{AB}}$.
However, to the best of our knowledge this has not been addressed yet and a rigorous justification in a general framework seems to be difficult.
We refer to Braides \cite{Braides2002} for some more details on the above mentioned models and an overview over different perimeter approximations.

In the present paper we consider an approximation that is in between the Cahn--Hilliard model and the Ising-type models just described.
It is motivated by an in general anisotropic two-variable energy studied by Solci and Vitali in \cite{SolciVitali2003}.
In the isotropic case the functional is characterized as
\begin{equation*}
  \calG_\eps(u,v) := \int_\Omega\Big(\frac{\eps}{2}|\nabla v|^2 + \frac{1}{2\eps}(u-v)^2 + \frac{1}{2\eps}W(u)\Big)\dx^n.
\end{equation*}
Such an energy does also appear in a one-dimensional model for the longitudinal
deformation of an elastic bar, proposed by Rogers and Truskinovsky
\cite{RogersTruskinovsky1997}, and can be connected to certain two-variable models
for phase separation processes, see the references in \cite{AmstutzvanGoethem2012}.
Solci and Vitali proved that $\calG_\eps$ does $\Gamma$--converge on $L^1(\Omega)^2$
towards a functional $\calG$ that is only finite on $\{(u,u)\in BV(\Omega;\{\pm 1\})^2\}$
with $\calG(u,u)=c_1\calP(u)$ for some $c_1>0$.

We follow here the analysis of Amstutz and van Goethem \cite{AmstutzvanGoethem2012} of
a gradient-free approximation of the perimeter functional that is obtained by considering
the marginal functional of $\calG_\eps$, where for given $u$ the variable $v$ is chosen as
minimizer of $\calG_\eps(u,\cdot)$. This leads to a representation of the optimal
$v_\eps=v_\eps[u]$ as solution of
\begin{align}
 -\eps^2\Delta v_\eps +v_\eps = u\quad\text{ in }\Omega,\qquad
 \nabla v_\eps\cdot\nu_\Omega = 0 \quad\text{ on }\partial\Omega,
 \label{eq:v_eps}
\end{align}
and to the functional
\begin{align}
  \calP^{\text{AG}}_\eps(u) &:= \inf\limits_{v\in H^1(\Omega)}\int_\Omega\Big( \frac{\eps}{2}|\nabla v|^2
  + \frac{1}{2\eps}(u-v)^2 + \frac{1}{2\eps}W(u)\Big)\dx^n\nonumber\\
  &=\int_\Omega\Big( \frac{\eps}{2}|\nabla v_\eps|^2
  + \frac{1}{2\eps}(u-v_\eps)^2 + \frac{1}{2\eps}W(u)\Big)\dx^n\nonumber\\
  &= \int_\Omega\frac{1}{2\eps}\Big(u(u-v_\eps) + W(u)\Big)\dx^n
  \label{eq:defdiffuseperimeterausminimiert}.
\end{align}

For the particular choice $W(r)=1-r^2$ in $[-1,1]$ and  locally constant linear growth outside $[-1,1]$ the $\Gamma$--convergence in $L^1(\Omega)$
towards a multiple of the perimeter functional was proved in \cite{AmstutzvanGoethem2012}.
This can be generalized to a bigger class of double-well potentials and in particular to the
potentials used below. For our analysis, however, it is important to have some smoothness and
quadratic behavior of $W$ in the wells, see \ref{rem:assumptionsW} below. This condition, on
the other hand, excludes the double well potential used in \cite{AmstutzvanGoethem2012}.

In \cite{AmstutzvanGoethem2012} numerical simulations of some topology optimization
problems were presented, where the gradient-free structure of the functional
(with respect to the variable $u$) proved to be advantageous.
We note that in the case $\Omega=\R^n$ the approximation $\calP^{\text{AG}}_\eps$ corresponds to
$\calP^{\text{AB}}_\eps$ in \eqref{eq:P_AB}, with a particular choice of $J$.

The solution operator induced by the PDE \eqref{eq:v_eps} is linear and self-adjoint.
The $L^2$-gradient of $\calP^{\text{AG}}_\eps$ therefore is given by
\begin{equation}
 H_\eps := \nabla_{L^2}\calP_\eps^{\text{AG}}(u)=\frac{1}{\eps}\big(u + \frac{1}{2}W'(u)-v_\eps\big).
 \label{eq:defH_eps}
\end{equation}
In analogy to the sharp interface situation $H_\eps$ can be seen as a {\em diffuse mean curvature}.
This suggests the formal Willmore energy approximation
\begin{align}
  \calW_\eps(u) \coloneqq\calW^{\text{AG}}_\eps(u)\coloneqq\int_\Omega
  \frac{1}{\eps^3}\Big(u + \frac{1}{2}W'(u)-v_\eps \Big)^2\dx^n,
    \label{eq:defwillmoreapproxag}
\end{align}
which is the main object of the current study. The additional factor $\eps^{-1}$ in $\calW_\eps$
again accounts for the small volume of the transition layer region. We remark that no gradients appear
(explicitly) in the functional and only the rather well-behaved solution operator $u\mapsto v_\eps$
associated to \eqref{eq:v_eps} enters the energy. This makes the above functional an interesting
candidate for numerical simulations.

Besides the theoretical interest in Willmore approximations
and its use in numerical simulations, the analysis of the functional $\calW_\eps$ is also important
for the understanding of the corresponding perimeter approximation \eqref{eq:defdiffuseperimeterausminimiert},
and in particular the associated $L^2$-gradient flow. In fact, the
$\Gamma$--convergence of $\calW_\eps$ is one of the properties
necessary to apply a general result
about convergence of gradient flows proved by Sandier
and Serfaty in \cite{SandierSerfaty2004}.

\medskip
The function $v_\eps$ that appears in the functional $\calP^{\text{AG}}_\eps$ and that is characterized by \eqref{eq:v_eps} represents a particular regularization of $u$.
Perimeter approximations based on other regularizations are possible as well (also the functional $\calP^{\text{AB}}_\eps$ can be represented this way) and the results in \cite{SolciVitali2003} could be used to prove its $\Gamma$--convergence, which also was central for the proof in \cite{AmstutzvanGoethem2012}.
This becomes less
clear when dealing with approximations of the Willmore energy.
Here we exploit the very convenient PDE characterization of $v_\eps$.
We suspect that the $\Gamma$--convergence towards the Willmore functional can be proved also for diffuse approximations based on other regularizations of $u$.
However, to the best of our knowledge no such results are currently available.
Developing a general theory therefore might be an interesting field for future research.

\medskip
Besides the static functionals we are also interested in $L^2$-type gradient flows.
In the sharp interface setting we therefore consider evolutions of phases
$(E(t))_{t\in(0,T)}$ and of the associated boundaries $\Gamma(t)=\partial E(t)$.
In case of the perimeter functional the formal $L^2$-gradient flow is given by
the mean curvature flow
\begin{equation}
	\calV= H, \label{eq:mcf}
\end{equation}
where $H$ and $\calV$ denote the scalar mean curvature and normal velocity of the evolution
in direction of the inner unit normal field associated to $E(t)$, $t\in (0,T)$.
Mean curvature flow is one of the most prominent geometric flows and has been studied
extensively of the past decades. We refer to \cite{Bellettini2013} for a proper introduction
to the subject.

\medskip
The Willmore flow is the formal $L^2$-gradient flow of $\calW$ and is given by
\begin{equation}
	\calV = -\Delta_{\Gamma}H+\half H^3-H|\FF|^2,
	\label{eq:willmoreflow}
\end{equation}
where $\displaystyle{\FF}(\cdot,t)$ denotes the second fundamental form of $\Gamma(t)$ and $\Delta_\Gamma$ the Laplace-Beltrami operator on $\Gamma(t)$.
The Willmore flow is a fourth order geometric evolution law, which introduces quite some additional challenges for the analysis of the flow.
We refer to the already mentioned fundamental contributions \cite{Simonett2001}, \cite{KuwertSchaetzle2001,KuwertSchaetzle2002,KuwertSchaetzle2004}.
For a derivation of the formula for the $L^2$-gradient of the Willmore functional see \cite[sections 7.4 - 7.5]{Willmore1993}.

\medskip
In an analogue way we associate $L^2$-gradient flows to the diffuse perimeter  $\calP^{\text{AG}}_\eps$ and Willmore energy approximations $\calW_\eps$ in \eqref{eq:defdiffuseperimeterausminimiert}, \eqref{eq:defwillmoreapproxag}.

The $L^2$-gradient of $\calP^{\text{AG}}_\eps$ has already been characterized in \eqref{eq:defH_eps} by $H_\eps=\frac{1}{\eps}\big(-v_\eps + u + \frac{1}{2}W'(u)\big)$.
Taking the variational derivative of $\calW_\eps$ we find
\begin{equation}
  \nabla_{L^2}\calW_\eps(u)=
  \frac{2}{\eps^2}\big(1+\frac{1}{2}W''(u)-(\Id-\eps^2\Delta)^{-1}\big)H_\eps.
\end{equation}
Appropriately rescaled this leads to the {\em diffuse mean curvature flow}
\begin{align}
 \eps\partial_t u_\eps = H_\eps,
 \label{eq:appMCF}
\end{align}
and the {\em diffuse Willmore flow}
\begin{align}
 \eps\partial_t u_\eps=
 -\frac{2}{\eps^2}\big(1+\frac{1}{2}W''(u)-(\Id-\eps^2\Delta)^{-1}\big)H_\eps.
 \label{eq:appWillflow}
\end{align}
We may expect that the diffuse flows converge in the sharp interface limit $\eps\to 0$ to mean curvature flow and Willmore flow, respectively.

\medskip
The goal of this paper is to provide some justification to the above mentioned formal
approximation properties. In particular, we will identify the candidate
for the $\Gamma$--convergence of the functionals $\calW_\eps$, which in
fact is proportional to the Willmore functional, with a specific
constant of proportionality that only depends on the choice of the
double well potential. For this
candidate we prove the corresponding $\Gamma$--limsup construction.
In addition, we give a rigorous lower bound in particular classes of phase field approximations that are described by suitable expansion properties.
Moreover, we justify by a formal
asymptotic expansion the convergence of the flow \eqref{eq:appWillflow}.
For the proof we basically follow the approach already used by Loreti and March
\cite{LoretiMarch2000} and Wang \cite{Wang2008}.
However, the operators that
define the gradient-free approximation are different from the standard
case and the derivation of the convergence property is much more involved,
in particular for the case of the approximate Willmore flow.

We finally use our approach for numerical simulations. We follow the implicit spectral
discretization scheme proposed by \cite{BretinMasnouOudet2015} for the standard diffuse
approximation and compare the new and the standard scheme.

\begin{prelim}\ \\
We collect some basic notations, definitions and assumptions that
we will use throughout the paper.

Let $n\in\N$ and $\Omega\subseteq\R^n$ be a bounded open domain with Lipschitz boundary.
We denote by $\calL^n$ the $n$-dimensional Lebesgue measure and by $\Ha^{n-1}$ the
$(n-1)$-dimensional Hausdorff measure. The space $BV(\Omega;\{\pm 1\})$ consists of
all function of bounded variation with values in $\{\pm 1\}$ almost everywhere.
$BV(\Omega;\{\pm 1\})$ can be identified with the sets of finite perimeter in $\Omega$,
where we associate to a set $E$ of finite perimeter the rescaled characteristic function
$u\coloneqq2\Chi_E-1$. The essential boundary of a set $E\subseteq\Omega$ with finite
perimeter in $\Omega$ is denoted by $\partial^*E$. For $u=2\Chi_E-1$ we then have
$|\nabla u|=2\Ha^{n-1}\mres \Gamma$.

With these notations we define the perimeter functional $\calP_\Omega:
L^1(\Omega)\longrightarrow[0,\infty]$ by
\begin{align}
	\calP_\Omega(u)\coloneqq
	\left\{
		\begin{array}{ll}\displaystyle
			\frac{1}{2}\int_\Omega |\nabla u|\dx^n,
			&\textrm{ if }u\in BV(\Omega,\{\pm1 \})\\\\
			+\infty,&\textrm{ else}.
		\end{array}
	\right.\label{eq:perimeterdef}
\end{align}
Note $\calP_\Omega(2\Chi_E-1)=\Ha^{n-1}(\partial^*E\cap\Omega)$ for $E\subseteq\Omega$
with finite perimeter.

Furthermore we define the Willmore functional
$\calW:L^1(\Omega)\longrightarrow[0,\infty]$,
\begin{align}
	\calW_\Omega(u)\coloneqq
		\begin{cases}\displaystyle
			\int_{\partial E\cap\Omega}H^2\ds^{n-1}, &\text{ if }u=2\Chi_E-1\text{ for }
			E\subseteq\Omega\text{ with }\partial E\cap\Omega\in C^2\\
			+\infty,&\text{ else,}\label{eq:willmoredef}
		\end{cases}
\end{align}
where the mean curvature $H$ is defined as the sum of the principle curvatures of
$\partial E$ (taken positive for convex $E$). We will drop the $\Omega$-index for
$\calP$ and $\calW$ as $\Omega$ is fixed.
When we use a Landau symbol $\calO$ in an
equation with multiple variables it will always mean, that the respective term is
uniformly bounded with respect to all variables.
If we need more precise information,
such as an exponential decay we will use a different notation, often $R_\eps$, and state the additional requirements.
As long as
the term keeps its properties we will not necessarily change notation and still write $R_\eps$
even though the term might have changed, just as in the Landau notation.
The same convention will be used for constants $\Lambda>0$. They might change from line
to line, however they will always be uniform in $\eps$.
\end{prelim}

\begin{assumptions}[on the double well potential $W$]\label{rem:assumptionsW}\ \\
To describe the diffuse approximations we fix a double well potential $W\in C^m(\R)$,
$m\in\N^{\geq 4}$, with $W\geq 0$, $\{W=0\}=\{\pm 1\}$, $W''(\pm 1)>0$, $1+\frac{1}{2}W''>0$
in $[-1,1]$ and at least linear growth at $\pm\infty$.
We associate to $W$ the mapping $f:\R\lra\R$, $f(r)\coloneqq r+\half W'(r)$
that appears in the computation of the $L^2$-gradient of $\Per_\eps$ above.
We have $f'(r)=1+\frac{1}{2}W''(r)\geq\gamma$ for some $\gamma>0$ and all $r\in[-1,1]$.
This implies that $f\in C^{m-1}([-1,1])$ is strictly increasing, we further obtain that
$f:[-1,1]\to [-1,1]$ is one-to-one and that $f$ has an inverse function
$f^{-1}\in C^{m-1}([-1,1])$, with
\begin{align}
	(f^{-1})'\leq\frac{1}{\gamma}
	\text{ in } [-1,1].
	\label{eq:finverseableiten}
\end{align}
If $W$ is an even function then $f$ and $f^{-1}$ are odd. This conditions cover
a large class of admissible double well potentials, such as the standard quartic
double well potential $W(r)=\frac{1}{4}(1-r^2)^2$ that is most often used in simulations.
See Section \ref{sec:numerics} for another convenient choice.

On the other hand the particular choice
$W(r)=1-r^2$, $r\in [-1,1]$, with locally constant linear growth outside $[-1,1]$ in \cite{AmstutzvanGoethem2012} is not allowed, since in this case $f$ would be constant in $(-1,1)$ and is not $C^1$-regular on $\R$.
\end{assumptions}

\section{One space dimension: The optimal profile problem}
For many diffuse approximations the study of the optimal transition between the pure phases on the real line is key for the understanding of its behavior, see for example \cite{Alberti2000} and the references therein.
As for the Cahn--Hilliard approximation $\calP_\eps^{\text{CH}}$ and Willmore
functional $\calW_\eps^{\text{CH}}$ we expect that typical small-energy configurations
for $\calP_\eps^{\text{AG}}$ and $\calW_\eps^{\text{AG}}$ are quasi one-dimensional and
can be constructed from an optimal transition profile and the rescaled signed distance
from the zero-level set.

To characterize the optimal profile associated to $\calP_\eps^{\text{AG}}$ we consider the
following minimization problem on the real line. We fix a suitable class of real functions
\begin{align*}
	\calM\coloneqq \{u\in L^\infty(\R)\,:\,
	\limsup_{r\to-\infty}u(r)<0,\,\liminf_{r\to\infty}u(r)>0\}.
\end{align*}
Moreover we define for $u\in\calM$, $v\in H^1_{\loc}(\R)$ with $\lim_{x\to\pm\infty}
v(x)=\pm 1$ the energies
\begin{align*}
	\calG_\eps^1(u,v) \coloneqq
	\int_\R \half \big(\eps(v')^2+\frac{1}{\eps}(u-v)^2+\frac{1}{\eps}W(u)\big)\dx^1,\\
	\calE_\eps(u) \coloneqq \inf \Big\{
	\calG_\eps^1(u,v)\,:\, v\in H^1_{\loc}(\R),\,\lim_{r\to\pm\infty}v(r)=\pm 1\Big\}.
\end{align*}
By rescaling we see that the minimization problems can be reduced to the case
$\eps=1$ and we write $\calE=\calE_1$ in what follows.

\begin{theorem}[Optimal profile]\label{thm:1dminimizer}
Consider $W$ as in Remark \ref{rem:assumptionsW}.
There exists a unique minimizer $\fraku_0\in \calM$ of $\calE$ that satisfies $u(0)=0$.
This minimizer satisfies $\fraku_0\in C^{m-1}(\R;(-1,1))$ and is determined by
\begin{align}
	\fraku_0=f^{-1}(\frakv_0),
	\label{eq:uv1dimbeziehung2}
\end{align}
where $\frakv_0\in C^{m+1}(\R;(-1,1))$ is the unique solution to
 \begin{equation}
	\frakv_0' = \sqrt{W\big(f^{-1}(\frakv_0)\big)+\frac{1}{4}
	W'\big(f^{-1}(\frakv_0)\big)^2}\text{ with }\frakv_0(0)=f(0).
	\label{eq:profileDGL}
 \end{equation}
We also have $\fraku_0(0)=0$, $\fraku_0'>0$, $\frakv_0'>0$ and
\begin{align}
	-\frakv_0'' + \frakv_0 = \fraku_0
	\quad\text{in}\quad\R.
	\label{eq:uv1dimbeziehung}
\end{align}
The functions $\fraku_0,\frakv_0$ converge exponentially fast to $\pm 1$ as
$x\to\pm\infty$. All derivatives of $\frakv_0$ up to order $m+1$ and all derivatives
of $\fraku_0$ up to order $m-1$ decay exponentially at $\pm\infty$.
\end{theorem}

\begin{remark}
We remark that we consider a priori functions with values in $\R$.
In the proof below that we can also restrict the minimization to functions with values
in $[-1,1]$ and obtain that the optimal profile takes it values only in $(-1,1)$.
Since the diffuse Willmore flow is of fourth order and does not satisfy a maximum principle,
we cannot guarantee that evolutions take values only in $(-1,1)$.
In particular, the behavior of the double well potential on $\R$ matters for the analysis below.
\end{remark}

\begin{proof}
We can project the values of any competitor onto $[-1,1]$, which will lower the value of all three
summands in the integral.
By the assumptions on $W$, for given value $v\in [-1,1]$ the
real function $r\mapsto (r-v)^2+W(r)$ has a unique strict minimum at $u=f^{-1}(v)$. We
therefore deduce that
\begin{align*}
	\calE(u) \geq \inf \Big\{
	\int_\R \half \big((v')^2+\frac{1}{4}W'(f^{-1}(v))^2+ W(f^{-1}(v))\big)\dx^1\,:\,\\
	v\in H^1_{\loc}(\R),\,\lim_{r\to\pm\infty}v(r)=\pm 1\Big\}.
\end{align*}
The right-hand side takes the form of the optimal profile problem for a standard
Cahn--Hilliard functional with double well potential
$W_\ast(r) = W(f^{-1}(r))+\frac{1}{4}W'(f^{-1}(r))^2$. We have
$W_\ast\in C^{m-1}([-1,1])$ and it satisfies $W_\ast\geq 0$ as well as
$\{W_\ast=0\}=\{\pm 1\}$.
For this problem it is well-known that a unique optimal profile exists and is given
by the solution to the ODE \eqref{eq:profileDGL},
see for example \cite[Section 3a]{Alberti2000}. Setting $\fraku_0=f^{-1}(\frakv_0)$
therefore achieves the optimal value of $\calE$ in $\calM$.

The regularity of $\fraku_0,\frakv_0$ and the exponential decay
properties follow from the combination of
the ODE \eqref{eq:profileDGL}, $\fraku_0=f^{-1}(\frakv_0)$ and standard ODE theory.
\end{proof}

To identify the candidate for the $\Gamma$--limit of the energies $\calP_\eps$,
$\calW_\eps$ we need to compute the energy of the optimal profile.

\begin{corollary}[Double well potential-depending constants]\ \\
The constants
\begin{align}
 	c_0\coloneqq \min_\calM \calE=\int_\R|\frakv_0'|^2\dx^1
 	\quad\text{and}\quad\sigma\coloneqq\frac{c_0}{\|\fraku_0'\|_{L^2}^2},
  \label{eq:defofc0}
\end{align}
are in terms of the double well potential characterized by
\begin{align}
	c_0&=\int_{-1}^1\Big(1+\frac{1}{2}W''\Big)
	\sqrt{W+\frac{1}{4}(W')^2}\dx^1\label{eq:calculateuvnorm}\quad\text{and}\quad
	\frac{c_0}{\sigma}&=\int_{-1}^1\frac{\sqrt{W+\frac{1}{4}(W')^2}}
	{1+\frac{1}{2}W''}\dx^1.
\end{align}
\end{corollary}
\begin{proof}
We get equation \eqref{eq:calculateuvnorm} from
\begin{align*}
	\|\frakv_0'\|_{L^2(\R)}^2=\int_{-1}^1W_\ast^{\frac{1}{2}}\dx^1
	=\int_{-1}^1f'\sqrt{W+\frac{1}{4}(W')^2}\dx^1
\end{align*}
and
\begin{align*}
	\|\fraku_0'\|_{L^2(\R)}^2&=\int_{\R}\frac{|\frakv_0'|^2}{|f'\circ f^{-1}|^2}\dx^1
	=\int_{-1}^1\frac{W_\ast^{\frac{1}{2}}}{f'}\dx^1
	=\int_{-1}^1\frac{\sqrt{W+\frac{1}{4}(W')^2}}{1+\frac{1}{2}W''}\dx^1.
\end{align*}
\end{proof}

Formula \eqref{eq:calculateuvnorm} characterizes $c_0,\sigma$ in terms of the double
well potential $W$. Hereafter we define the operator
\begin{align}
	\frakA_0(w)=J\ast w\quad\text{with}\quad J(x)\coloneqq\half e^{-|x|},
	\label{eq:def1dsolutionoperator}
\end{align}
which maps $u$ to a solution $v$ of
\begin{equation*}
	-v'' + v = u\quad\text{in}\quad\R,
\end{equation*}
see \cite[Thm. 6.23]{LiebLoss2001}. This is well-defined for $w\in L^2(\R)+L^\infty(\R)$.
More details and properties are given in the Appendix \ref{lemma:expdecaysolution6}.
In Chapter \ref{sec:3} we need the following properties of the
linearization $\frakL_0$ of $\nabla_{L^2(\R)}\calE$ at $\fraku_0$.

\begin{lemma}[$\frakL_0$ is Fredholm and has a one-dimensional kernel]
 \label{lemma:1dkernel}\ \\
The operator
\begin{align}
	\frakL_0:L^2(\R)\longrightarrow L^2(\R),~\frakL_0(w)\coloneqq(f'(\fraku_0)-\frakA_0)w
	\label{2eq:def-L0}
\end{align}
is a Fredholm operator with index 0. More precisely,
\begin{align}
	\ker(\frakL_0)=\spann(\fraku_0')
	\label{eq:1dimkernel}
\end{align}
and $\frakL_0:\{\fraku_0'\}^\perp\to\{\fraku_0'\}^\perp$ is a bijection.
\end{lemma}
The proof uses a clever splitting of $\frakL_0$, which was first introduced in
\cite[Lemma 5.3]{BatesFifeRenWang1997}. The idea is to write $\frakL_0$ as the sum
of an isomorphism and a compact operator.
We adapt the method from $L^{\infty}(\R)$ to $L^2(\R)$.
\begin{proof}
We choose an arbitrary $q\in C(\R)$ such that
$\lim_{x\to\pm\infty}q(x)=f'(\pm 1)$
at an exponential rate of convergence
and $m\leq q(x)\leq M$ for some $1<m<M$ and all $x\in\R$.
We use the convolution representation of $\frakA_0$ and get
for arbitrary $w\in L^2(\R)$
\begin{align*}
	\frakL_0={}&f'(\fraku_0)\big(\frakL^{(1)}(w)+\frakL^{(2)}(w)\big)\quad\text{with}\\
	\frakL^{(1)}(w)\coloneqq{}& w-\frac{1}{q}J\ast w\quad\text{and}\quad
	\frakL^{(2)}(w)\coloneqq\Big(\frac{1}{q}-\frac{1}{f'(\fraku_0)}\Big)J\ast w.
\end{align*}
Owing to $0<\gamma\leq f'(\fraku_0)\leq C<\infty$ the operator $\frakL_0$ is Fredholm if
$\frakL^{(1)}+\frakL^{(2)}$ is. We first address $\frakL^{(1)}$. Assume
$\|w\|_{L^2(\R)}=1$, then
\begin{align*}
	\Big\|\frac{1}{q}J\ast w\Big\|_{L^2(\R)}\leq\frac{1}{m}
	\|J\|_{L^1(\R)}\|w\|_{L^2(\R)}=\frac{1}{m}<1.
\end{align*}
Therefore $\frakL^{(1)}:L^2(\R)\to L^2(\R)$ is a bijection with continuous inverse operator,
given by the corresponding Neumann series,
hence $\frakL^{(1)}$ is an isomorphism. If we can prove that $\frakL^{(2)}:L^2(\R)\to
L^2(\R)$ is a compact operator we are done. We can write $\frakL^{(2)}$ as an
integral operator and calculate the Hilbert-Schmidt norm
\begin{align*}
	\int_{\R}\bigg|\frac{1}{f'(\fraku_0(x))}-\frac{1}{q(x)}\bigg|^2
	\int_{\R}\frac{1}{4}e^{-2|x-y|}\d y\d x
	\leq\frac{1}{4m^2\gamma^2}\int_{\R}|f'(\fraku_0)-q|^2\d x<\infty
\end{align*}
owing to the exponential decay. So $\frakL^{(2)}$ is Hilbert-Schmidt and hence compact.

It follows from equations \eqref{eq:uv1dimbeziehung2} and \eqref{eq:uv1dimbeziehung}
that $\frakL_0(\fraku_0')=0$. Following the proof of Lemma 5.3 in
\cite{BatesFifeRenWang1997} we obtain $\ker(\frakL_0)=\spann(\fraku_0')$.
\end{proof}

We include a further preparation of the analysis below that deals with properties of one dimensional profiles.

\begin{lemma}\label{lem:u1}
The equation
\begin{equation}
	\frakL_0(w)-\frakA_0(\frakv_0')
	=-\tilde\sigma\fraku_0'
	\label{eq:eulerlagrangeu1-pre}
\end{equation}
has a solution $w\in L^2(\R)$ if an only if $\tilde\sigma=\sigma$ as defined in \eqref{eq:defofc0}.
For $\tilde\sigma=\sigma$ there exists a unique solution $\fraku_1$ of \eqref{eq:eulerlagrangeu1-pre} with $\fraku_1(0)=0$.
All other solutions are given by $\fraku_1+\alpha\fraku_0'$, $\alpha\in\R$.

The profile functions $\fraku_1,\frakv_1:=\frakA_0 \big(\fraku_1+\frakv_0'\big)$ and $\frakv_1'$ all decay at an exponential rate to $0$ at $\pm\infty$.
\end{lemma}

\begin{proof}
  By Lemma \ref{lemma:1dkernel} the range of $\frakL_0$ is $\{\fraku_0'\}^\perp$.
  This yields as a necessary and sufficient condition for solvability
	\begin{align*}
		0=\langle\fraku_0'|-\tilde\sigma\fraku_0'+\calA_0(\frakv_0')\rangle
		\Longleftrightarrow\tilde\sigma=\frac{\|\frakv_0'\|^2}{\|\fraku_0'\|^2}=\sigma.
	\end{align*}
	In this case, since $\frakL_0:\{\fraku_0'\}^\perp\longrightarrow\{\fraku_0'\}^\perp$
	is a bijection, there is a unique solution $\fraku_1$ in $\{\fraku_0'\}^\perp$ and the solution space in $L^2(\R)$ is given by $\fraku_1+\ker(\frakL_0)$.
  Since $\fraku_0'(0)>0$ we find a unique
	solution $\fraku_1\in L^2(\R)$ with $\fraku_1(0)=0$. For the exponential decay
	we need a standard ODE argument, this will be explained in the appendix in Lemma
	\ref{lemma:expconvergenceforall}
\end{proof}

\section{The $\Gamma$--limit candidate and the $\Gamma$--limsup construction}
\label{sec:3}
In this section we identify the candidate for the $\Gamma$--limit, prove a $\liminf$--inequality for sequences that are representable by certain asymptotic expansion (see Assumption \ref{ass:expansion} below), and show the $\limsup$--estimate needed for the $\Gamma$--convergence.

Achieving a general $\liminf$--estimate seems currently out of reach.
One key obstacle is the lack of a suitable concept of approximating varifolds that allows for an associated first variation formula, which was crucial in the analysis of $\calW_\eps^{\text{CH}}$ from \cite{RoegerSchaetzle2006}.

\subsection{Identification of the $\Gamma$--limit}\label{sec:3ident}
In this section we consider the approximate Willmore energy $\calW_\eps$ as defined in equation \eqref{eq:defwillmoreapproxag}.
For a given phase boundary $\Gamma$ we identify energy-optimal approximations in the class of phase fields representable in terms of an asymptotic expansion, see below for the precise definitions. We therefore identify a candidate for the $\Gamma$--limit of $(\calW_\eps)_{\eps>0}$ and also prepare the proof of the corresponding $\Gamma$--limsup estimate.

In the following we denote by $\calA_\eps:L^2(\Omega)\longrightarrow L^2(\Omega)$ the solution operator $u\longmapsto v_\eps$ to the problem \eqref{eq:v_eps}.
With this notation we can reformulate the {\em diffuse mean curvature} as
\begin{align*}
	H_\eps(u)=\frac{1}{\eps}(f(u)-\calA_\eps( u)).
\end{align*}

\begin{assumptions}\label{ass:Gamma}
We consider a fixed open subset $E\ssubset\Omega$ with $C^4$-boundary
$\Gamma\coloneqq\partial E$, inner normal $\nu$ and principal curvatures
$k_1,\dots,k_{n-1}$ in direction of $\nu$. We denote the signed distance
function to $\Gamma$ by $d\coloneqq\dist(\cdot,\Omega\setminus E)-\dist(\cdot,E)$.
Let $\delta>0$ be chosen such that $d$ is $C^4$-regular and the projection
$y(x)\coloneqq\Pi_\Gamma(x)$ on $\Gamma$ is well-defined in $\omega\coloneqq\{|d|<5\delta\}$.
\end{assumptions}

\begin{remark}[Choice of coordinates and notations]
\label{defi:newcoordinates}
We define a parametrization
$\Psi_\eps:(-\frac{5\delta}{\eps},\frac{5\delta}{\eps})\times\Gamma\longrightarrow \omega$ of $\omega$ by
\begin{align}
	\Psi_\eps(z,y)=y+\eps z\nu(y).
	\label{eq:newcoordinatesnotime}
\end{align}
By Assumption \ref{ass:Gamma} the map $\Psi_\eps$ is bijective and $C^4$-regular, with inverse characterized by
\begin{align}
  x=\Psi_\eps\Big(\frac{d(x)}{\eps},y(x)\Big)
  \quad\text{for any}\quad x\in \omega.
  \label{eq:newcoordinatesnotime-x}
\end{align}
Possibly lowering the value of $\delta>0$ we can assume that in $(-\frac{5\delta}{\eps},\frac{5\delta}{\eps})\times\Gamma$
\begin{equation}
  2\eps \geq \det(D\Psi_\eps(z,y)) = \eps -\eps^2zH(y) +\eps^3z^2 R_\eps(z,y)\geq \frac{\eps}{2},
  \label{eq:determinant}
\end{equation}
where $R_\eps:\R\times\Gamma\to\R$ is uniformly bounded.

We can represent a function $u:\omega
\longrightarrow\R$ as $\tilde u:(-5\delta/\eps,5\delta/\eps)
\times\Gamma\longrightarrow\R,$ $\tilde u(z,y)\coloneqq u(\Psi_\eps(z,y))$.
Often it is convenient to extend $\tilde u$ to a function $U$ on $(-5\delta/\eps,5\delta/\eps)\times\omega$ that is constant in normal directions, i.e.~$U(z,x)=\tilde u(z,y)$ for all $z\in (-5\delta/\eps,5\delta/\eps)$, $x\in\omega$ with $\Pi_\Gamma(x)=y$.

We then write $U'=\partial_z U$ for the $z$-derivative, and $\nabla U$ and $\Delta U$ for the gradient and Laplace-Operator with respect to the $x$-variable.
From \cite{LoretiMarch2000, Wang2008} we recall that
\begin{align}
	\nabla u&=\frac{1}{\eps} U'\nabla d+\nabla U,\quad
  \Delta u=\frac{1}{\eps^2}U''
  +\frac{\Delta d}{\eps}U'+\Delta U,
  \label{eq:Delta_in_zy}\\
 	|\nabla d| &=1,\quad\nabla d\cdot\nabla_{\Gamma}=0,\quad
  \Delta d(x)=H(y)-\eps z|\FF|^2(y)+\eps^2 |z|^2R_\eps^H(x),
 	\label{eq:exp_Delta_d}
\end{align}
where $|R_\eps^H|\leq C(\Gamma,\delta)$ in $\omega$.

In the following we will often write $u(x,z)$ instead of $U(x,z)$ if it is clear from the context what is meant.
\end{remark}

We next introduce a class of phase field approximations of $\Gamma$ that can be represented
by an inner asymptotic expansion and an outer expansion given by the rescaled characteristic
function $2\Chi_E-1$. Since the nonlocal operator $\calA_\eps$ acts on functions defined on
the whole domain we cannot only consider the inner region.
In the following we
prepare the definition of an appropriate class of phase field functions.
We will fix two different
modifications of the signed distance function and one cut-off function.

\begin{assumptions}[Modified distance and cut-off functions]
\label{ass:modified_distances}
Choose an odd and increasing function $\phi_1\in C^\infty(\R)$ with
\begin{align*}
	\phi_1'(0)=0,\quad 0 &< \phi_1(z)\leq \frac{9}{10}z,\quad
  0\leq\phi_1'(z)\leq 1 \qquad\text{ for all }z\in (0,\infty),\\
	\phi_1(z)&=
	\begin{cases}
		z-\frac{1}{4},\quad &\text{ if } z\in(\frac{1}{2},2)\\
		2,\quad &\text{ if } z\in(\frac{5}{2},\infty).
	\end{cases}
\end{align*}
We set $\phi_\delta(z):=\delta\phi_1\big(\frac{z}{\delta}\big)$ for $z\in\R$ and define a
modification of the signed distance function (being constant outside $\{|d|<4\delta\}$) by
\begin{equation*}
  d_\delta:=\phi_\delta\circ d\in C^4(\Omega).
\end{equation*}
Finally, choose an even and on $(0,\infty)$ decreasing function
$\eta_1\in C^\infty_c(\R)$ with
\begin{equation*}
	0\leq \eta_1\leq 1,\quad |\eta_1'|\leq 2,
	\quad
	\eta_1 =
	\begin{cases}
		1 \quad&\text{ in }[0,3],\\
		0 \quad&\text{ in }[4,\infty)
	\end{cases}
\end{equation*}
and define the cut-off function
\begin{equation*}
  \eta_\delta(x):=\eta_1\Big(\frac{d(x)}{\delta}\Big)\quad\text{ for all }x\in\Omega.
\end{equation*}
We remark that $\eta_\delta\in C^4(\Omega)$ since $\eta_\delta$ has
support in $\{|d|<4\delta\}$.
\end{assumptions}

We next define spaces of functions that are exponentially small away from $\Gamma$.
\begin{definition}
For $\Lambda,\mu>0$ we consider
\begin{align*}
	X_\delta^{\mu,\Lambda}(\Omega):=\Big\{w\in L^\infty(\Omega)\,\big|\,
	\esssup_{x\in\Omega}|e^{\mu |d_\delta(x)|}w(x)|\leq\Lambda\Big\}
\end{align*}
and
\begin{align*}
	X^{\mu,\Lambda}(\R;\Gamma) := \Big\{w\in L^\infty(\R\times\omega)\,\big|\,
	&\esssup_{(z,x)\in\R\times\omega}|e^{\mu |z|}w(z,x)|\leq\Lambda,\\
  &w(z,\cdot)\text{ is constant in normal direction}\}
\end{align*}
and set $X(\R;\Gamma):=\bigcup_{\mu,\Lambda>0}X^{\mu,\Lambda}(\R;\Gamma)$.
\end{definition}

Note that for $0<\mu_1<\mu_2$ and any $\Lambda>0$ we have
$X_\delta^{\mu_2,\Lambda}(\Omega)\subseteq X_\delta^{\mu_1,\Lambda}(\Omega)$.
Next we define a suitable class of phase field approximations that have an inner
and outer expansion (the latter given by $\pm 1$) up to order $\eps^K$.
\begin{assumptions}
\label{ass:expansion}
Let $\Gamma$, $\delta\in(0,1)$ as in Assumption \ref{ass:Gamma} and $K\in\N_0$ be given.
Consider a family $(u_\eps)_{0<\eps<\eps_0}$ that can be represented as follows:
There exist $\mu\in (0,1)$, $\Lambda>0$,
and profile functions $u_j$, $j=0,\dots,K$, such that for all $0<\eps<\eps_0$
\begin{align}
	u_\eps &= \eta_\delta \uei +
	(1-\eta_\delta)\sgn(d)+\eps^{K+1}R_\eps^u\quad\text{ in }\Omega,\\
  \uei &= \Big(\sum_{j=0}^K \eps^ju_j\Big)\circ \Psi_\eps^{-1}
  \quad\text{ in }\{|d|<3\delta\},
\end{align}
and such that the following properties hold:
\begin{enumerate}
	\item The profile functions
	$u_j\in C^0(\R\times\omega)$,	$u_j=u_j(z,x)$ are $C^4$-regular with
	respect to the $x$-variable and satisfy
	\begin{align*}
		u_0-\sgn &\in X(\R;\Gamma),\ u_j\in X(\R;\Gamma)\text{ for }1\leq j\leq K,\\
    |\nabla u_j|,\,\Delta u_j &\in X(\R;\Gamma) \text{ for }0\leq j\leq K.
	\end{align*}
	\item The remainder satisfies
  $R_\eps^u\in X^{\frac{\mu}{\eps},\Lambda}_\delta(\Omega)$ for all $0<\eps<\eps_0$.
\end{enumerate}
Finally we assume that there are height functions $h_j$, $j=0,\dots,K-1$ such that
\begin{equation}
   y\mapsto y+\eps\Big(\sum_{j=0}^{K-1}\eps^jh_j(y)+\eps^KR^h_\eps(y)\Big)\nu(y)\,,\,y\in\Gamma
   \label{eq:ueps=0}
\end{equation}
is a $C^4$-diffeomorphism onto $\{u_\eps=0\}$ with  $\sup_{\eps>0}\|R^h_\eps\|_{C^4(\Gamma)}<\infty$.
\end{assumptions}

\begin{lemma}[Convergence towards $u$]
\label{lemma:u_epstou}
Let $E,\Gamma$ be as in Assumption \ref{ass:Gamma} and set $u:=2\Chi_E-1$.
Consider $K\in\N_0$ and $(u_\eps)_{0<\eps<\eps_0}$
as in Assumption \ref{ass:expansion}.
Then we have $u_\eps\to u$ in $L^p(\Omega)$ for all $1\leq p<\infty$.
\end{lemma}
\begin{proof}
Note that $2\Chi_E-1=\sgn(d)$. We fix $\Lambda,\mu,\delta>0$ such that
$R^u_\eps\in X^{\frac{\mu}{\eps},\Lambda}_\delta(\Omega)$ and $u_0-\sgn,\, u_j\in X^{\mu,\Lambda}(\R;\Gamma)$ for all $1\leq j\leq K$.
Using the representation of $u_\eps$ and \eqref{eq:determinant} we have
\begin{align*}
	\int_\Omega|u_\eps-u|\dx^n&=
	\int_\Omega\big|\eta_\delta (\uei-\sgn(d))+\eps^{K+1}R_\eps^u\big|\dx^n\\
	&\leq\int_{\{|d|<4\delta\}}\big|\uei-\sgn(d)|\dx^n
	+\eps^{K+1}\int_\Omega\big|R_\eps^u\big|\dx^n\\
	&\leq\int_\Gamma\int_{-\frac{2\delta}{\eps}}^{\frac{4\delta}{\eps}}
	2\eps\Big|\sum_{j=0}^K \eps^ju_j-\sgn\Big|\dx^1\ds^{n-1}+\eps^{K+1}\Lambda
	|\Omega|.
\end{align*}
We use the bounds for the profile functions to further estimate the right-hand side and deduce that
\begin{align*}
	\int_\Omega|u_\eps-u|\dx^n&\leq 2\calH^{n-1}(\Gamma)\eps
	\int_\R e^{-\mu|t|}\Lambda\dx^1(t)\\
	&+2\calH^{n-1}(\Gamma)\sum_{j=1}^K\eps^{j+1}\int_\R e^{-\mu|t|}\Lambda\dx^1(t)+\eps^{K+1}\Lambda
	|\Omega|\\
	&\leq\frac{4\Lambda\calH^{n-1}(\Gamma)}{\mu}\sum_{j=0}^K\eps^{j+1}+\eps^{K+1}\Lambda
	|\Omega|.
\end{align*}
Since $(u_\eps)_{\eps>0}$ is uniformly bounded in $L^\infty(\Omega)$ the convergence follows for all $1\leq p<\infty$.
\end{proof}

Below we will only need orders $K\leq 2$. The key observation at this point is that the solution operator $\calA_\eps$ conserves the expansion properties.

\begin{proposition}\label{pro:expansion-v}
Consider $K\leq 2$ and $(u_\eps)_{0<\eps<\eps_0}$ as in Assumption \ref{ass:expansion}.
Then the family $(v_\eps)_{0<\eps<\eps_0}$, $v_\eps=\calA_\eps u_\eps$
has an analogue representation, meaning that there exist $\tilde\Lambda>0$ and profile functions
$v_j\in C^2(\R\times\omega)$, $j=0,\dots,K$ such that for all $0<\eps<\eps_0$
\begin{align}
	v_\eps &= \eta_\delta \vei +
	(1-\eta_\delta)\sgn(d)+\eps^{K+1} R_\eps^v\quad\text{ in }\Omega,
  \label{eq:expansion-v-1}\\
  \vei &= \bigg(\sum_{j=0}^K\eps^jv_j\bigg)\circ \Psi_\eps^{-1}
  \quad\text{ in }\{|d|<4\delta\},
  \label{eq:expansion-v-2}
\end{align}
with $R_\eps^v\in X^{\frac{\mu}{\eps},\tilde\Lambda}_\delta(\Omega)$
for all $0<\eps<\eps_0$.

The profile functions are given by
\begin{align}
  v_0(z,x)&=\frakA_0 u_0(z,x) \label{eq:veps-0}\\
  v_1(z,x)&=\frakA_0 \big(u_1(z,x)+H(y)v_0'(z,x)\big) \label{eq:veps-1}\\
  v_2(z,x)&=\frakA_0 \big(u_2(z,x)+H(y)v_1'(z,x)+(\Delta-z|\FF|^2(y)\partial_z)v_0(z,x)\big) \label{eq:veps-2}
\end{align}
for $z\in\R$, $x\in \omega$, $y=\Pi_\Gamma(x)$.
Moreover,
\begin{align*}
	v_0-\sgn,\, v_1, v_2 \in X(\R;\Gamma)\quad\text{ and }\quad
	\partial_z v_j,\,|\nabla v_j|,\,\Delta v_j \in X(\R;\Gamma) \,\text{ for }j=0,1,2.
\end{align*}
\end{proposition}
Following the proof of Lemma \ref{lemma:u_epstou} we deduce $v_\eps\to 2\Chi_E-1$ in $L^1(\Omega)$.
\begin{proof}

The full proof is given in Appendix \ref{sec:A-expansion}. Here we only formally derive the
characterization \eqref{eq:veps-0}-\eqref{eq:veps-2} of the profile functions.

Assume $K=2$. Using the representation for $\Delta d$ from \eqref{eq:exp_Delta_d} we obtain in $\{|d|<3\delta\}$
\begin{align*}
	-\eps^2\Delta+\Id&=-\partial_z^2+\Id-\eps \Delta d\partial_z-\eps^2\Delta\\
	&=-\partial_z^2+\Id-\eps H\partial_z-\eps^2(\Delta-z|\FF|^2\partial_z)-\eps^3|z|^2R_\eps^H\partial_z,
\end{align*}
and
\begin{align}
	&\Big(\uei-(-\eps^2\Delta+\Id)\vei \Big)\circ\Psi_\eps^{-1}\nonumber\\
  &\qquad = u_0+\eps u_1+\eps^2 u_2-(-\partial_z^2+\Id)v_0-\eps\big((-\partial_z^2+\Id)v_1-H\partial_z v_0\big)\nonumber\\
  &\qquad\quad -\eps^2\big((-\partial_z^2+\Id)v_2-H\partial_z v_1-(\Delta
	-z|\FF|^2\partial_z)v_0\big)+\eps^3R_\eps^v,
  \label{eq:uei-vei}
\end{align}
where $R_\eps^v$ is uniformly bounded in $\omega$ (see also the representation \eqref{eq:pf3.7-1} in the Appendix).
The equations \eqref{eq:veps-0}-\eqref{eq:veps-2} are then equivalent to the property, that  the expression in \eqref{eq:uei-vei}
vanishes up to order $\bigO(\eps^2)$.
\end{proof}

Also the application of $f$ retains the properties of an asymptotic expansion.
\begin{lemma}\label{lem:expansion-f}
Consider $K\leq 2$ and $(u_\eps)_{0<\eps<\eps_0}$ as in Assumption \ref{ass:expansion}.
Then the family $\big(f(u_\eps)\big)_{0<\eps<\eps_0}$ can be represented as
\begin{align}
  f(u_\eps) = \eta_\delta f(u_\eps)^{\inn} +(1-\eta_\delta)\sgn d + \eps^{K+1}R_\eps
  \label{eq:expansion-f}
\end{align}
with $f(u_\eps)^{\inn}\circ \Psi_\eps = \sum_{j=0}^K \eps^j F_j$,
\begin{equation*}
  F_0 =f(u_0),\quad F_1=f'(u_0)u_1,\quad F_2=\frac{1}{2}f''(u_0)(u_1)^2+f'(u_0)u_2.
\end{equation*}
Moreover, $F_0-\sgn,F_1,F_2 \in X(\R;\Gamma)$, $R_\eps\in X^{\frac{\mu}{\eps},C(f)\Lambda}_\delta(\Omega)$ holds.
\end{lemma}

\begin{proof}
We proof the claim for $K=2$.
Choose $\mu,\Lambda$ such that $u_0,u_1,u_2\in X^{\mu,\Lambda}(\R)$, $R_\eps\in X_\delta^{\frac{\mu}{\eps},\Lambda}$.
We first obtain
\begin{equation}
  f(u_\eps) = f\big(\eta_\delta u_\eps^{\inn} +(1-\eta_\delta)\sgn d + \eps^3R_\eps^u\big) = f\big(\eta_\delta u_\eps^{\inn} +(1-\eta_\delta)\sgn d \big)+ \eps^3R_\eps^{(1)},
  \label{eq:pf3.8-1}
\end{equation}
with $|R_\eps^{(1)}|\leq C(f)|R_\eps^u|$.
Since $f\in C^2([-1,1])$ a Taylor expansion (see Lemma \ref{lem:calc}
in the Appendix) yields
\begin{equation}
  \big|f\big(\eta_\delta \uei + (1-\eta_\delta)\sgn d\big)
  - \eta_\delta f(\uei) - (1-\eta_\delta)\sgn d\big|
  \leq C(f)\eta_\delta (1-\eta_\delta)(\uei -\sgn d)^2.
  \label{eq:pf3.8-2}
\end{equation}
Another Taylor expansion implies that in $\{\eta_\delta>0\}$
\begin{equation}
  \Big|f(\uei) - \Big(f(u_0)+ \eps f'(u_0)u_1+ \eps^2
  \Big(\frac{1}{2}f''(u_0)(u_1)^2+f'(u_0)u_2\Big)\Big)\circ\Psi_\eps^{-1}\Big|
  \leq \eps^3 C(f)R_\eps^{(2)}
  \label{eq:pf3.8-3}
\end{equation}
for $R_\eps^{(2)}$ with $R_\eps^{(2)}\leq (|u|_1+|u_2|+|u_3|)^2\circ\Psi_\eps^{-1}$.
From \eqref{eq:pf3.8-1}-\eqref{eq:pf3.8-3} we conclude the desired representation \eqref{eq:expansion-f} with
\begin{equation*}
  |R_\eps|\leq C(f)\Big[R_\eps^{(1)}+\eta_\delta (1-\eta_\delta)(\uei -\sgn d)^2+(|u|_1+|u_2|+|u_3|)^2\circ\Psi_\eps\Big].
\end{equation*}
Since $|\uei -\sgn d|\leq \big(|u_0-\sgn|+|u_1|+|u_2|\big)\circ\Psi_\eps$ we deduce that $R_\eps\in X_\delta^{\frac{\mu}{\eps},\Lambda}$.
Finally, $F_0-\sgn,F_1,F_2 \in X(\R;\Gamma)$ follows from $|f(u_0)-\sgn|=|f(u_0)-f(\sgn)|\leq C(f)|u_0-\sgn|$ and the assumptions on $u_0,u_1,u_2$.
\end{proof}

The next theorem identifies a candidate for the $\Gamma$--limit of $(\calW_\eps)_{\eps}$
in $\Gamma$ and proves a corresponding lower bound estimate in the just introduced class of Ansatz functions.

\begin{theorem}\label{thm:liminf-ansatz}
Let $E,\Gamma$ satisfy Assumption \ref{ass:Gamma} and denote by $u:=2\Chi_E-1$
the rescaled characteristic function of $E$.
For any sequence $(u_\eps)_{\eps}$ as in Assumption \ref{ass:expansion} there holds
\begin{equation}
  c_0\sigma\calW(u)\leq\liminf_{\eps\to 0}\calW_\eps(u_\eps).
  \label{eq:liminf-expansion}
\end{equation}
\end{theorem}

\begin{proof}
Let $v_\eps=\calA_\eps(u_\eps)$.
By Assumption \ref{ass:expansion}, Proposition \ref{pro:expansion-v} and
Lemma \ref{lem:expansion-f} we deduce
\begin{align}
  \eps H_\eps &= \big(f(u_\eps)-v_\eps\big)
  \nonumber\\
  &= \eta_\delta f(u_\eps)^{\inn} +(1-\eta_\delta)\sgn d
  -\big(\eta_\delta \vei + (1-\eta_\delta)\sgn(d)\big)+\eps^{K+1} R_\eps\nonumber\\
  &= \eta_\delta \big(f(u_\eps)^{\inn}- \vei\big)+\eps^{K+1} R_\eps,
  \label{eq:Heps-K}
\end{align}
and in particular for $K=0$
\begin{align*}
  \eps H_\eps &= \eta_\delta\cdot\big( f(u_0)-\frakA_0 u_0\big)
  \circ\Psi_\eps^{-1}+\eps R_\eps.
\end{align*}

Together with \eqref{eq:determinant} we deduce
\begin{align}
  \eps^2 \calW_\eps(u_\eps) &=
  \int_\Omega \frac{1}{\eps}\eta_\delta^2\big(f(u_0)
  -\frakA_0 u_0\big)^2\circ\Psi_\eps^{-1}\dx^n\nonumber\\
  &\qquad + 2\int_\Omega \eta_\delta \big(f(u_0)-\frakA_0 u_0\big)\circ\Psi_\eps^{-1}
  R_\eps\dx^n +\eps \int_\Omega R_\eps^2\dx^n \nonumber\\
  &\geq \int_\Gamma \int_{-\frac{3\delta}{\eps}}^{\frac{3\delta}{\eps}}
  \frac{1}{2}\big(f(u_0)-\frakA_0 u_0\big)^2\dx^1\ds^{n-1}
  -C\eps\int_\Gamma \int_{-\infty}^{\infty}
  \big|f(u_0)-\frakA_0 u_0\big|\dx^1\ds^{n-1} \nonumber\\
  &\geq \int_\Gamma \int_{-\frac{3\delta}{\eps}}^{\frac{3\delta}{\eps}}
  \frac{1}{2}\big(f(u_0)-\frakA_0 u_0\big)^2\dx^1\ds^{n-1}-C\eps,
  \label{eq:limsup-1}
\end{align}
where we have used, that $f(u_0)-\sgn$ and $\frakA_0u_0-\sgn$ both decay exponentially at $\pm \infty$.

In order to prove \eqref{eq:liminf-expansion} it is sufficient to consider the case $\liminf_{\eps\to 0}\eps^2\calW_\eps(u_\eps)=0$.
This implies $f(u_0)=\frakA_0 u_0$, that is $u_0(\cdot,x)=\fraku_0(\cdot-z_0(y))$ and $v_0(\cdot,x)=\frakv_0(\cdot-z_0(y))$ with $y=\Pi_\Gamma(x)$.
The condition \eqref{eq:ueps=0} implies that $z_0(y)=h_0(y)$.

\medskip
With $u_0(\cdot,x)=\fraku_0(\cdot-h_0(y))$ and $v_0(\cdot,x)=\frakv_0(\cdot-h_0(y))$ we deduce from Proposition \ref{pro:expansion-v}, Lemma \ref{lem:expansion-f} and \eqref{eq:Heps-K} with $K=1$
\begin{align}
  H_\eps(x) &= \eta_\delta(x) \Big(f'\big(\fraku_0(z-h_0(y))\big)u_1(z,x)-\frakA_0 \big(u_1(z,x)+H(y)\frakv_0'(z-h_0(y))\big)\Big)+\eps R_\eps(x).
    \label{eq:Heps-K=1}
\end{align}
and by similar calculations as above
\begin{align}
  &\liminf_{\eps\to 0}\calW_\eps(u_\eps)\nonumber\\
  &\geq \int_\Gamma \int_{-\frac{3\delta}{\eps}-h_0(y)}^{\frac{3\delta}{\eps}-h_0(y)}
  \big|\big(f'(\fraku_0)-\frakA_0\big)
	u_1(\cdot+h_0(y),x)-H(y)\frakA_0(\frakv_0')\big|^2 \dx^1\ds^{n-1}(y)\\
  &\geq\int_\Gamma \int_{-\infty}^{\infty}
  \big|\big(f'(\fraku_0)-\frakA_0\big)
	u_1(\cdot+h_0(y),x)-H(y)\frakA_0(\frakv_0')\big|^2 \dx^1\ds^{n-1}(y).
  \label{eq:limsup-2}
\end{align}

If $H(y)=0$ the inner integral is minimized by $u_1(\cdot+h_0(y),x)\equiv 0$. Therefore,
to prove a lower bound, we can assume $u_1(z+h_0(y),x)=H(y)\tilde u_1(z,x)$ and compute
\begin{align}
  \liminf_{\eps\to 0}\calW_\eps(u_\eps)
  &\geq \int_\Gamma|H|^2\ds^{n-1}\inf_{w\in L^2(\R)}\Xi(w),
	\label{eq:formalerechnungzwischenschrittinersterordnung}
\end{align}
with the functional $\Xi:L^2(\R)\to [0,\infty]$
\begin{align*}
	\Xi(w)\coloneqq\int_{-\infty}^{\infty}\big|\frakL_0(w)-\frakA_0(\frakv_0')\big|^2\dx^1.
\end{align*}
From Lemma \ref{lemma:1dkernel} we get that any minimizer $w$ of
$\Xi$ is characterized by the fact that $\frakL_0(w)$ is equal to the
$L^2$-projection of $\frakA_0(\frakv_0')$ onto
$\range(\frakL_0)=\{\fraku_0'\}^\perp$, in particular
\begin{align*}
	\langle \frakL_0(\phi)|\frakL_0(w)-\frakA_0(\frakv_0')
	\rangle_{L^2}=0\quad\text{ for all }\phi\in L^2(\R),
\end{align*}
which implies
\begin{align}
	\frakL_0\big(\frakL_0(w)-\frakA_0(\frakv_0')\big)=0.
	\label{eq:u1}
\end{align}
Using Lemma \ref{lemma:1dkernel} we conclude
\begin{equation}
	\frakL_0(w)-\frakA_0(\frakv_0')
	=-\tilde \sigma\fraku_0'.
	\label{eq:eulerlagrangeu1}
\end{equation}
By Lemma \ref{lem:u1} we obtain that $\tilde\sigma=\sigma$ and that
$w=\fraku_1+\alpha \fraku_0'$ for some $\alpha\in\R$.
This implies
\begin{align}
	\inf_{w\in L^2(\R)}\Xi(w)=\Xi(\fraku_1+\alpha \fraku_0')=\Xi(\fraku_1)=c_0\sigma.
	\label{eq:ergebnisminimierenersteordnung}
\end{align}
Together with \eqref{eq:formalerechnungzwischenschrittinersterordnung} this proves \eqref{eq:liminf-expansion}.

Finally, we can determine $\alpha$ from condition \eqref{eq:ueps=0}, which implies
\begin{align*}
  0 &= u_\eps\big(y+\eps(h_0(y)+\eps h_1(y)+\eps^2R_\eps^h(y))\nu(y)\big)\\
  &= \fraku_0\big(\eps h_1(y)+\eps^2R_\eps^h(y)\big)
  + \eps\big(\fraku_1\big(\eps h_1(y)+\eps^2R_\eps^h(y)\big)+\alpha\fraku_0'\big(\eps h_1(y)+\eps^2R_\eps^h(y)\big)\big)+\bigO(\eps^2)\\
  &= \eps \fraku_0'(0)\big(h_1(y)+\alpha\big)+\bigO(\eps^2)
\end{align*}
and therefore $\alpha=-h_1(y)$.
\end{proof}

The proof shows that equality in \eqref{eq:liminf-expansion} can only be attained if $u_0(z,x)=\fraku_0(\cdot-h_0(y))$, $u_1(z,x)=H(y)\fraku_1(z)-h_1(y)\fraku_0'(z)$.
By Theorem \ref{thm:1dminimizer} we have that $\fraku_0-\sgn,\frakv_0-\sgn$, $\fraku_0'$ all decay exponentially at $\pm\infty$.
Lemma \ref{lem:expansion-f}, Lemma \ref{lemma:expconvergenceforall} and Lemma \ref{lemma:expdecaysolution6} in the Appendix show that $\fraku_1$, $\frakv_0'$ also decay exponentially at $\pm\infty$.

We therefore obtain as a candidate for a recovery sequence $(u_\eps)_{\eps>0}$
\begin{equation}
  u_\eps = \eta_\delta \uei + (1-\eta_\delta)\sgn(d),\quad
  \uei(\cdot,x)= \fraku_0+\eps H(y)\fraku_1.
  \label{eq:ansatz-limsup}
\end{equation}

\subsection{$\Gamma$--limsup construction}
This section is devoted to the constructive part necessary to prove the
$\Gamma$--convergence of $\calW_\eps$ to $c_0\sigma\calW$.
We use the previous computations and the candidate \eqref{eq:ansatz-limsup}.

\begin{theorem}[$\limsup$ estimate for Willmore approximation]
\label{thm:limsupamstutzvangoethem}\ \\
Let $E,\Gamma$ satisfy Assumption \ref{ass:Gamma} and denote by $u:=2\Chi_E-1$ the rescaled characteristic function of $E$. Then there exists a sequence $(u_\eps)_{\eps>0}$
such that $u_\eps\to u$ in $L^1(\Omega)$ and
\begin{align*}
	\lim\limits_{\eps\to 0}\calW_\eps(u_\eps)=c_0\sigma\calW(u).
\end{align*}
\end{theorem}

\begin{proof}
The convergence towards $u$ was already shown in Lemma \ref{lemma:u_epstou}.
We will use the Ansatz \eqref{eq:ansatz-limsup} and let
\begin{equation*}
  u_\eps(x) = \eta_\delta(x) \big(\fraku_0(z)+\eps H(y)\fraku_1(z)\big) + \big(1-\eta_\delta(x)\big)\sgn(d(x)),
\end{equation*}
where $x=\Psi_\eps(z,y)$.
We deduce from Proposition \ref{pro:expansion-v} that $v_\eps := \calA_\eps u_\eps$ can be represented as
\begin{equation*}
  v_\eps(x) = \eta_\delta(x)\big(\frakv_0(z) + \eps \frakv_1(z,x)\big) + (1-\eta_\delta(x))\sgn z + \eps^2 R_\eps(x)
\end{equation*}
with $\frakv_0 =\frakA_0 \fraku_0$ as in \eqref{eq:veps-0} and
\begin{align}
  v_1(\cdot,x)&= H(y)\frakv_1\,\text{ for any }x\in\omega,\,y=\Pi_\Gamma(x),
  \quad \frakv_1=\frakA_0 \big(\fraku_1+\frakv_0'\big),
  \label{3eq:def-v1}
\end{align}
as introduced in Lemma \ref{lem:u1}.
Moreover, we have $\sup_{\eps>0}\sup_{x\in\Omega}|R_\eps(x)|\leq C$ and $\frakv_1,\frakv_1'$
decay exponentially at $\pm\infty$ by Lemma \ref{lemma:expdecaysolution6}. We deduce from
equations \eqref{eq:Heps-K=1} and \eqref{eq:eulerlagrangeu1}.
\begin{align*}
  H_\eps(x) &= \eta_\delta(x) H(y) \big(f'(\fraku_0)\fraku_1-\frakA_0
  (\fraku_1+\frakv_0')\big)(z)+\eps R_\eps(x)
  \nonumber\\
  &= -\eta_\delta(x) \sigma H(y) \fraku_0'(z)+\eps R_\eps(x).
\end{align*}
By similar calculations as above this implies
\begin{align*}
  \calW_\eps(u_\eps)&=
  \int_\Omega \frac{1}{\eps}\eta_\delta^2 (\sigma \fraku_0' )^2H^2\circ\Psi_\eps^{-1}\dx^n
  + 2\int_\Omega \eta_\delta (\sigma \fraku_0' )H\circ\Psi_\eps^{-1} R_\eps \dx^n
  +\eps \int_\Omega R_\eps^2\dx^n \nonumber\\
  &\leq \sigma^2 \int_\Gamma H^2\ds^{n-1}\,\int_{-\infty}^{\infty}
  (\fraku_0')^2\dx^1\nonumber
  +C\calH^{n-1}(\Gamma)\eps \int_{-\infty}^{\infty}
  \big|\fraku_0'\big|(z)\dx^1+\eps C|\Omega|\\
  &\leq c_0\sigma \int_\Gamma H^2\ds^{n-1} +C\eps.
\end{align*}
This yields $\limsup_{\eps\downarrow 0}\calW_\eps(u_\eps)\leq c_0\sigma \calW(u)$ and
together with \eqref{eq:liminf-expansion} the recovery sequence property.
\end{proof}

\section{Convergence towards the Mean Curvature Flow and Willmore flow: asymptotic expansion}
In this section we consider the rescaled $L^2$-gradient
flows \eqref{eq:appMCF}, \eqref{eq:appWillflow} of the gradient-free diffuse
Willmore and Perimeter energy.

We recall the the first is given by
\begin{equation}
	\eps\partial_tu_\eps = -H_\eps\label{eq:timeevperimeter}
\end{equation}
and the second by
\begin{equation}
	\eps\partial_tu_\eps = -\frac{2}{\eps^2}L_\eps (H_\eps ),
	\label{eq:timeev}
\end{equation}
where $H_\eps$ denotes the diffuse mean curvature defined in equation \eqref{eq:defH_eps} and
\begin{align*}
	L_\eps \coloneqq f'(u_\eps  )\Id-\calA_\eps
\end{align*}
with $\calA_\eps$ as in equation \eqref{eq:v_eps}.

We refer to \eqref{eq:mcf} and \eqref{eq:willmoreflow} for the formulation of mean curvature and Willmore flow.

\medskip
Our goal is to justify by an asymptotic expansion that the evolution
\eqref{eq:timeevperimeter} approximates the
MCF and evolution \eqref{eq:timeev} approximates the Willmore flow.

The justification of phase field approximations of geometric evolution laws has a long history.
Our analysis closely follows the formal analysis in \cite{LoretiMarch2000}, see also
\cite{Wang2008,BretinMasnouOudet2015} and \cite{RaetzRoeger2021}.
We formulate here assumptions under which the derivation is rigorous.
This however does not give a general convergence proof, since the assumed properties need to be verified for a phase field approximation.
Complete convergence proofs based on asymptotic expansion techniques are known for the standard diffuse approximation of mean curvature and Willmore flow, see \cite{MoSc90} and \cite{FeiLiu2021}.

\begin{assumptions}[Set evolution]\label{ass:evolution}
Consider an evolution of open sets $(E(t))_{t\in [0,T]}$ in $\Omega$ with associated phase boundaries $\Gamma(t):=\partial E(t)$, $t\in [0,T]$ and signed distance function $d:\Omega\times [0,T]\to \R$, $d(\cdot,t)=
\dist(\cdot,\Omega\setminus E(t))-\dist(\cdot,E(t))$.

We choose $\delta>0$ sufficiently small such that for all $t\in [0,T]$ the projections $\Pi_{\Gamma(t)}:\{|d(\cdot,t)|<5\delta\}\to \Gamma(t)$ are well defined and set
\begin{equation*}
  \omega_T := \{(x,t)\in\Omega\times [0,T]\,:\, |d(x,t)|<5\delta\}.
\end{equation*}

We assume the following properties:
\begin{enumerate}
  \item $\Gamma(t)$ is a $C^4$-regular hypersurface for all $t\in [0,T]$.
  \item $\bigcup_{t\in [0,T]} E(t)\ssubset \Omega$.
  \item $d\in C^1_b(\omega_T)$ and $D^\gamma_xd\in C^0_b(\omega_T)$ for all $\gamma\in \N_0^n$ with $|\gamma|\leq 4$.
\end{enumerate}
Let $\Psi_\eps(\cdot,t)$, $t\in [0,T]$ denote the parametrization that are defined according to \eqref{eq:newcoordinatesnotime} with $\Gamma$ replaced by $\Gamma(t)$.
\end{assumptions}

We extend the definition of functions that are exponentially decaying to the time-dependent case and set
\begin{align*}
	X_\delta^{\mu,\Lambda}(\Omega_T):=\Big\{u\in L^\infty(\Omega_T)\,:\,
	\esssup_{x\in\Omega_T}|e^{\mu |d_\delta(x)|}u(x,t)|\leq\Lambda\Big\}
\end{align*}
and
\begin{align*}
	X(\R;\Gamma_T) := \Big\{u\in L^\infty(\R\times\omega_T)\,\big|\,\exists \Lambda,\mu>0:
	&\esssup_{(z,x,t)\in\R\times\omega_T}|e^{\mu |z|}u(z,x,t)|\leq\Lambda,\\
  &u(z,\cdot,t)\text{ is constant in normal direction}\}.
\end{align*}

We consider the modified distance functions $d_\delta$ and the cut-off functions $\eta_\delta$ as defined in Assumption \ref{ass:modified_distances} and introduce classes of phase field evolutions that we will consider in the following.

\begin{assumptions}[Phase field evolution]\label{ass:pf-evolution}
Let an evolution $(E(t))_{t\in [0,T]}$ of sets in $\Omega$, the signed distance function $d$ and $\delta>0$ as in Assumption \ref{ass:evolution} be given.

Consider an evolution of smooth phase fields $(u_\eps)_{0<\eps<\eps_0}$.
We assume that there exist $\mu\in (0,1)$, $\Lambda>0$,
and profile functions $u_j:\R\times\omega_T\to\R$, $j=0,\dots,K$, such that for all $0<\eps<\eps_0$
\begin{align}
	u_\eps &= \eta_\delta \uei +
	(1-\eta_\delta)\sgn(d)+\eps^{3}R_\eps\quad\text{ in }\Omega,
  \label{4eq:expansion-u}\\
  \uei &= \Big(\sum_{j=0}^2 \eps^ju_j\Big)\circ \Psi_\eps^{-1}
  \quad\text{ in }\{|d|<4\delta\},
  \label{4eq:uei}
\end{align}
and such that the following properties hold:
\begin{enumerate}
	\item The profile functions $u_j\in C^0(\R\times\omega_T)$, $u_j=u_j(z,x,t)$ satisfy $u_j(z,\cdot)\in C^1_b(\omega_T)$, $D^\gamma_xu_j(z,\cdot)\in C^0_b(\omega_T)$ for all $\gamma\in \N_0^n$ with $|\gamma|\leq 4$
	\begin{equation*}
		u_0-\sgn \in X(\R;\Gamma_T),\, u_j,\,|\nabla u_j|,\,\Delta u_j\in X \in X(\R;\Gamma_T)\text{ for }1\leq j\leq K.
	\end{equation*}
	\item The remainder satisfies
  $R_\eps\in X^{\frac{\mu}{\eps},\Lambda}_\delta(\Omega_T)$ for all $0<\eps<\eps_0$.
\end{enumerate}
Moreover, we assume that
\begin{equation}\label{4eq:ueps=0}
  \{u_\eps(\cdot,t)=0\} = \Gamma(t)\quad\text{ for all }t\in [0,T],\,0<\eps<\eps_0
\end{equation}
and that
\begin{equation}
  \calW_\eps(u_\eps(\cdot,0))+\calP^{\text{AG}}_\eps(u_\eps(\cdot,0))\leq C
  \label{4eq:uinit}
\end{equation}
for all $0<\eps<\eps_0$.
\end{assumptions}

We have chosen in \eqref{4eq:ueps=0} for a more restrictive setting than in the static case.
We could also have allowed for an offset between the zero level set of $u_\eps(\cdot,t)$ and $\Gamma(t)$ as in \eqref{eq:ueps=0}.
For simplicity we restrict ourselves to \eqref{4eq:ueps=0} but allow an additional contribution of order $\bigO(\eps)$ in the gradient flow equations, see \eqref{4eq:dWillFl} and \eqref{4eq:dMCF} below.

\begin{theorem}[Convergence towards the mean curvature/Willmore flow]\ \\
Consider a sequence of evolutions of smooth phase fields $(u_\eps)_{0<\eps<\eps_0}$ as in Assumption \ref{ass:pf-evolution}, satisfying an asymptotic expansion \eqref{4eq:expansion-u} with respect to an evolution $(E(t))_{t\in [0,T]}$ of sets in $\Omega$.

Assume that $u_\eps:\Omega_T\longrightarrow\R$ satisfies
\begin{align}
	\eps\partial_tu_\eps &=
  -\frac{2}{\eps^2}\Big(f'(u_\eps)\Id-\calA_\eps\Big)\big(f(u_\eps)-\calA_\eps(u_\eps)\big) + \eps R_\eps,
	\label{4eq:dWillFl}
\end{align}
with $\sup_{0<\eps<\eps_0}\|R_\eps\|_{C^0(\overline{\Omega_T})}\leq C$.

Then $\Gamma(t)$ evolves by the (rescaled) Willmore flow
\begin{align}
	\calV = -\kappa\big(\Delta_{\Gamma}H+H|\FF|^2-\half H^3\big)
	\label{eq:willmoreflowrescaled}
\end{align}
with $\kappa=2\sigma^2$ and $\sigma$ as defined in \eqref{eq:defofc0}.

\medskip
If $(u_\eps)_\eps$ instead satisfies
\begin{align}
	\eps\partial_tu_\eps &= -\frac{1}{\eps}\big(
  f(u_\eps)-\calA_\eps(u_\eps)\big) + \eps R_\eps,
	\label{4eq:dMCF}
\end{align}
then $\Gamma(t)$ evolves by the (rescaled) mean curvature flow
\begin{align}
	\calV=\sigma H.
\end{align}
\end{theorem}

\begin{proof}
We show the convergence of the diffuse Willmore flow to the Willmore flow and briefly comment on the convergence of the diffuse mean curvature flow at the end of the proof.

We will expand both sides of \eqref{4eq:dWillFl} and evaluate the identity order by order.
To identify the evolution law in the limit $\eps\to 0$ it is sufficient to consider the region $\{|d|<2\delta\}$, in which $\eta_\delta\equiv 1$.
We in particular use that the right-hand side of equation \eqref{4eq:dWillFl} is in this region to the relevant orders already determined by the inner expansion with respect to $\eps$ of $u_\eps$:
Even though $v_\eps=\calA_\eps u_\eps$ and $\calA_\eps H_\eps$ depend on the values of $u_\eps$ in the whole of $\Omega$, applying Proposition \ref{pro:expansion-v} and Lemma \ref{lem:expansion-f} shows that we only need the inner expansion of $u_\eps$ to determine the contributions in  $\{|d|<2\delta\}$.

\medskip
To expand the right-hand side of \eqref{4eq:dWillFl} we first consider $H_\eps$.
Since under the flow \eqref{4eq:dWillFl} the energy $\calW_\eps$ decreases with time and by \eqref{4eq:uinit} we obtain that $\calW_\eps(u_\eps(\cdot,t))$ is uniformly bounded.
By the calculations in the proof of Theorem \ref{thm:liminf-ansatz}, see \eqref{eq:limsup-1} we therefore deduce that
\begin{equation*}
  u_0(z,x,t) = \fraku_0(z)\quad\text{ for all }(x,t)\in\omega_T.
\end{equation*}
We deduce from Proposition \ref{pro:expansion-v}, Lemma \ref{lem:expansion-f} and \eqref{eq:Heps-K} that in  $\{|d|<2\delta\}$
\begin{align}
	H_\eps(x,t)&=H_0(z,x,t)+\eps H_1(z,x,t)+\eps^2 H_2(z,x,t) +\eps^3R^H_\eps(x,t)
	\label{eq:diffusemeancurvatureexpansion}
\end{align}
with $R^H_\eps\in X^{\frac{\mu}{\eps},\Lambda}_\delta(\Omega_T)$
and $H_0,H_1$ characterized as follows: Firstly, by \eqref{eq:Heps-K=1}
\begin{align}
  H_0 &= f'\big(\fraku_0\big)u_1-\frakA_0 \big(u_1+H\frakv_0'\big)
  = \frakL_0(u_1)-\scrA_1(\fraku_0), \label{eq:identifyh0h1}
\end{align}
where $\scrA_1=H \frakA_0^2\partial_z$ and where here and below $H_\eps, R_\eps$ are evaluated in $(x,t)$, $\fraku_0$ in $z$, $u_j$ in $(z,x,t)$ and $H$ in $(y,t)$ with $y=\Pi_{\Gamma(t)}x$.

Secondly, we derive from \eqref{eq:Heps-K} with $K=2$ and Proposition \ref{pro:expansion-v}, Lemma \ref{lem:expansion-f}
\begin{align}
	H_1 &= \Big(\frac{1}{2}f''(u_0)(u_1)^2+f'(u_0)u_2\Big)-\frakA_0 \big(u_2+Hv_1'+(\Delta-z|\FF|^2\partial_z)v_0\big)
  \nonumber\\
  &=\frakL_0(u_2)+\half f''(u_0)u_1^2-\scrA_1(u_1)-\scrA_2(u_0),
	\label{eq:diffusemeancurvatureexpansion-defs}
\end{align}
where
\begin{align}
	\scrA_2= \Delta\frakA_0^2
	-|\FF|^2\frakA_0z\partial_z\frakA_0+H^2\frakA_0^3\partial_z^2.
  \label{eq:solutionoperatorexpansion}
\end{align}
The function $H_2$ in \eqref{eq:diffusemeancurvatureexpansion} belongs to $X(\R;\Gamma_T)$ and is $C^2$-regular with respect to the $x$ variable.
We will see below that this term is not relevant for the identification of the evolution law, thus we will not include its precise characterization.

\medskip
We next consider the action of $L_\eps=f'(u_\eps)\Id-\calA_\eps$ on $H_\eps$.
Since we are only interested in the values in the region $\{|d|<2\delta\}$ we can use a Taylor expansion and the representation of $u_\eps$ in this region for the local term $f'(u_\eps)\Id$.

For the application of $\calA_\eps$ to $H_\eps$ we use analogue arguments as in Proposition \ref{pro:expansion-v} with the following  differences:
Firstly, we only have $C^2$-regularity with respect to the $x$-variable of the profile functions that represent $H_\eps$.
Therefore we obtain only $C^0$-regularity with respect to the $x$-variable for the profile functions representing $\calA_\eps H_\eps$.
Secondly, the representation of $H_\eps$ has no contribution $(1-\eta_\delta)\sgn$ that was present in the case of $u_\eps$.
This just yields that also in $\calA_\eps H_\eps$ the same term is missing.
Since we only evaluate expressions in $\{|d|<2\delta\}$ this difference does not matter.

Therefore, following the analogue computation as in the motivation of the results just after Proposition \ref{pro:expansion-v} we deduce
\begin{equation}
 L_\eps(H_\eps)= L_0(H_0)+\eps\big(L_1(H_0)
 +L_0(H_1)\big)+\eps^2\big(L_2(H_0)+ L_1(H_1)+L_0(H_2)\big)+\eps^3R_\eps
  \label{4eq:expLeps}
\end{equation}
with $L_0=\frakL_0=f'(\fraku_0)\Id-\frakA_0$ as defined in equation \eqref{2eq:def-L0},
\begin{align}
	L_1=f''(\fraku_0)u_1\Id-\scrA_1,\quad
	L_2=f''(\fraku_0)u_2\Id +\half f'''(\fraku_0)u_1^2\Id-\scrA_2.
	\label{eq:defL12}
\end{align}
and $R_\eps\in X^{\frac{\mu}{\eps},\Lambda}_\delta(\Omega_T)$.

We now expand the evolution \eqref{4eq:dWillFl}.
For the left-hand side we obtain in $\{|d|<2\delta\}$
\begin{align}
  \eps\partial_t u_\eps &= \eps\sum_{j=0}^2\big(\partial_t u_j +\frac{1}{\eps}\partial_z u_j\partial_td\big)+\calO(\eps)
	\nonumber\\
  &= \partial_z u_0\partial_td + \bigO(\eps)
  \,=\,\fraku_0'\calV+\bigO(\eps).
	\label{eq:evol-lhs}
\end{align}

We next consider the right-hand side of evolution \eqref{4eq:dWillFl}.
By \eqref{eq:diffusemeancurvatureexpansion}-\eqref{eq:diffusemeancurvatureexpansion-defs} we obtain
\begin{align}
	-\eps^{-2} L_\eps (H_\eps)&=
	-\eps^{-2} \frakL_0( H_0)+\eps^{-1}\big(\frakL_0( H_1)+L_1( H_0)\big)
  \nonumber\\
	&\qquad -\big(\frakL_0(H_2)+L_1( H_1)
	+L_2(H_0)\big)+\calO(\eps).
	\label{eq:nullteordnungzeitentwicklung}
\end{align}
To order $\eps^{-2}$ we deduce from equations \eqref{4eq:dWillFl} and \eqref{eq:evol-lhs}
that $\frakL_0( H_0)=0$, which is equivalent to
\begin{align*}
	0 = \frakL_0(\frakL_0(u_1) -H\frakA_0^2(\fraku_0')).
\end{align*}
In addition we have the condition $u_1(0)=0$. Comparing this to Sections \ref{sec:3ident}
and \eqref{eq:u1} we deduce $H_0=-H \sigma \fraku_0'$ and $u_1=H\fraku_1$. In particular,
\begin{align*}
	L_1 = H\frakL_1,~\quad\quad\quad\quad\quad\quad\quad &\scrA_1=H\frakA_1,\\
	\frakL_1=f''(\fraku_0)\fraku_1\Id-\frakA_0^2\partial_z,
	\quad &\frakA_1=\frakA_0^2\partial_z,
\end{align*}
where $\frakL_1$ and $\frakA_1$ only depend on $z$.

By evolution \eqref{4eq:dWillFl} also the $\calO(\eps^{-1})$ contribution of the right-hand side vanishes, thus
\begin{align}
	0 = L_1 (H_0)+L_0 (H_1)=-\sigma H^2 \frakL_1(\fraku_0') + \frakL_0(H_1).
	\label{eq:l1mu0gleichung}
\end{align}
We will now proceed to the crucial order $\calO(1)$ in equation \eqref{4eq:dWillFl}.
We test the corresponding equation with $\fraku_0'$
and integrate with respect to the variable $z$. We get by formulas \eqref{4eq:dWillFl} and \eqref{eq:evol-lhs}
\begin{align}
	-\half\|\fraku_0'\|_{L^2}^2\calV=\int_{\R }\fraku_0'(L_1 (H_1)+L_2( H_0))\dx^1.
	\label{eq:firstfloweq}
\end{align}

For the second term on the right-hand side of equation \eqref{eq:firstfloweq} we use
\begin{align}
	L_2 (H_0) &=(f''(\fraku_0)u_2+\half f'''(\fraku_0)u_1^2-\scrA_2)(-H \sigma \fraku_0')\nonumber\\
	&=-\sigma Hf''(\fraku_0)\fraku_0'u_2 -\frac{\sigma}{2}H^3f'''(\fraku_0)\fraku_0'\fraku_1^2 + \sigma \scrA_2 (H\fraku_0').
	\label{eq:L2H0}
\end{align}
We calculate a commutator
\begin{align*}
	[\partial,\frakL_0](w)&=
	 \big(\frakL_0(w)\big)' - \frakL_0(w')= f''(\fraku_0)\fraku_0'w
	 \quad\text{ for all }w\in L^2(\R)\cap C^1(\R)
\end{align*}
and rewrite equation \eqref{eq:diffusemeancurvatureexpansion-defs}
\begin{align}
	\frakL_0(u_2)&= w,\qquad w\coloneqq H_1 -H^2\half f''(\fraku_0)\fraku_1^2
	+H^2\frakA_1(\fraku_1)+\scrA_2(\fraku_0).\label{eq:L0u2}
\end{align}
The commutator helps us to generate $\frakL_0$ in front of $u_2$ in the right-hand side of
equation \eqref{eq:L2H0} so we can apply equation \eqref{eq:L0u2}
\begin{align*}
	\int_\R \fraku_0'f''(\fraku_0)\fraku_0'u_2\dx^1
	&= \int_\R\fraku_0'\big(w'-\frakL_0(u_2')\big)\dx^1=\int_\R\fraku_0'w'\dx^1.
\end{align*}
Together with equation \eqref{eq:L2H0} we obtain
\begin{align}
	\int_\R \fraku_0' L_2 (H_0)\dx^1
	&=-\sigma H \int_\R \fraku_0'w'\dx^1-\sigma \int_\R\fraku_0'
	\Big(\frac{1}{2}H^3 f'''(\fraku_0)
	\fraku_0'\fraku_1^2-\scrA_2 (H\fraku_0')\Big)\dx^1\nonumber\\
	&= -\sigma H \int_\R \fraku_0'\Big(H_1-H^2\half f''(\fraku_0)\fraku_1^2
	+H^2\frakA_1(\fraku_1)+\scrA_2(\fraku_0)\Big)'\dx^1\nonumber\\
	&\quad-\sigma \int_\R\fraku_0'\Big(\frac{1}{2}H^3 f'''(\fraku_0)
	\fraku_0'\fraku_1^2-\scrA_2 (H\fraku_0')\Big)\dx^1.\nonumber
\end{align}
By calculating $w'$ many terms cancel out.
\begin{align}
	\int_\R \fraku_0' L_2 (H_0)\dx^1
	&= \sigma H\int_\R \Big(\fraku_0'' H_1  +\fraku_0'H^2f''(\fraku_0)\fraku_1\fraku_1'
	-H^2\fraku_0'\frakA_1(\fraku_1')\Big)\dx^1\label{eq:2}\\
	&\quad+\sigma\int_\R\fraku_0'\big[\frakA_2,H\partial_z\big](\fraku_0)\dx^1\nonumber
\end{align}
Calculating the commutator yields for any $w\in L^2(\R)\cap C^1(\R)$
\begin{align*}
	\big[\frakA_2,H\partial_z\big](w)=\frakA_1(w)(\Delta_{\Gamma} H-H|\FF|^2).
\end{align*}
Using this, $\int_\R|\frakA_0(\fraku_0')|^2\dx^1=c_0$ and the definition of $\frakL_1$ we get
\begin{align}
	\int_\R \fraku_0' L_2 (H_0)\dx^1
	&= \sigma H\int_\R\fraku_0'' H_1 \dx^1 +\sigma H^3
	\int_\R\fraku_0'\frakL_1(\fraku_1')\dx^1
	+c_0\sigma(\Delta_{\Gamma}H+H|\FF|^2)\label{eq:rhs2}.
\end{align}
By differentiating formula \eqref{eq:eulerlagrangeu1} we have
\begin{align}
	-\sigma \fraku_0'' &= f''(\fraku_0)\fraku_0'\fraku_1 + f'(\fraku_0)\fraku_1'
	-\frakA_0(\fraku_1') -\frakA_1(\fraku_0')= \frakL_1(\fraku_0') + \frakL_0(\fraku_1').
	\label{eq:sigmau0}
\end{align}
By multiplying the defining equation for $\scrA_0(h)$ with $z$ and comparing it
with the equation for $\scrA_0(zh)$ for smooth $h$ we can identify the commutator
$\big[\frakA_0,z\big]=2\frakA_1$. We conclude
\begin{align}
	0 &= z\frakL_0(\fraku_0')=f'(\fraku_0)z\fraku_0'-z\frakA_0(\fraku_0')
	=\frakL_0(z\fraku_0')+2\frakA_1(\fraku_0').
	\label{eq:weiterel0formel}
\end{align}
Before moving to the final calculations we need the antisymmetrical
part of $\frakL_1$. We have for $w_1,w_2\in L^2(\R)$
\begin{align}
	\int_\R\big(w_1\frakL_1(w_2)-w_2\frakL_1(w_1)\big)\dx^1
	=-\int_\R2w_2\frakL_1^{\text{as}}(w_1)\dx^1=\int_\R2w_2\frakA_1(w_1)\dx^1.
	\label{eq:l1as}
\end{align}
We consider the sum of the contributions from the first terms of the right-hand side
of formulas \eqref{eq:firstfloweq} and \eqref{eq:rhs2} to deduce
\begin{align}
	H\int_\R\Big(\fraku_0' \frakL_1(H_1) +\sigma\fraku_0'' H_1 \Big)\dx^1
	&\stackrel{\mathclap{\mathrm{\eqref{eq:sigmau0}}}}=
	H\int_\R\Big(\fraku_0' \frakL_1(H_1)-H_1
	\big(\frakL_1(\fraku_0')+\frakL_0(\fraku_1')\big)\Big)\dx^1\nonumber\\
	&\stackrel{\mathclap{\mathrm{\eqref{eq:l1as}}}}
	=H\int_\R H_1\Big(2\frakA_1(\fraku_0')-\frakL_0(\fraku_1')\Big)\dx^1\nonumber\\
	&\stackrel{\mathclap{\mathrm{\eqref{eq:weiterel0formel}}}}
	= H \int_\R H_1  \frakL_0(-z\fraku_0'-\fraku_1')\dx^1\nonumber\\
	&\stackrel{\mathclap{\mathrm{\eqref{eq:l1mu0gleichung}}}}
	=-\sigma H^3 \int_\R(z\fraku_0'+\fraku_1')\frakL_1(\fraku_0')\dx^1.
	\label{eq:guteszwischenergebnis}
\end{align}
Plugging the results from equations \eqref{eq:guteszwischenergebnis}
and \eqref{eq:rhs2} into the identity \eqref{eq:firstfloweq} we obtain
\begin{align*}
	-\calV=2\sigma^2\big(\Delta_{\Gamma}H + H|\FF|^2+\frac{\kappa_1}{c_0}H^3\big),
\end{align*}
with
\begin{align*}
	\kappa_1 &= \int_\R\Big(\fraku_0'\frakL_1(\fraku_1')
	-(z\fraku_0'+\fraku_1')\frakL_1(\fraku_0')\Big)\dx^1\\
	&\stackrel{\mathclap{\mathrm{\eqref{eq:l1as}}}}
	= \int_\R\Big(2\fraku_1'\frakA_1(\fraku_0')-z\fraku_0'\frakL_1(\fraku_0')\Big)\dx^1\\
	&\stackrel{\mathclap{\mathrm{\eqref{eq:weiterel0formel}}}}
	=-\int_\R z\fraku_0'\big(\frakL_1(\fraku_0')+\frakL_0(\fraku_1')\big)\dx^1
	\stackrel{\mathclap{\mathrm{\eqref{eq:sigmau0}}}}= -\frac{c_0}{2},
\end{align*}
which proves the Willmore-flow equation.

\medskip
Let us now consider the evolution \eqref{4eq:dMCF}. The $\calO(\eps^{-1})$ order
gives $f(u_0)=\calA_0(u_0)$ and thus $u_0=\fraku_0$ as before.
We expand the right-hand side of \eqref{4eq:dMCF} and get in $\{|d|<2\delta\}$
\begin{align*}
	-\frac{1}{\eps}\big(f(u_\eps)-\calA_\eps(u_\eps)\big)
	&=-f'(\fraku_0)u_1+\frakA_0(u_1)+\frakA_1(\fraku_0)+\calO(\eps)\\
	&=-\frakL_0(u_1)+\frakA_1(\fraku_0)+\calO(\eps).
\end{align*}
The expansion of the left-hand side of \eqref{4eq:dMCF} has already be computed in \eqref{eq:evol-lhs}.
Equating the two and testing with $\fraku_0'\in\ker(\frakL_0)$ yields
\begin{align*}
	\calV\int_\R\big|\fraku_0'\big|^2\dx^1 &=\int_\R\fraku_0'\frakA_1(\fraku_0)\dx^1
\end{align*}
Taking the definition \eqref{eq:defofc0} into account we get the evolution by mean curvature
\begin{align*}
	\calV=\sigma H.
\end{align*}
\end{proof}

\section{Numerical simulations}\label{sec:numerics}
The goal of this section is to use the new approximation scheme for numerical simulations.
We consider gradient flows for linear combinations of perimeter and Willmore
functional in the plane. For $\lambdaCH_1,\lambdaCH_2\geq 0$ the corresponding
evolution law is given by (see \eqref{eq:mcf},\eqref{eq:willmoreflow})
\begin{equation*}
 \calV=-\lambdaCH_1\big(\Delta_{\Gamma(t)}H + H|\FF|^2-\frac{1}{2}H^3\big)+\lambdaCH_2 H.
\end{equation*}
The standard Cahn--Hilliard/De Giorgi approximation for this flow,
with double well potential $\WCH$, is given by
\begin{align}
  \eps \partial_t \uCH \,&=\, -\frac{\lambdaCH_1}{\eps}\Big(-\eps\Delta
  + \frac{1}{\eps}\WCH''(\uCH)\Big) \big(-\eps\Delta \uCH+
  \frac{1}{\eps}\WCH'(\uCH)\big) -\lambdaCH_2\big(-\eps\Delta \uCH+
  \frac{1}{\eps}\WCH'(\uCH)\big), \label{eq:Wf-ac}
\end{align}
whereas the equation for the corresponding approximation with
respect to the gradient-free model reads
\begin{align}
 \eps\partial_tu = -\frac{\lambda_1}{\eps^2}\Big(f'(u)\Id-\scrA_\eps\Big)
 \frac{1}{\eps}\big(-\scrA_\eps u + f(u)\big) - \frac{\lambda_2 }{\eps}
 \big(-\scrA_\eps u + f(u)\big),
 \label{eq:Wf-gf}
\end{align}
with $\lambda_1=\frac{\lambdaCH_1}{\sigma^2}$, $\lambda_2=\frac{\lambdaCH_2 }{\sigma}$,
$f(r)=r + \frac{1}{2}W'(r)$, and $\scrA_\eps=(-\eps^2\Delta + \Id)^{-1}$
as in \eqref{eq:v_eps}.

We use an implicit spectral discretization and follow the implemention
proposed by Bretin, Masnou and Oudet \cite{BretinMasnouOudet2015}.
We see that the new scheme leads to results that are in most cases quite
similar to the standard diffuse approximation.
However, a different behavior can be observed in certain cases.
Depending on the context and the available numerical methods the new
approach might offer advantages.

For a meaningful comparison we need to adjust the choice of the double
well potential in both models.
We therefore require that the optimal profile $q$ for the standard diffuse
model with respect to a double well potential $\WCH$ coincides with the
optimal profile $\fraku_0$ of the gradient-free model with double well potential $W$.

One possibility is to choose for the standard diffuse approximation the
quartic potential $r\mapsto \frac{1}{8}(1-r^2)^2$,
which leads to the optimal profile $\tanh(\cdot/4)$.
Using the characterization of $\fraku_0$ in Theorem \ref{thm:1dminimizer}
a double well potential $W$ for the gradient-free model can be computed
such that $\fraku_0=\tanh(\cdot/4)$ holds.
The resulting double well potential $W$ belongs to $C^\infty((-1,1))\cap C^2([-1,1])$
and can be extended to a $C^2$-function with quadratic growth.

Other choices lead to a compatible pair of smooth double well potentials.
In our simulations we choose the optimal profile $q=\fraku_0$ that is given by
\begin{align*}
  q(x) = \frakv_0(x) - \frakv_0''(x)
\end{align*}
with
\begin{align*}
  \frakv_0(z) \coloneqq \tanh\Big(\frac{1}{4} z\Big).
\end{align*}
The corresponding double well potentials $\WCH$, $W$ for the standard and the gradient-free approximation then belong to $C^\infty([-1,1])$ and can be extended to satisfy all assumptions that were imposed in the previous sections.
Furthermore, the double well potential can be computed from
\begin{equation*}
 W'(r) = 2\frakv_0''(q^{-1}(r)),\qquad
 \WCH'(r) = q''(q^{-1}(r)).
\end{equation*}
For the constants $c_0,\sigma$ in equation \eqref{eq:defofc0} we find
\begin{equation*}
 c_0 = \frac{1}{3},\quad \sigma = \frac{69}{560}.
\end{equation*}

\begin{figure}
\hspace{10cm}
\includegraphics[scale=0.4]{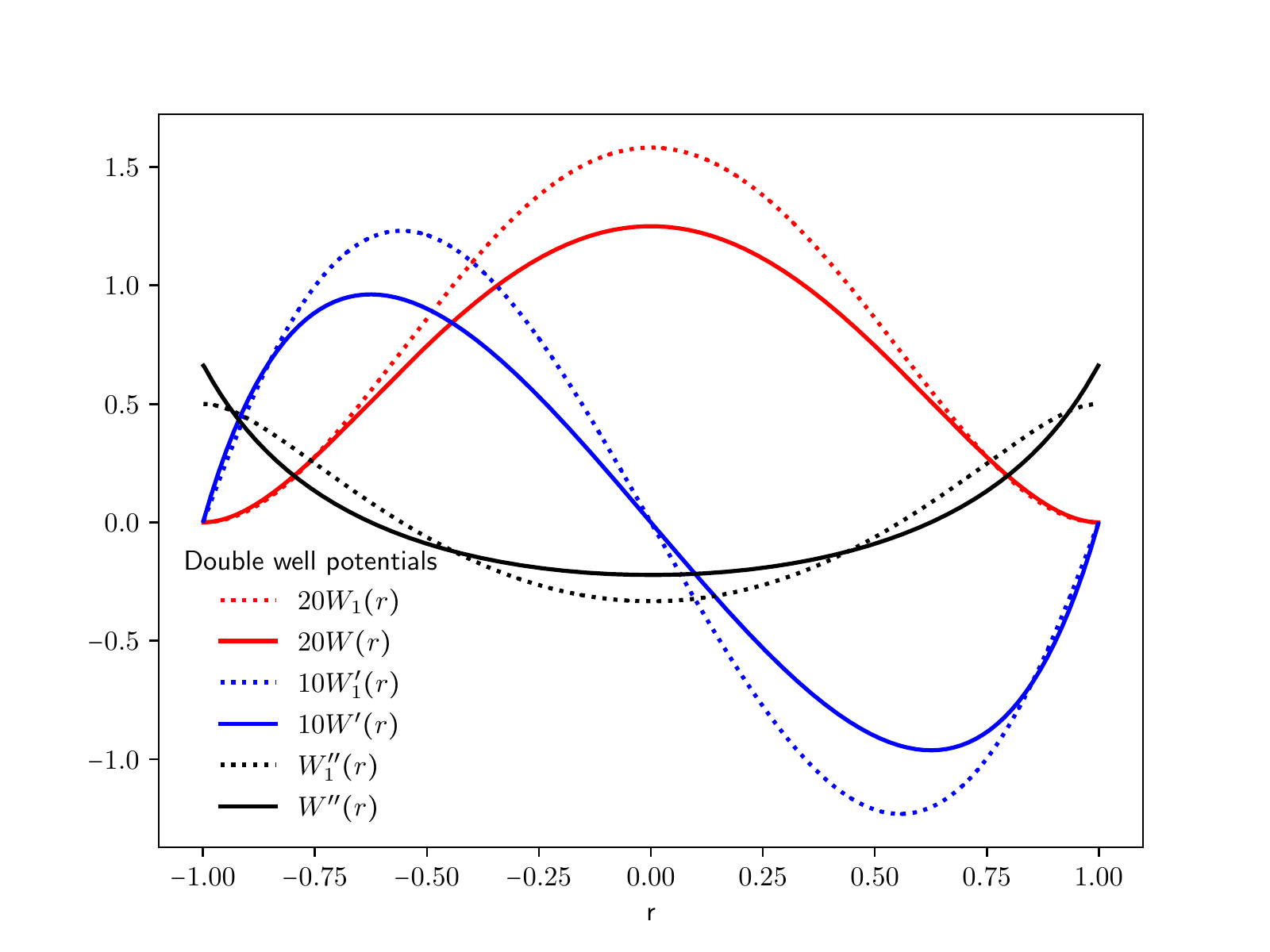}
 \caption{Comparison of the standard DWP $W_1$ and the new function $W$}
 \label{fig:plot_new_W}
\end{figure}

\subsection{Implicit spectral discretization}

We adapt the scheme of \cite{BretinMasnouOudet2015} based on the splitting
\begin{align*}
  \eps\partial_t u &= -\frac{\lambda_1}{\varepsilon^2} \Big(f'(u)\Id-\scrA_\eps\Big) H_\eps - \lambda_2 H_\eps
  = -\frac{\lambda_1}{\varepsilon^2} \Big(B_\eps +\frac{1}{2}W''(u)\Id\Big) H_\eps - \lambda_2 H_\eps \\
  H_\eps &= \frac{1}{\eps}\big(-\scrA_\eps u + u + \frac{1}{2}W'(u)\big)
  = \frac{1}{\eps}\big(B_\eps u + \frac{1}{2}W'(u)\big),
\end{align*}
with $B_\eps\coloneqq\Id-\scrA_\eps$.

We assume periodic boundary conditions for $u$ and $H_\eps$ and approximate the functions by finite Fourier sums, i.e.\ $u = \sum_{\mathbf{k} \in \Z_N^n} \widehat{u}_{\mathbf{k}}  e_{\mathbf{k}}$ where $(e_{\mathbf{k}})_{\mathbf{k} \in \Z^n}$ is an orthogonal basis of $L^2(\T^n)$ given by the eigenvectors of $-\Delta$ and $\Z_N^n \coloneqq \{ \mathbf{k} \in \Z^n \mid |\mathbf{k}|_\infty \le N \}$.

An implicit discretization in time with step size $\delta_t > 0$ yields
\begin{align*}
  u^{n+1} &= -\delta_t \left( \frac{\lambda_1}{\varepsilon^3} B_\eps H_\eps^{n+1} + \frac{\lambda_1}{2\varepsilon^3} W''(u^{n+1}) H_\eps^{n+1} + \frac{\lambda_2}{\varepsilon} H_\eps^{n+1} \right) + u^n \\
  &\eqqcolon -\frac{\delta_t \lambda_1}{\varepsilon^3} B_\eps H_\eps^{n+1} + E \\
    \eps H_\eps^{n+1} &= B_\eps u^{n+1} + \frac{1}{2} W'(u^{n+1})\eqqcolon B_\eps u^{n+1}+F
\end{align*}
which is equivalent to
\begin{align*}
  u^{n+1} + \frac{\delta_t \lambda_1}{\varepsilon^3} B_\eps \left(B_\eps u^{n+1} + F \right) &= E \\
  \eps H_\eps^{n+1} - B_\eps \left( -\frac{\delta_t \lambda_1}{\varepsilon^3} B_\eps H_\eps^{n+1} + E \right) &= F.
\end{align*}
This can be simplified to
\begin{align*}
  \left(\Id +  \frac{\lambda_1\delta_t}{\varepsilon^4}B_\eps^2\right) u^{n+1} &= E - \frac{\delta_t \lambda_1}{\varepsilon^4} B_\eps F \\
  \eps\left(\Id +  \frac{\lambda_1\delta_t}{\varepsilon^4}B_\eps^2\right) H_\eps^{n+1} &= F + B_\eps E.
\end{align*}

After introducing the mapping
\begin{align*}
  \phi(u, \eps H_\eps) \coloneqq \left(\Id + \frac{\lambda_1 \delta_t}{\varepsilon^4} B_\eps^2\right)^{-1}
  \begin{pmatrix}
   \Id & -\frac{\delta_t \lambda_1}{\varepsilon^4} B_\eps \\ B_\eps & \Id
  \end{pmatrix}                                                   \begin{pmatrix}
   - \frac{\lambda_1\delta_t}{2\varepsilon^4}W''(u)H_\eps - \frac{\lambda_2\delta_t}{\varepsilon^2}H_\eps + u^n \\ \frac{1}{2} W'(u)
  \end{pmatrix}
\end{align*}
the above equation is equivalent to
\begin{align*}
  \begin{pmatrix} u^{n+1} \\ \eps H_\eps^{n+1} \end{pmatrix} &= \phi(u^{n+1}, \eps H_\eps^{n+1}).
\end{align*}
In each timestep this equation can be solved with a fixed point iteration, because $\phi$ is a contraction on a suitable subset of $L^2 \times L^2$ under appropriate assumptions on $\delta_t$, $\varepsilon$ and $N$.

We exploit that the operator $B_\eps$ is a multiplication operator in Fourier space, i.e. $(B_\eps u^{n+1})^{\widehat{}}_{\mathbf{k}} = \widehat{(B_\eps)}_{\mathbf{k}} \widehat{u}_{\mathbf{k}}$. More precisely, we have
\begin{align*}
  (B_\eps)^{\widehat{}}_{\mathbf{k}} &=  1 - \frac{1}{(-\varepsilon^2\Delta + \Id)^{\widehat{}}_{\mathbf{k}}}
  =1 - \frac{1}{1 + 4 \pi^2 \varepsilon^2 |\mathbf{k}|^2}
  = \frac{4\pi^2\varepsilon^2 |\mathbf{k}|^2}{1 + 4\pi^2 \varepsilon^2 |\mathbf{k}|^2}.
\end{align*}
Note that
\begin{align*}
  \left(\left(\Id + \frac{\lambda_1 \delta_t}{\varepsilon^4} B_\eps^2\right)^{-1}\right)^{\widehat{}}_{\mathbf{k}} &= \frac{1}{1 + \frac{16\pi^4 \lambda_1 \delta_t |\mathbf{k}|^2}{\left(1 + 4\pi^2\varepsilon^2 |\mathbf{k}|^2\right)^2}}.
\end{align*}

\subsubsection{Constrained problem}
A volume constraint can be naturally implemented in the Fourier approach. In addition, we propose a new numerical scheme to enforce a constant perimeter.

For the volume constraint, we use the orthogonal projection $\proj : L^2(\T^n) \to  \{u \mid \widehat{u}_{\mathbf{0}} = 0 \}$ to keep the first Fourier coefficient and therefore the volume constant.
This leads to the following modification $\tilde{\phi}$ of $\phi$
\begin{align*}
  \tilde{\phi}(u, \eps H_\eps) \coloneqq \left(\Id + \frac{\lambda_1 \delta_t}{\varepsilon^4} B_\eps^2\right)^{-1}
  \begin{pmatrix}
   \proj & -\frac{\delta_t \lambda_1}{\varepsilon^4} B_\eps \\ B_\eps & \Id
  \end{pmatrix}                                                   \begin{pmatrix}
   - \frac{\lambda_1\delta_t}{2\varepsilon^4}W''(u)H_\eps - \frac{\lambda_2\delta_t}{\varepsilon^2}H_\eps + u^n \\ \frac{1}{2} W'(u)
 \end{pmatrix}.
\end{align*}

To enforce the perimeter constraint, we use the solution $u^{n+1}$ of the fixed-point iteration.
Starting with this function and the initial diffuse perimeter $c \coloneqq \calP^{\text{AG}}_\eps(u^0)$, we solve the diffuse approximation of the mean curvature flow
\begin{align*}
  \mathfrak{s} \eps \partial_t \mathfrak{u} = \frac{1}{\sigma \eps}\scrA_\eps \mathfrak{u} - \frac{1}{\sigma \eps} \mathfrak{u} - \frac{1}{2\sigma \eps} W'(\mathfrak{u})
\end{align*}
with $\mathfrak{u}(0) \coloneqq u^{n+1}$ and $\mathfrak{s} \coloneqq \sgn ( \calP^{\text{AG}}_\eps(u^{n+1}) - c)$. Note that this equation is ill-posed for $\mathfrak{s} = -1$, i.e. when the diffuse perimeter of $u^{n+1}$ is smaller than the diffuse perimeter of the initial condition. In this situation, the equation becomes a diffuse approximation of the time-reversed mean curvature flow.

Let $\bar{u}^{n+1} \coloneqq \mathfrak{u}(\mathfrak{t})$ where $\mathfrak{t} \coloneqq \inf \{ \tau \ge 0 \mid \calP^{\text{AG}}_\eps(\mathfrak{u}(\tau)) = c \}$.
Formally, $\bar{u}^{n+1}$ is an approximation of the orthogonal projection of $u^{n+1}$ to the set $\{ v \mid \calP^{\text{AG}}_\eps(v) = c \}$ of functions with diffuse perimeter equal to the diffuse perimeter of the initial condition.

In order to simulate $\mathfrak{u}$ we use a semi-implicit scheme in Fourier space for the diffuse mean curvature flow. With a binary search we find $\mathfrak{t}$ within in a tolerance of $|\calP^{\text{AG}}_\eps(\mathfrak{u}({\mathfrak{t}})) - c| < 10^{-7}$.

This approximation is in general not well-posed. However, the numerical experiments are stable in two and three dimensions, and yield the expected shapes.

\subsection{Results}
Here we present some numerical simulations for the standard diffuse approximation and the gradient-free model, where we use the scheme introduced by \cite{BretinMasnouOudet2015} for the former and the modified scheme as described above for the latter.

As the computational domain we choose the cube $[-\frac{1}{2},\frac{1}{2}]^n$ with periodic boundary conditions.
Initial data are chosen as $u^0 = q(\frac{\sdist}{\varepsilon})$ where $\sdist$ is the signed distance to the boundary of the set we are interested in, taken positive inside the set.

The timestep size is chosen of order $\eps^2$ for simulations of the mean curvature flow and of order $\eps^4$ for all other simulations.

\subsubsection{Verifying the theoretical results}
\begin{figure}
\hspace{10cm}
 \resizebox{\textwidth}{!}{\import{}{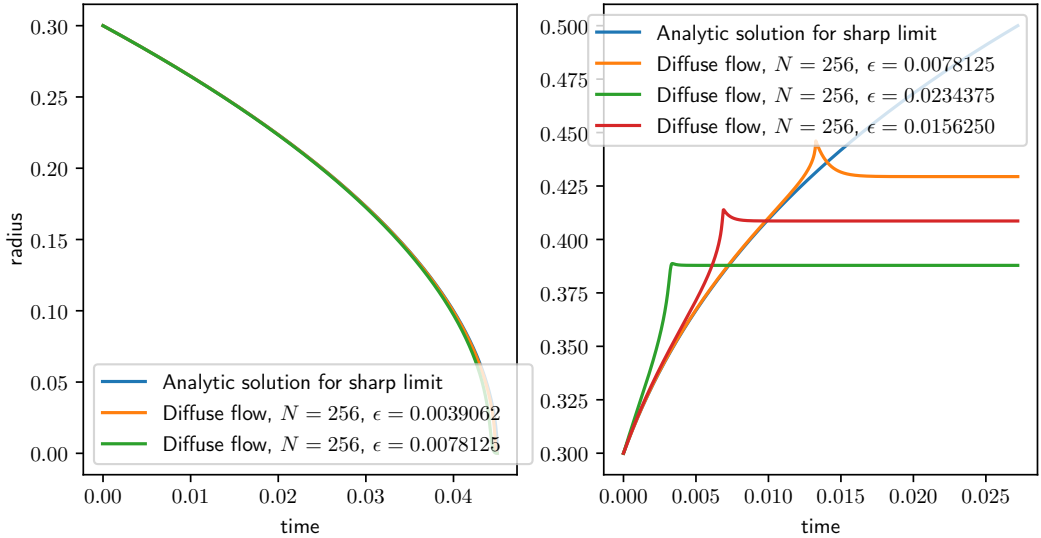}}
 \caption{Mean curvature and Willmore flow for the standard approximation
  }
 \label{fig:bench1a}
\end{figure}

\begin{figure}\centering
 \resizebox{\textwidth}{!}{\import{}{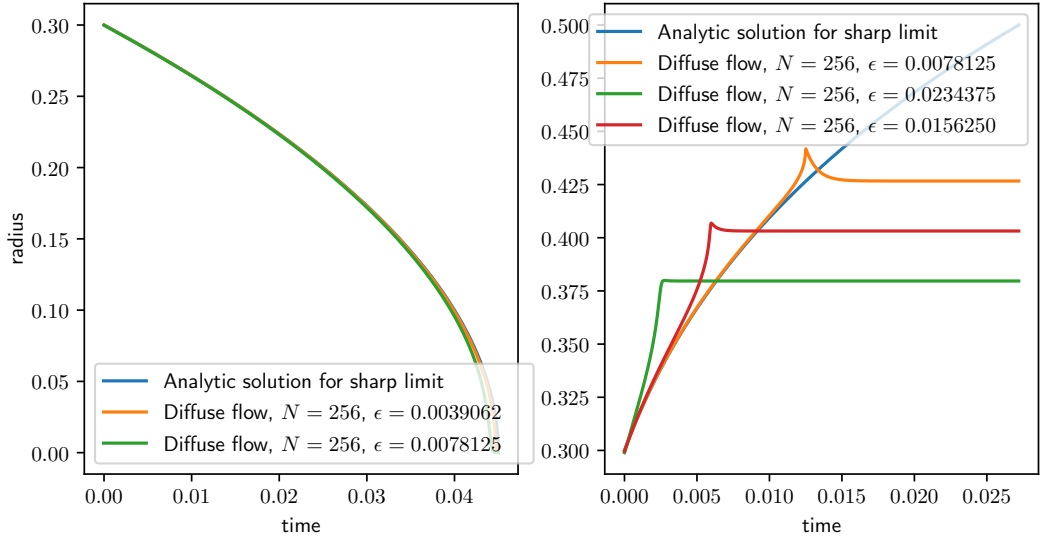}}
 \caption{Mean curvature and Willmore flow for the new approximation
  }
 \label{fig:bench1b}
\end{figure}

In order to verify the theoretical results we start with a sphere of prescribed radius and compare the simulations to the known sharp interface evolution. In order to compute the radius of the diffuse approximation, we compute the diffuse perimeter and solve for the radius.
In \Cref{fig:bench1a}, \Cref{fig:bench1b} we have plotted the radius (vertical axis) as a function of time (horizontal axis) for the diffuse mean curvature and diffuse Elastica flow, using the standard and the new approximation.
We see a very good agreement with the analytic solution and hardly any differences between the two diffuse approximations.

%
%

\begin{figure}\centering
 \begin{tikzpicture}[trim left=-10.6em,scale=1]
  \node (f1) at (0,0)
  {\includegraphics*[height=0.3\textheight]%
      {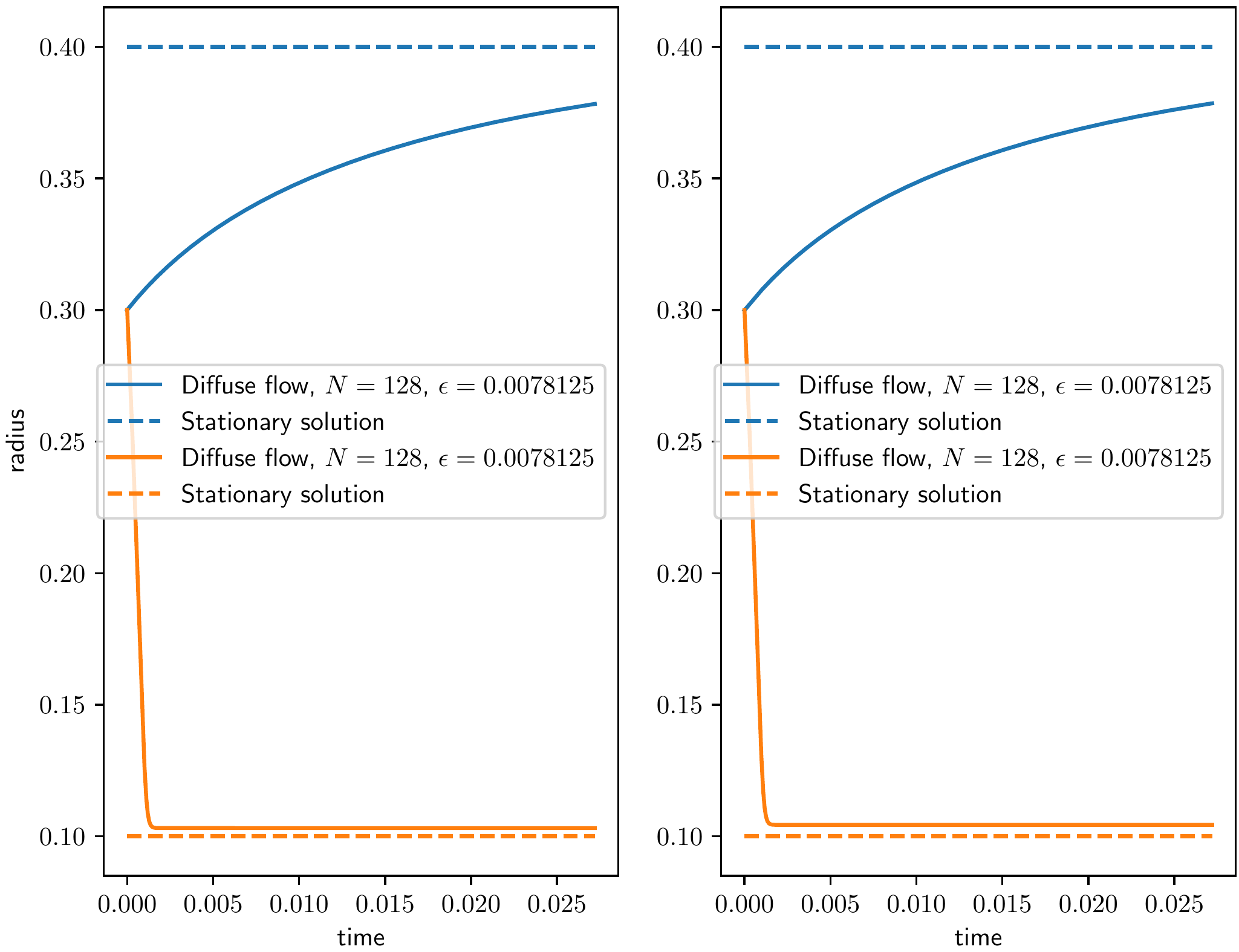}
  };
  \node (f3) [right=0.6ex of f1]
  {\includegraphics[height=0.293\textheight]{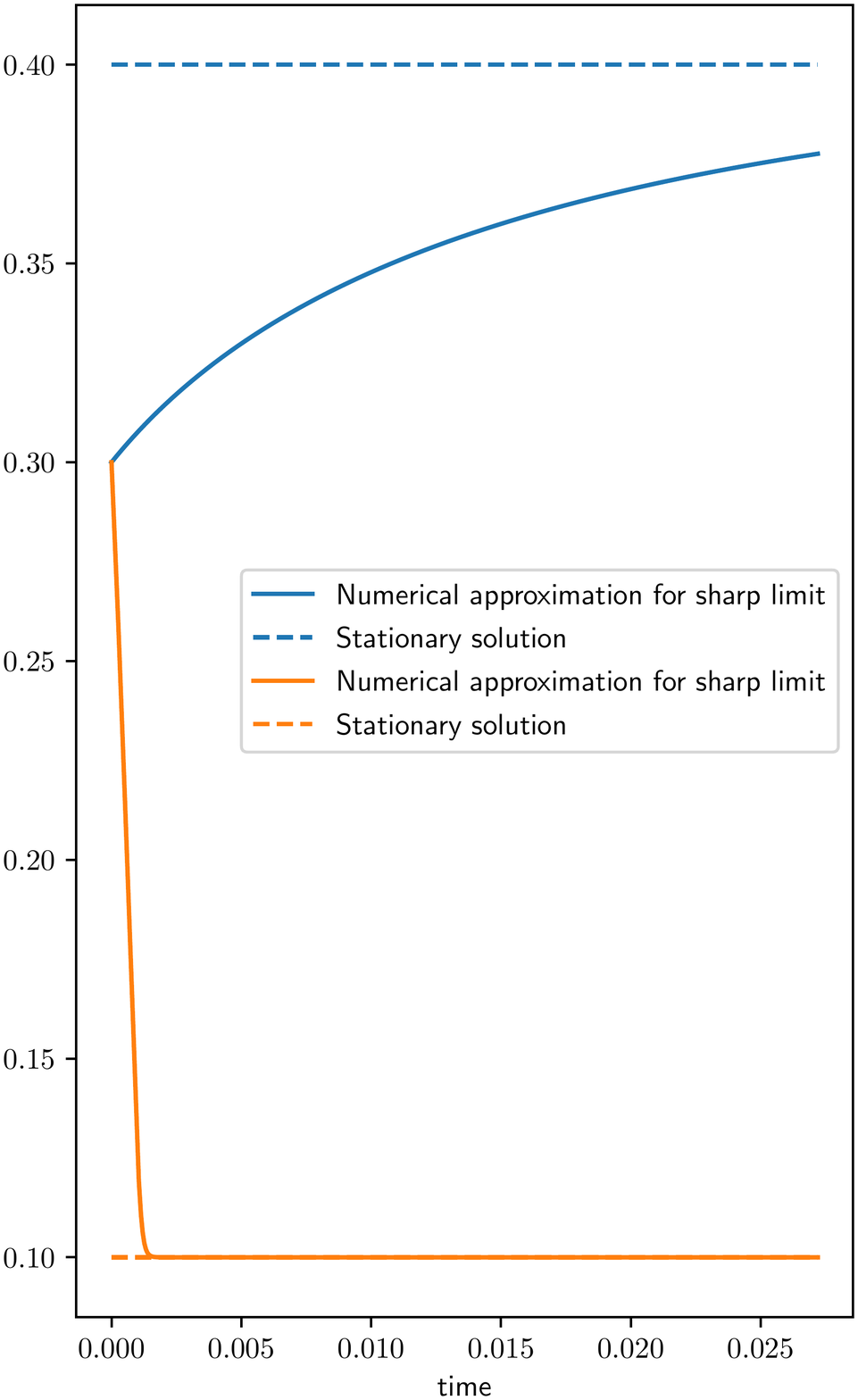}
  };
 \end{tikzpicture}
 \caption{Diffuse flow for $\lambdaCH_1=1$ and $\lambdaCH_2=50$, $\lambdaCH_2=\frac{25}{8}$. From left to right: The standard approximation, the new approximation, and a numerical approximation of the sharp interface ODE.}
 \label{fig:bench2}
\end{figure}

For the next benchmark we again consider a radially symmetric situation.
Here we consider two different linear combinations of Willmore and perimeter
functional such that respective stationary states have larger and smaller radius than
the initial states. The sharp interface evolution is governed by an ODE.
In \Cref{fig:bench2} we show corresponding approximations by the standard
and gradient-free diffuse approximation and a numerical approximation of the
sharp interface ODE. Again we plot the radius as a function of time and see a very
good agreement of both diffuse approximations with the sharp interface dynamics.

\subsection{Numerical experiments}
Here we use the gradient-free model for numerical simulations for the Elastica flow and particular initial conditions.
In some of the experiments we in addition consider constraints on volume and/or perimeter.
In the simulations we have chosen $N=256$ and $\eps=1/N$, except for the fourth example below, where $N=512$ and $\eps=1/(2N)$.

\begin{figure}
 \centering
   \begin{tikzpicture}[trim left=-4.2em]
    \node (f1) at (0,0) {\includegraphics[scale=0.9, width=0.2\textwidth]{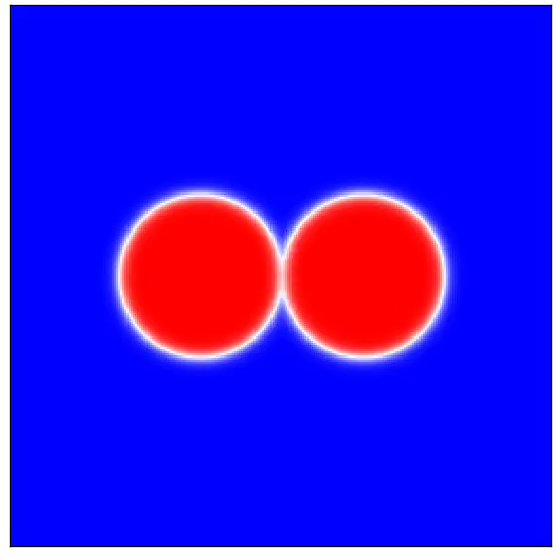}
    };
    \node (f2) [right=0.1ex of f1]{\includegraphics[scale=0.9, width=0.2\textwidth]{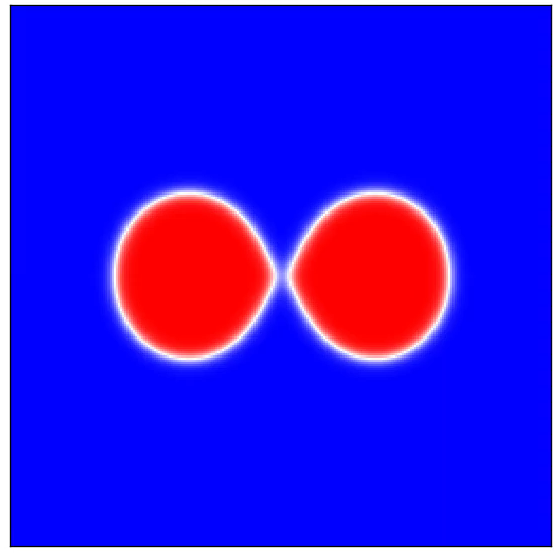}
    };
    \node (f3) [right=0.1ex of f2] {\includegraphics[scale=0.9, width=0.2\textwidth]{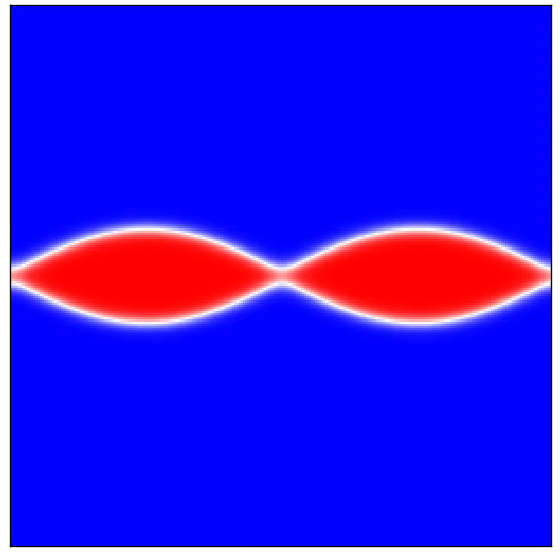}};
    \node (f4) [right=0.1ex of f3]{\includegraphics[scale=0.9, width=0.2\textwidth]{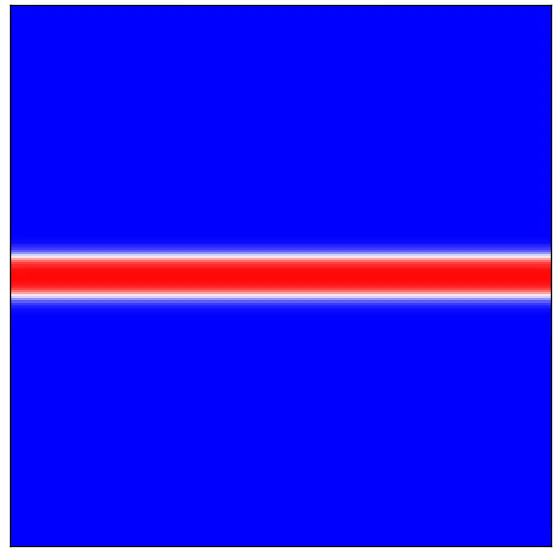}
    };
    \node (f5) [below=0.1ex of f1] {\includegraphics[scale=0.9, width=0.2\textwidth]{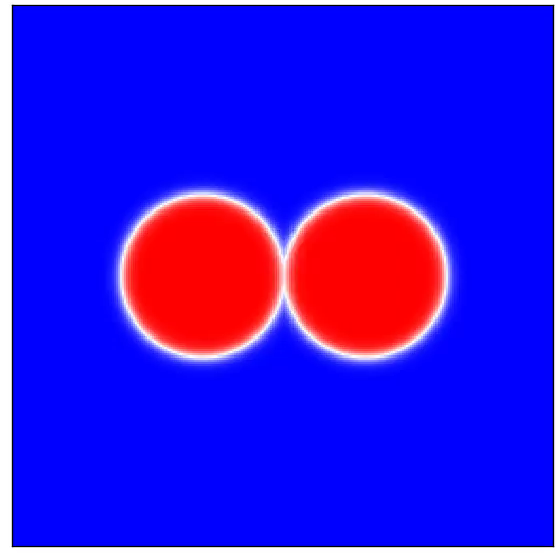}
    };
    \node (f6) [right=0.1ex of f5]{\includegraphics[scale=0.9, width=0.2\textwidth]{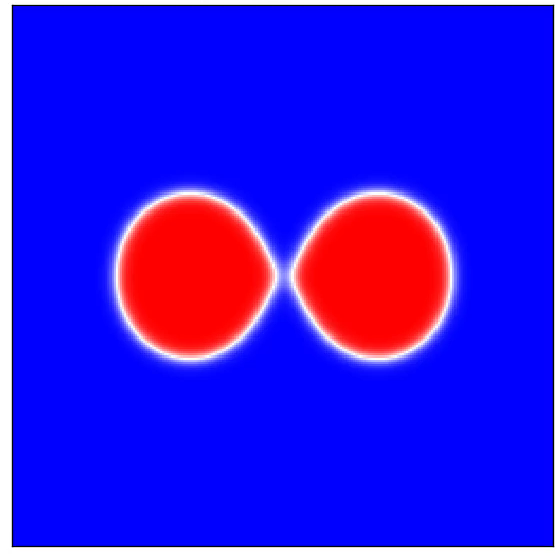}
    };
    \node (f7) [right=0.1ex of f6] {\includegraphics[scale=0.9, width=0.2\textwidth]{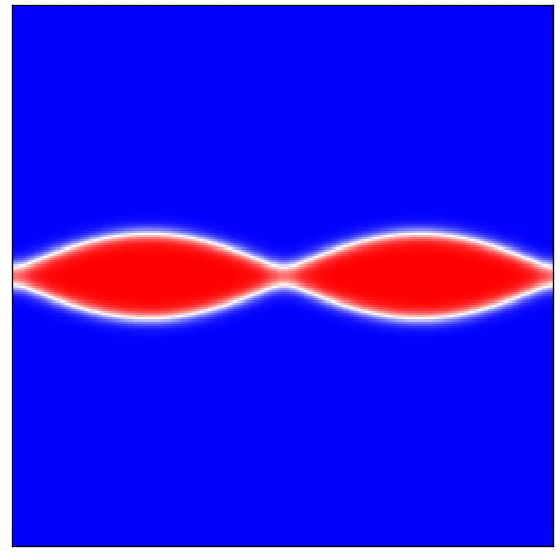}};
    \node (f8) [right=0.1ex of f7]{\includegraphics[scale=0.9, width=0.2\textwidth]{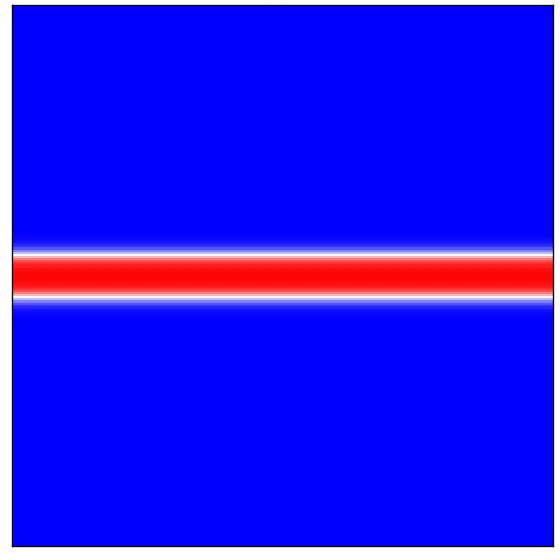}
    };
   \end{tikzpicture}
    \caption{Evolution by diffuse Elastica flow. Discrete
    phase-fields $u(\cdot,t)$ at times 0, 0.000046, 0.001839 and 0.003679. The top row shows the standard diffuse and the bottom row the new diffuse approximation}
    \label{fig:new-touch-circ-a}
\end{figure}

The first example considers the Elastica flow, i.e.~$\lambdaCH_1=1$, $\lambdaCH_2=0$ for two touching circles (see Figure \ref{fig:new-touch-circ-a}). The touching point stays fixed for some time, a transversal intersection develops, followed by a subsequent flattening of the structure and finally an evolution towards a strip with flat boundaries, which clearly is a global minimizer of the Elastica energy.

\begin{figure}
 \centering
  \begin{tikzpicture}[trim left=-3.4em]
   \node (f0) at (0,0) {\includegraphics[scale=0.9, width=0.15\textwidth]{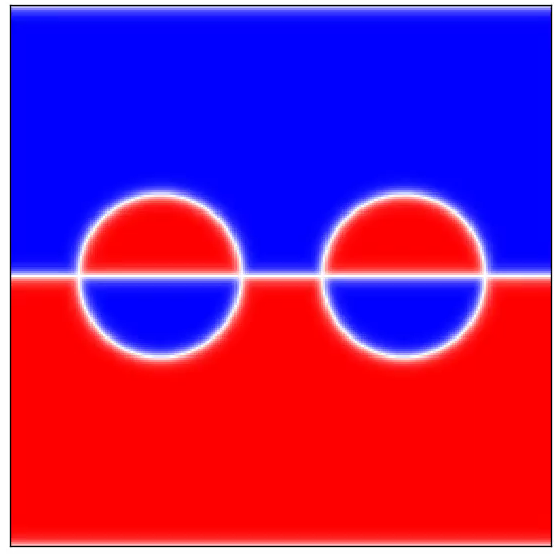}
   };
   \node (f1) [right=0.05ex of f0] {\includegraphics[scale=0.9, width=0.15\textwidth]{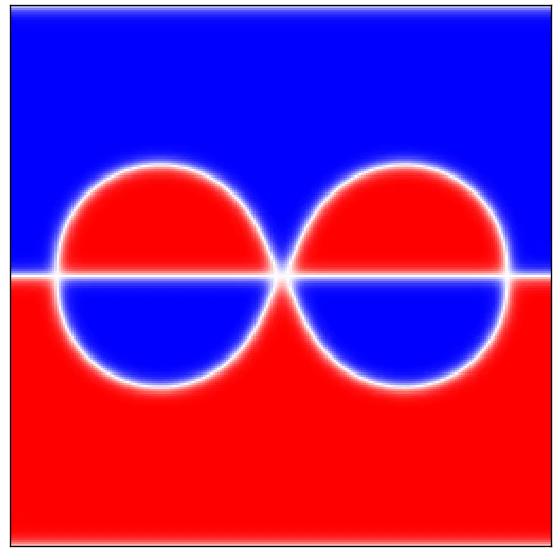}
   };
   \node (f3) [right=0.05ex of f1]{\includegraphics[scale=0.9, width=0.15\textwidth]{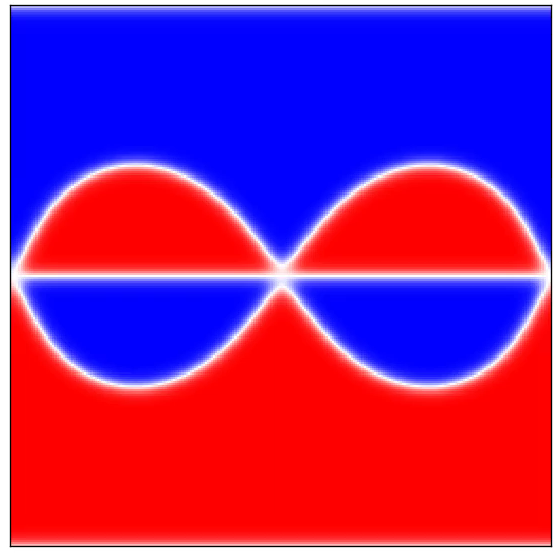}
   };
   \node (f4) [right=0.05ex of f3] {\includegraphics[scale=0.9, width=0.15\textwidth]{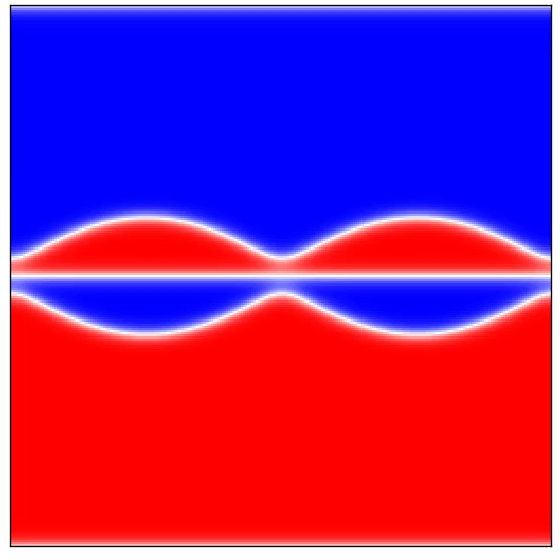}};
   \node (f5) [right=0.05ex of f4]{\includegraphics[scale=0.9, width=0.15\textwidth]{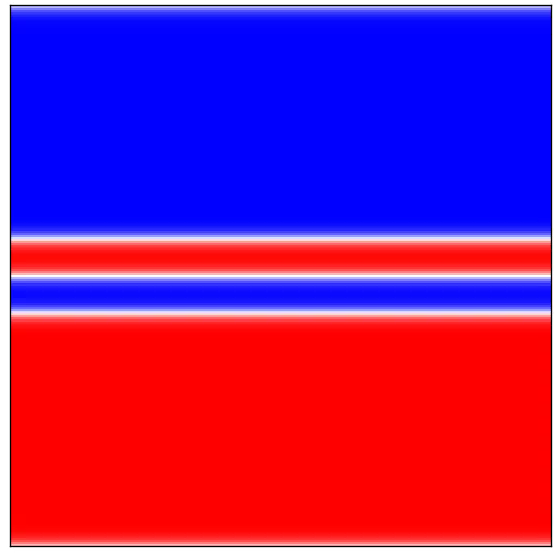}
   };
   \node (f7) [below=0.1ex of f0]{\includegraphics[scale=0.9, width=0.15\textwidth]{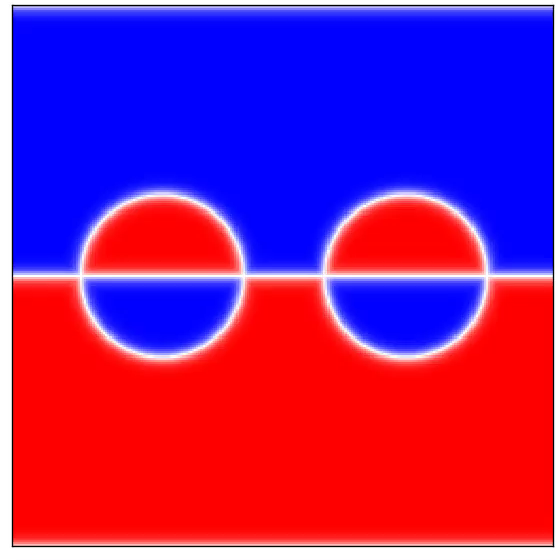}
   };
   \node (f8) [right=0.05ex of f7]{\includegraphics[scale=0.9, width=0.15\textwidth]{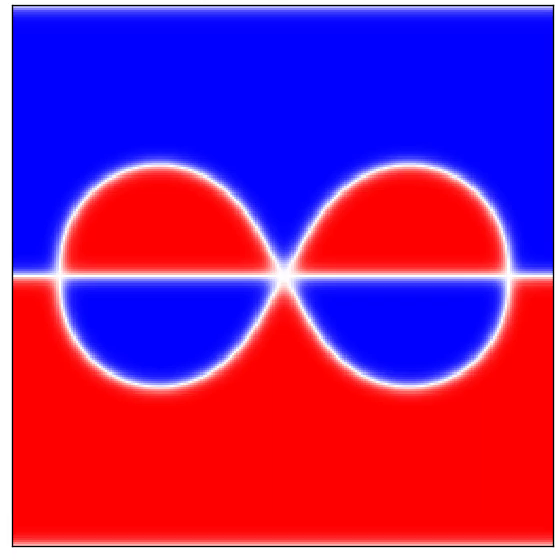}
   };
   \node (f10) [right=0.05ex of f8]{\includegraphics[scale=0.9, width=0.15\textwidth]{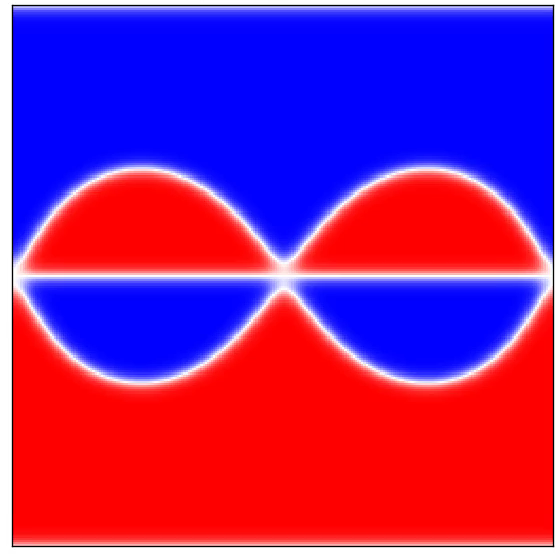}
   };
   \node (f11) [right=0.05ex of f10] {\includegraphics[scale=0.9, width=0.15\textwidth]{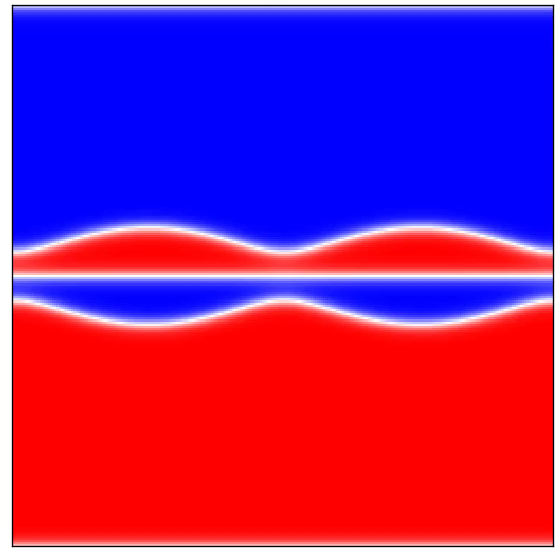}};
   \node (f12) [right=0.05ex of f11]{\includegraphics[scale=0.9, width=0.15\textwidth]{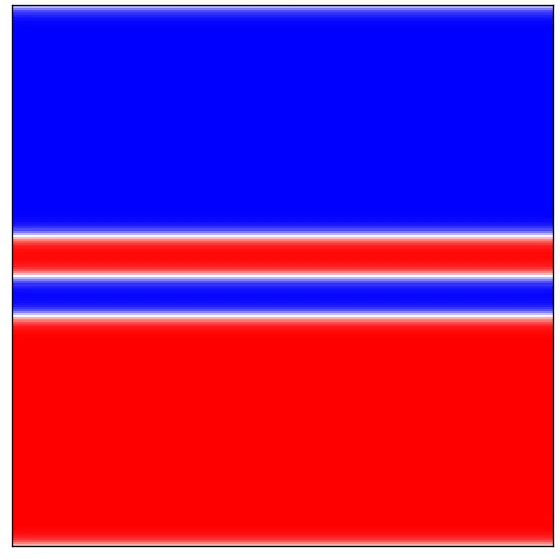}
   };
  \end{tikzpicture}
    \caption{Evolution by the diffuse Elastica flow. Phase-fields $u(\cdot,t)$ at times 0, 0.000414, 0.000920, 0.001839, 0.003679. The top row shows the standard diffuse and the bottom row the new diffuse approximation}
    \label{fig:new-inverted_1-b}
\end{figure}

The second example shows the Elastica flow for two touching circles with inverted phases across the horizontal axis (see Figure \ref{fig:new-inverted_1-b}).
Again we see an evolution into a stripe pattern, this time with multiple connected components for both phases.

\begin{figure}
 \centering
  \begin{tikzpicture}[trim left=-3.2em]
    \node (f1) at (0,0) {\includegraphics[scale=0.9, width=0.15\textwidth]{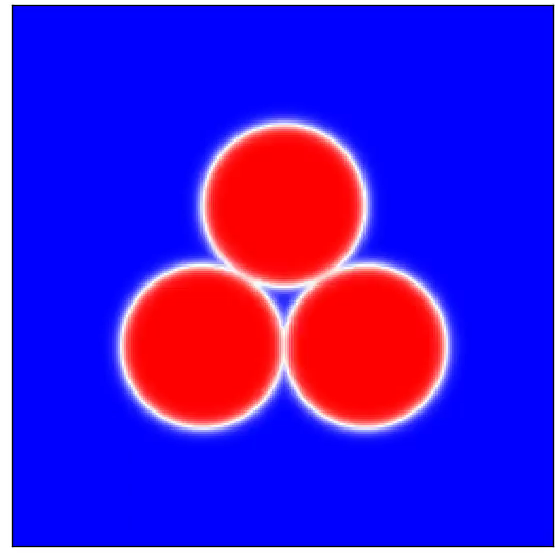}
    };
    \node (f2a) [right=0.1ex of f1]{\includegraphics[scale=0.9, width=0.15\textwidth]{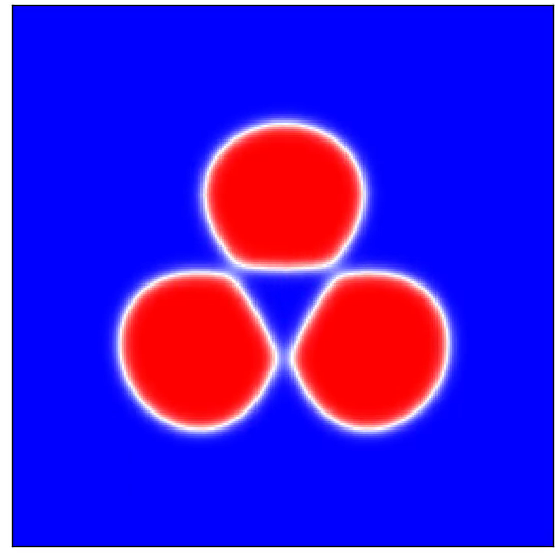}
    };
    \node (f2) [right=0.05ex of f2a] {\includegraphics[scale=0.9, width=0.15\textwidth]{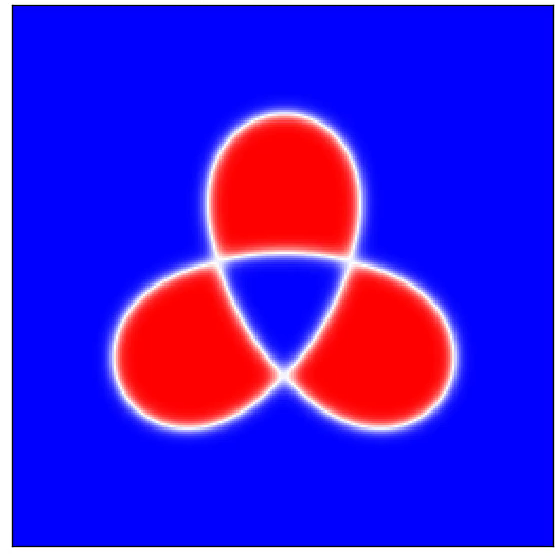}};
    \node (f3) [right=0.1ex of f2] {\includegraphics[scale=0.9, width=0.15\textwidth]{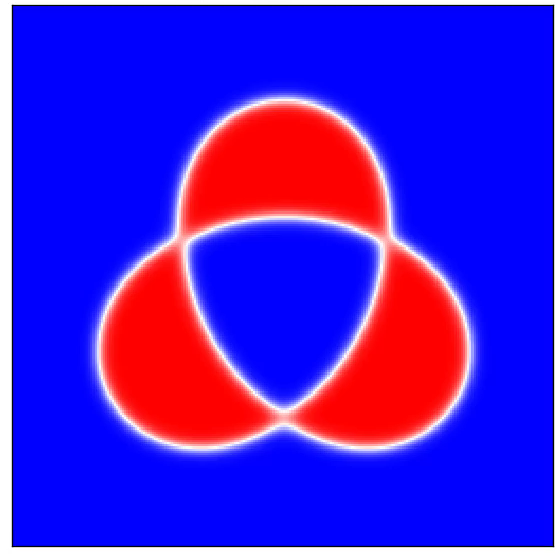}};
    \node (f4) [right=0.1ex of f3]{\includegraphics[scale=0.9, width=0.15\textwidth]{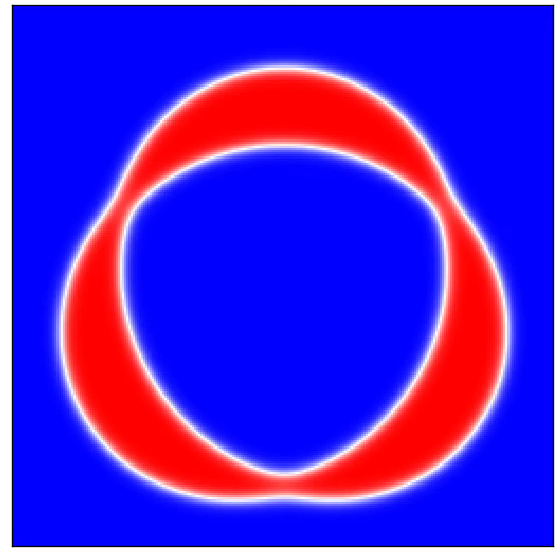}
    };
  \end{tikzpicture}
  \caption{Evolution by the new diffuse Elastica flow. Phase-fields $u(\cdot,t)$ at times 0, 0.000009, 0.000092, 0.000828, 0.004966.}
  \label{fig:new-threecircles-a}
\end{figure}

The third example shows the Elastica flow of three touching circles (see Figure \ref{fig:new-threecircles-a}), with the centers forming an equilateral triangle. A triangular symmetry is kept during the evolution and the circles deform in a ring-type structure that grows until boundary effects occur.
This behavior can also be observed for the standard diffuse approximation (not shown here).

\begin{figure}
 \centering
  \begin{tikzpicture}[trim left=-4.2em] 
    \node (f1) at (0,0)
    {\includegraphics[scale=0.9,width=0.2\textwidth]{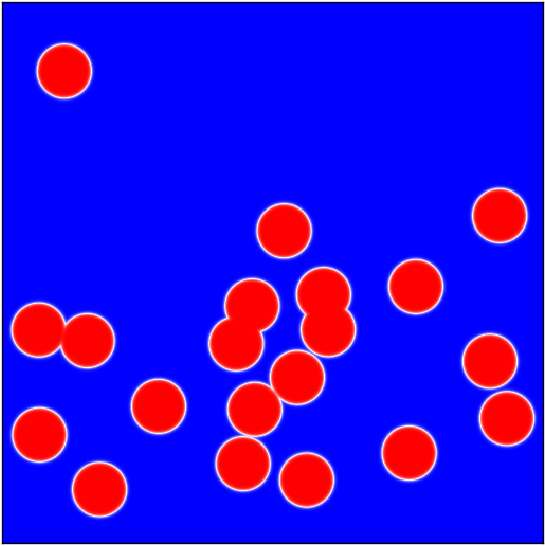}
    };
    \node (f2) [right=0.1ex of f1]
    {\includegraphics[scale=0.9,width=0.2\textwidth]{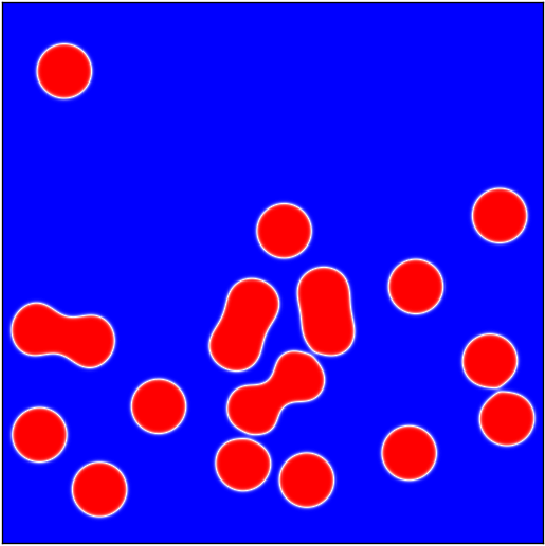}
    };
    \node (f3) [right=0.1ex of f2]
    {\includegraphics*[scale=0.9,width=0.2\textwidth]%
      {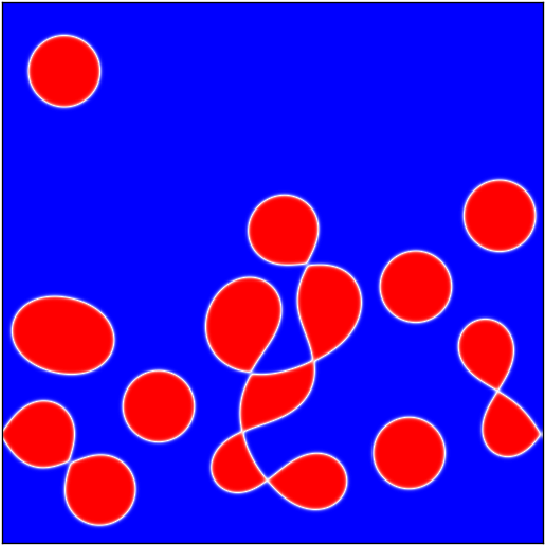}
    };
    \node (f4) [right=0.1ex of f3]
    {\includegraphics[scale=0.9,width=0.2\textwidth]{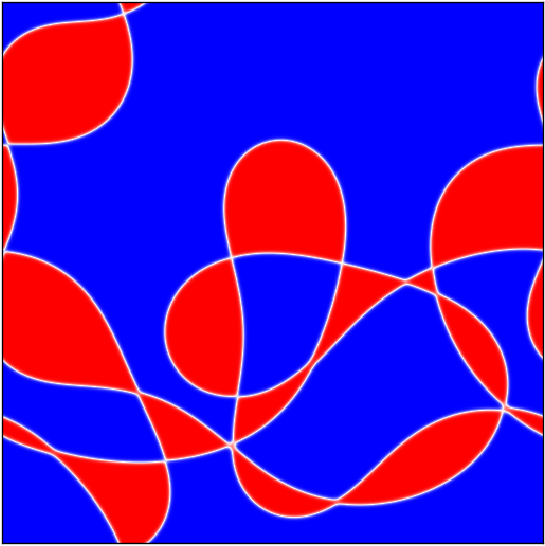}
    };
  \end{tikzpicture}
    \caption{The new diffuse Elastica flow. Evolution of the discrete
    phase-fields.} 
    \label{fig:new-noise}
\end{figure}

The fourth example shows the Elastica flow, starting from randomly positioned intersecting spheres. Some simulation snapshots are shown in Figure \ref{fig:new-noise}. We observe the occurrence of transversal intersections, which is a well-known behavior for diffuse approximations of the Elastica flow, see for example
\cite{EsedogluRaetzRoeger2014,BretinMasnouOudet2015}.
We recall that $N=512$ and $\eps=1/(2N)$ in this example, which leads to thinner transition regions in the corresponding pictures, compared to the other examples.

\begin{figure}
 \centering
  \begin{tikzpicture}[trim left=-5em]
   \node (f1) at (0,0) {\includegraphics[scale=0.9,width=0.2\textwidth]{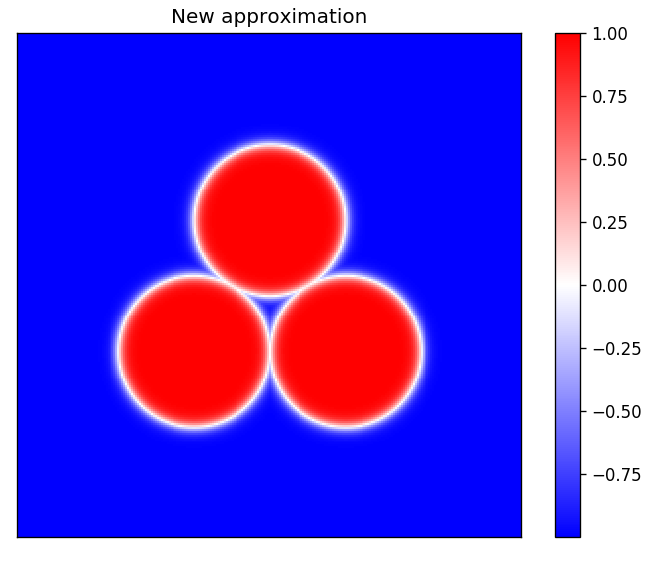}   };
   \node (f2) [right=0.1ex of f1]{\includegraphics[scale=0.9,width=0.2\textwidth]{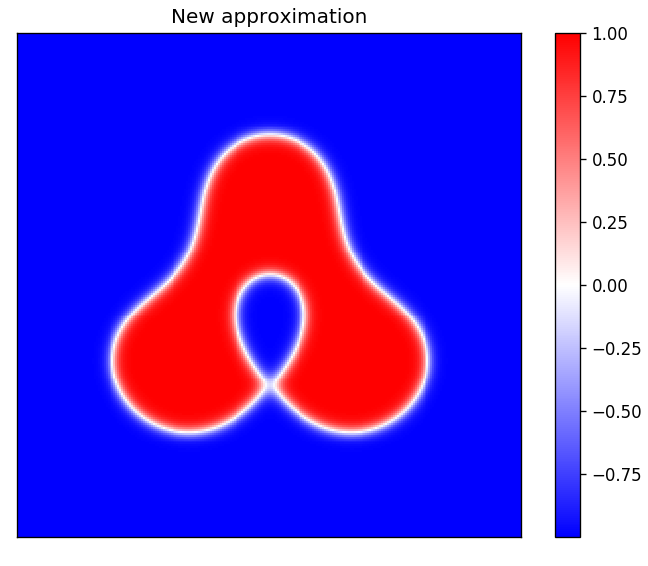}   };
   \node (f3) [right=0.1ex of f2] {\includegraphics[scale=0.9,width=0.2\textwidth]{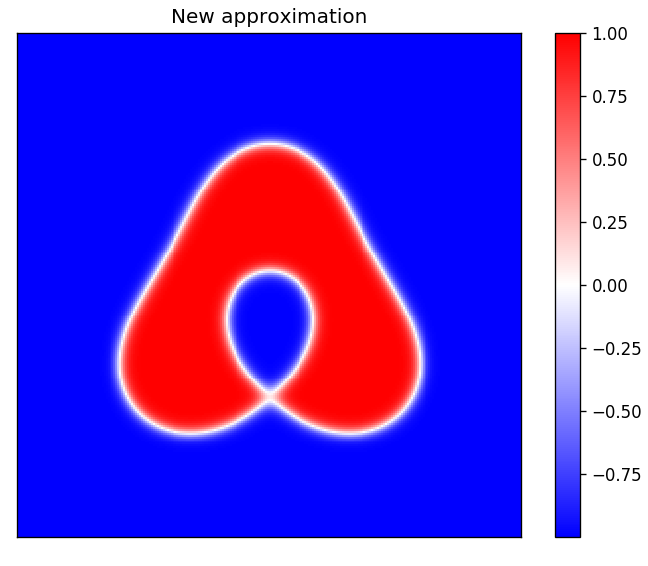}};
   \node (f4) [right=0.1ex of f3]{\includegraphics[scale=0.9,width=0.2\textwidth]{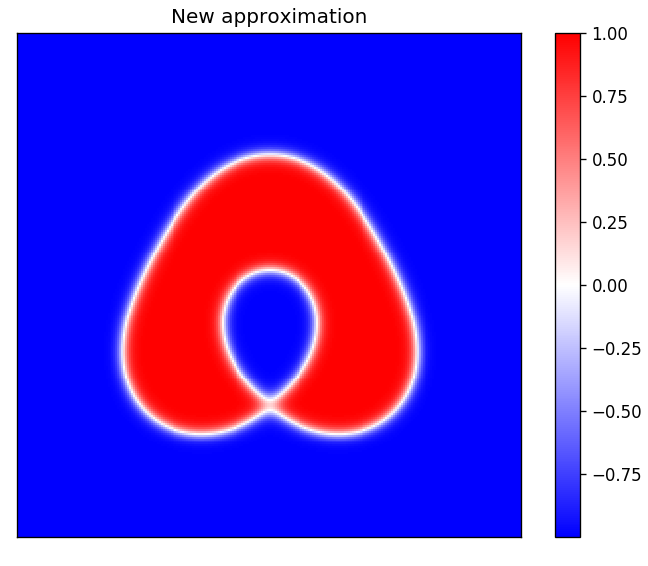} };
  \end{tikzpicture}
    \caption{The new diffuse approximation of Elastica flow with volume and perimeter constraint at times 0, 0.000019, 0.000074, 0.000139}
    \label{fig:new-constrainedthreecircles}
\end{figure}

We next consider the evolution by Elastica flow for three touching circles subject to constraints on volume and perimeter (see Figure \ref{fig:new-constrainedthreecircles}).
This time the triangular symmetry is broken. We again only present the evolution based on the new approximation.
For the standard approximation the processes are slightly slower but follow the same pattern.

\begin{figure}
 \centering
  \begin{tikzpicture}[trim left=-4.8em]
   \node (f1) at (0,0) {\includegraphics[scale=0.9,width=0.2\textwidth]{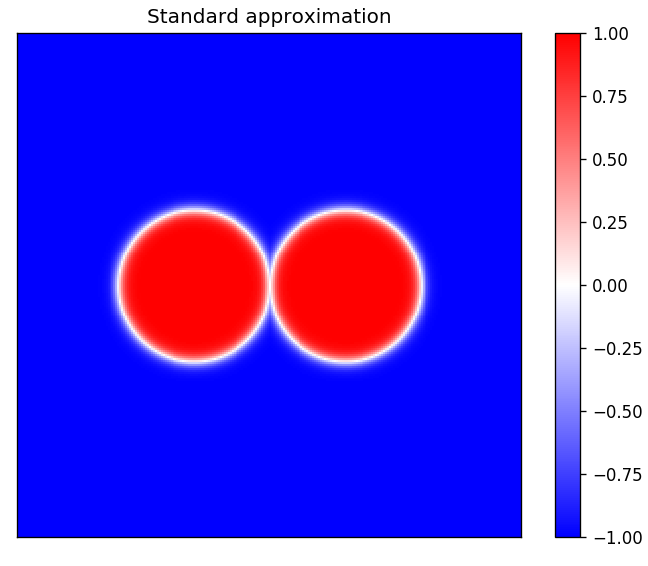}
   };
   \node (f2) [right=0.1ex of f1]{\includegraphics[scale=0.9,width=0.2\textwidth]{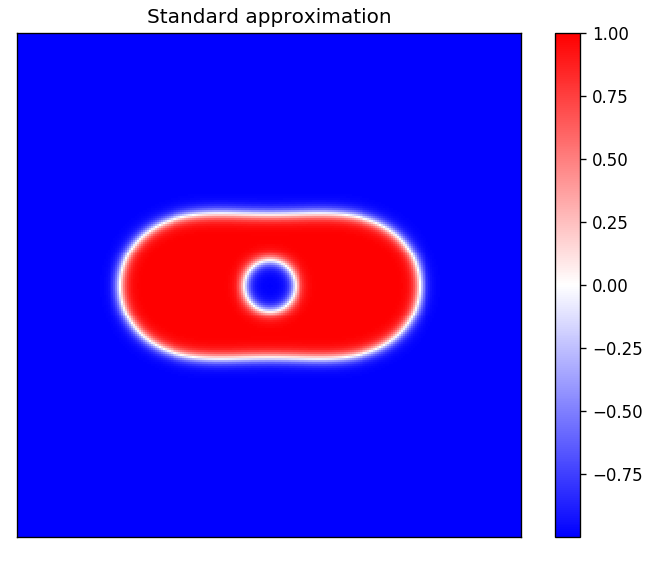}
   };
   \node (f3) [right=0.1ex of f2] {\includegraphics[scale=0.9,width=0.2\textwidth]{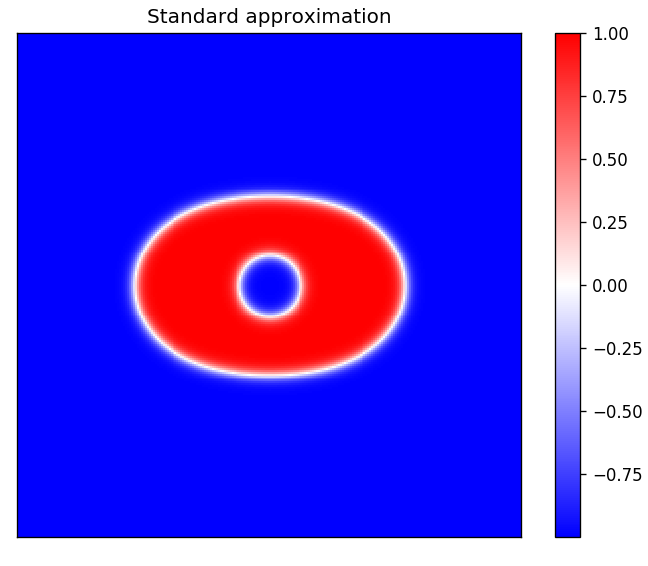}};
   \node (f4) [right=0.1ex of f3]{\includegraphics[scale=0.9,width=0.2\textwidth]{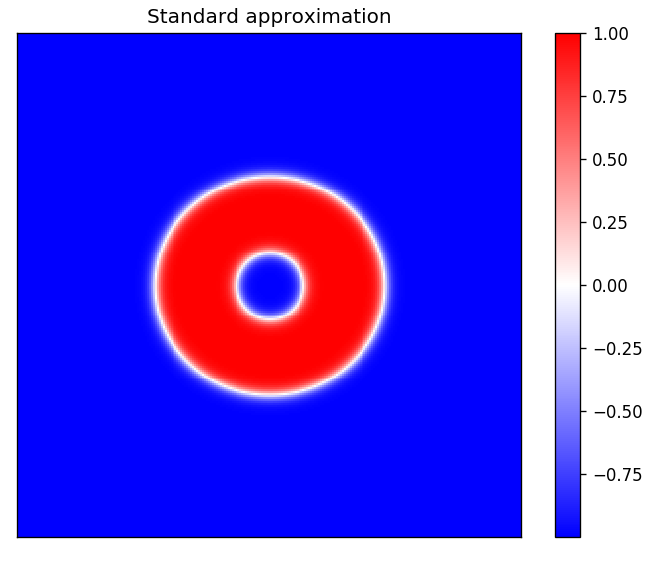}
   };
   \node (f5) [below=0.1ex of f1] {\includegraphics[scale=0.9,width=0.2\textwidth]{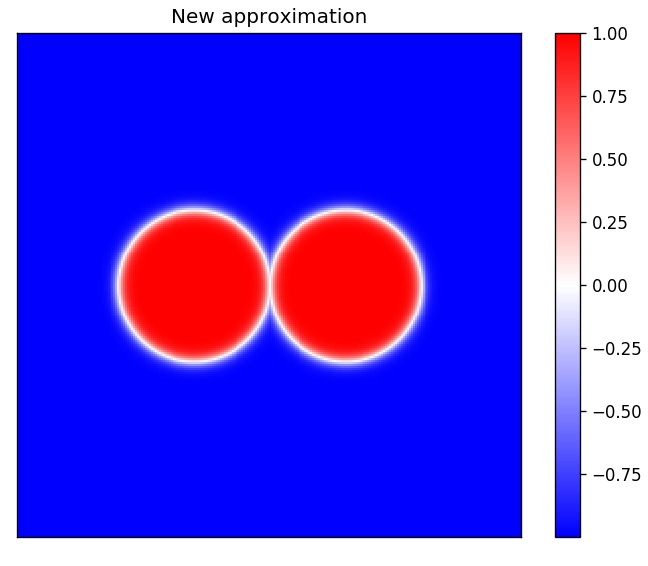}
   }; 
   \node (f6) [right=0.1ex of f5]{\includegraphics[scale=0.9,width=0.2\textwidth]{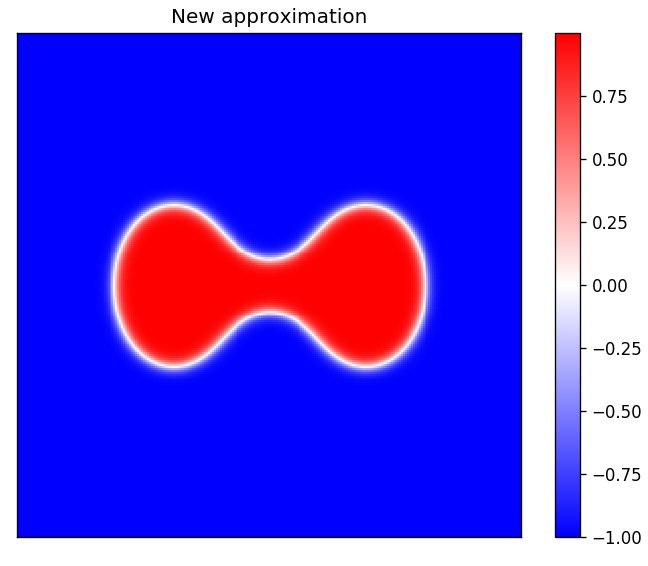}
   };
   \node (f7) [right=0.1ex of f6] {\includegraphics[scale=0.9,width=0.2\textwidth]{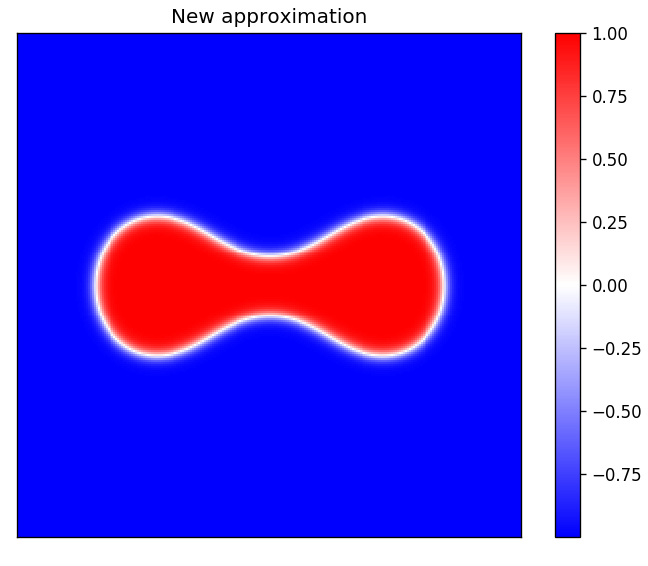}};
   \node (f8) [right=0.1ex of f7]{\includegraphics[scale=0.9,width=0.2\textwidth]{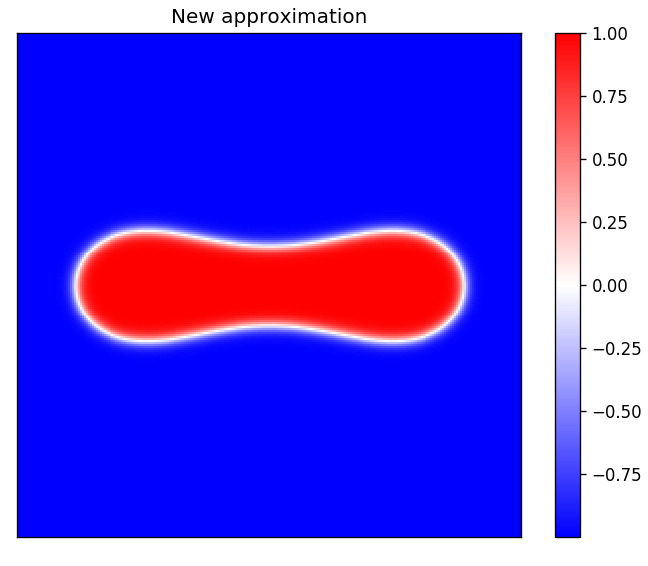}
   };
  \end{tikzpicture}
    \caption{Diffuse approximation of Elastica flow with volume and perimeter constraint at times 0, 0.000019, 0.000074, 0.000309. The top row shows the standard diffuse and the bottom row the new diffuse approximation.}
    \label{fig:new-constrainedtwocircles}
\end{figure}

The next example concerns the evolution by Elastica flow for two touching circles subject to volume and perimeter constraints (see Figure \ref{fig:new-constrainedtwocircles}). This time the evolution for the two diffuse approximations is different and approaches stationary states with different topology.
The respective structures can both be expected to be local minima of the constrained Elastica functional.

\begin{figure}\centering
 \begin{tikzpicture}[trim left=-8.5em,scale=1]
  \node (f1) at (0,0)
  {\includegraphics*[scale=0.9, width=0.45\textwidth]%
      {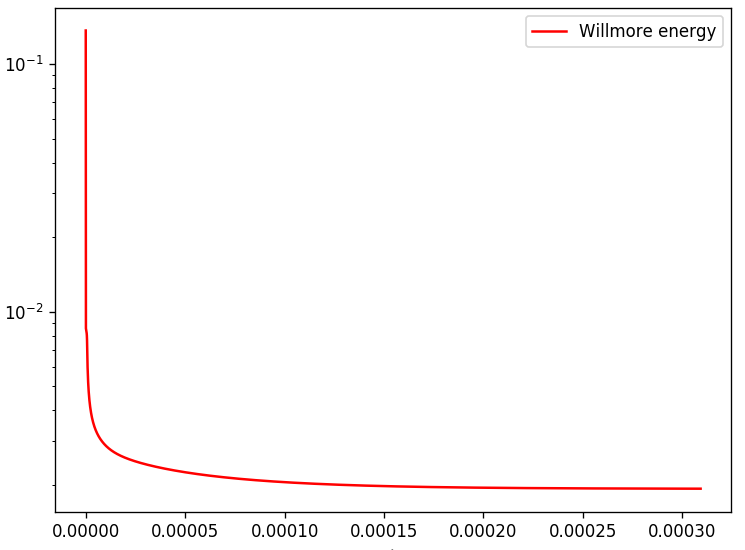}
  };
  \node (f3) [right=0.6ex of f1]
  {\includegraphics[scale=0.9, width=0.45\textwidth]{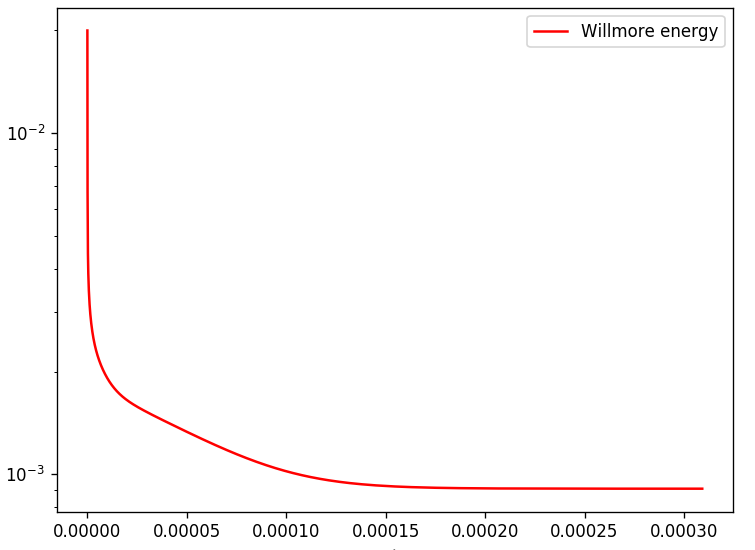}
  };
 \end{tikzpicture}
 \caption{Evolution of the discrete diffuse Willmore energy of the discrete approximations for the standard (left) and new approach (right).}
 \label{fig:energy_c2tc}
\end{figure}

An energy plot \ref{fig:energy_c2tc} shows that  the new approximation reaches a (numerically) stationary state with lower diffuse Willmore energy,
compared to the corresponding state in the standard evolution.
This example demonstrates that in general the evolution beyond the occurrence of non-smooth configurations does depend on the choice of the diffuse approximation.
The Willmore energies are computed for the discrete approximate solutions by discretizing the energies and using (inverse) Fourier transforms.

The last example is a three-dimensional numerical experiment of a cuboid evolution.
We consider the area and volume constrained Willmore flow for the new diffuse approximation.
We see an evolution to a biconcave discoid shape that resembles the shape of red blood cell is known appear as a local minimizer of constrained Willmore energies.

\begin{figure}\centering
   \begin{tikzpicture}[trim left=-4.2em, scale=1, every node/.style={scale=1}]
	   \node (f1) at (0,0) {\includegraphics[width=0.4\textwidth]{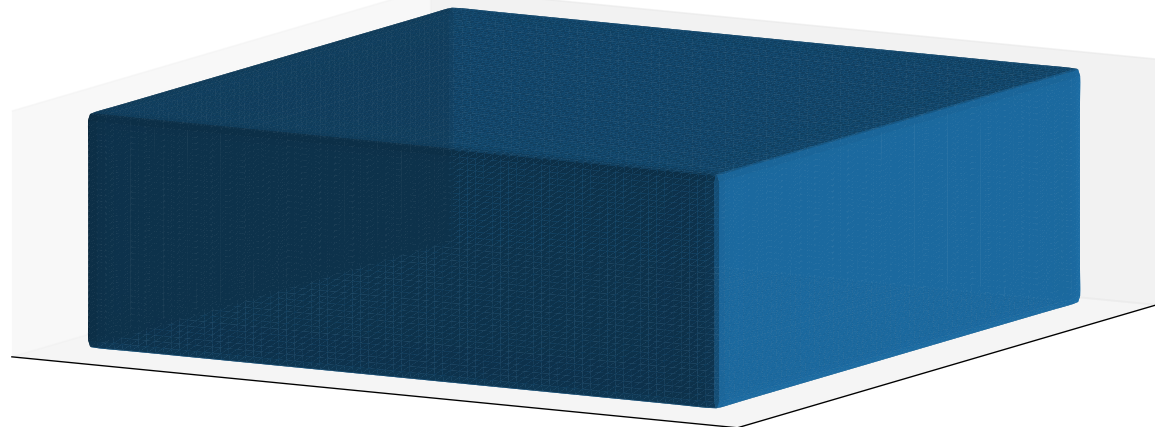}    };
     \node (f2) [right=0.1ex of f1]{\includegraphics[width=0.4\textwidth]{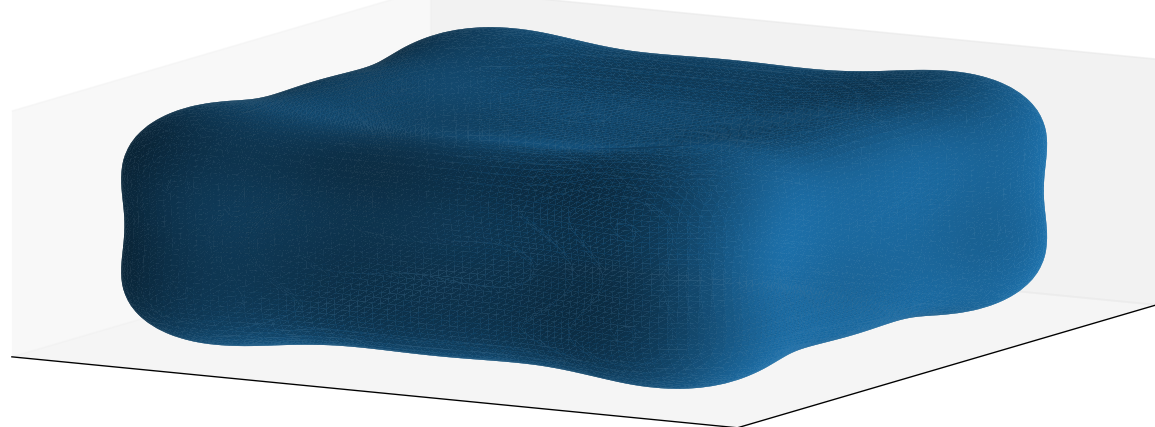}    };
    \node (f3) [below=0.1ex of f1] {\includegraphics[width=0.4\textwidth]{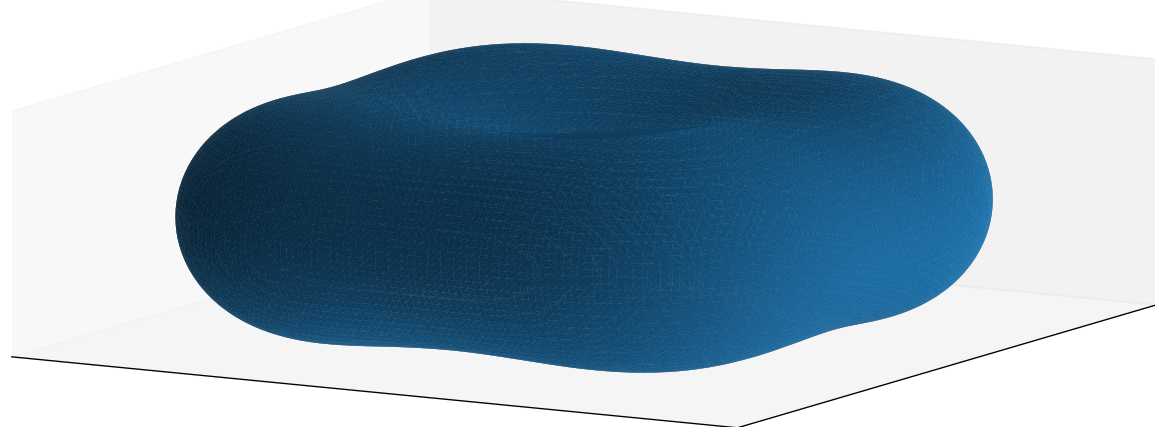}    };
    \node (f34) [right=0.05ex of f3]{};
    \node (f4) [right=0.1ex of f3]{\includegraphics[width=0.4\textwidth]{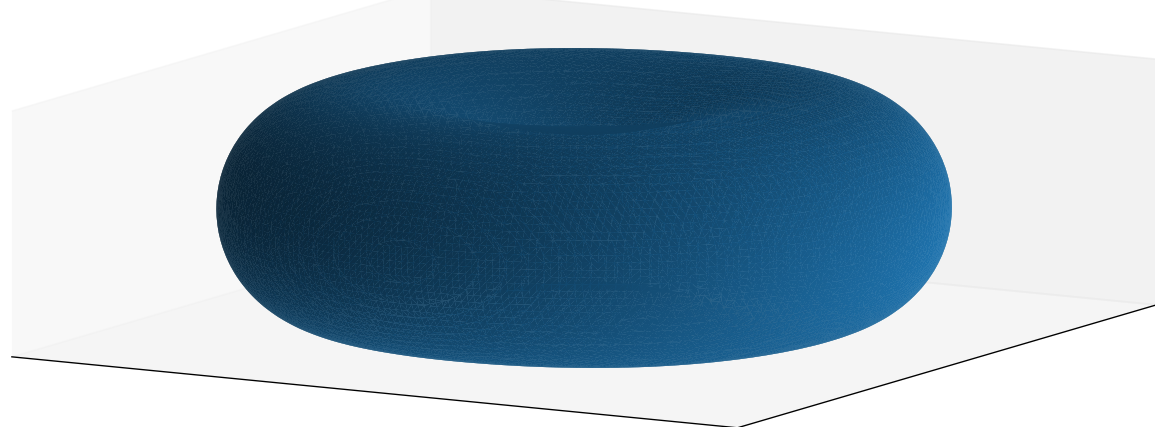}    };
    \node (f5) [below=of f34, ]{\includegraphics[width=0.8\textwidth]{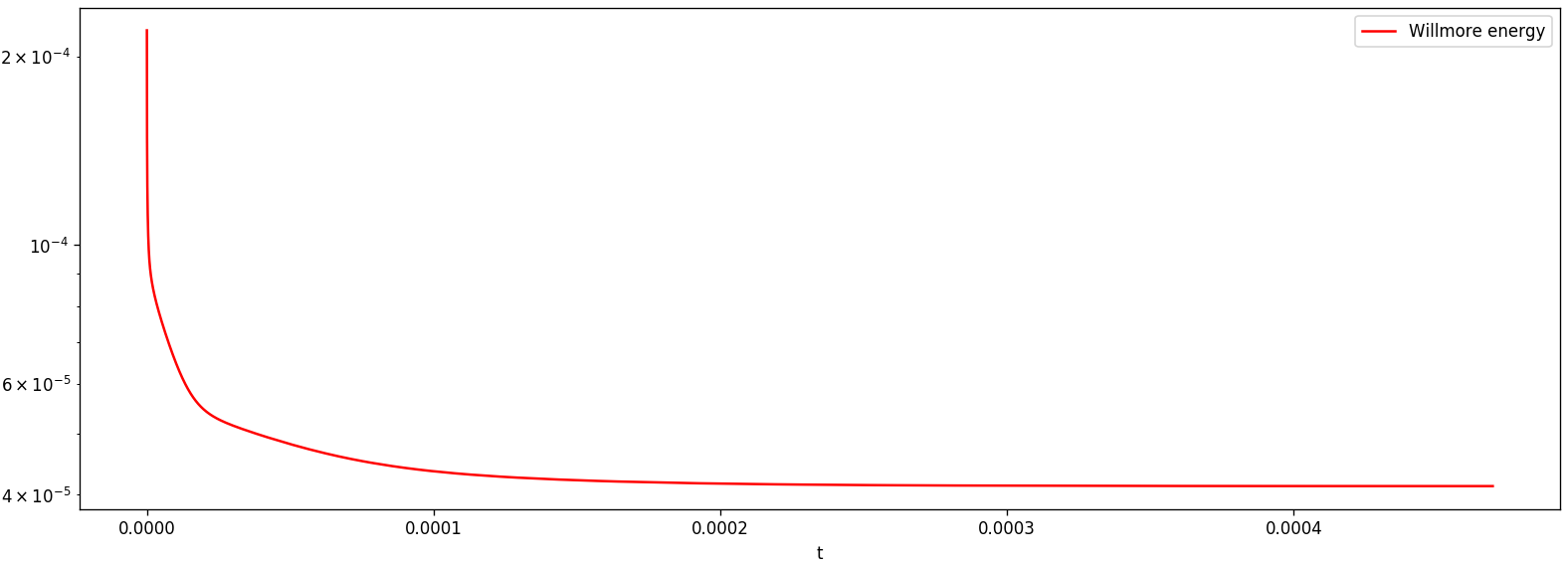}    };
   \end{tikzpicture}
 \caption{Evolution of the constrained Willmore flow for the new approximation. Snapshots and evolution of the discrete diffuse Willmore energy.}
 \label{fig:3dexp}
\end{figure}

\appendix
\section{A calculus lemma}
\begin{lemma}\label{lem:calc}
Let $f\in C^2(\R)$, $a<b$ and define for $0\leq\lambda\leq 1$
\begin{equation*}
  g(\lambda):= f\big((1-\lambda)a+\lambda b\big)-\big((1-\lambda)f(a)+\lambda f(b)\big).
\end{equation*}
Then
\begin{equation}
  |g(\lambda)|
  \leq C\|f''\|_{C^0([a,b])}\lambda(1-\lambda)(a-b)^2
  \label{eq:calc-lem}
\end{equation}
holds.
\end{lemma}
\begin{proof}
We have $g(0)=g(1)=0$. Taylor expansions give
\begin{align*}
  g(\lambda) &= f'(a)(b-a)\lambda + \frac{1}{2}f''(\xi_1)(b-a)^2\lambda^2 -\lambda\big(f(b)-f(a)\big)\\
  &=f'(a)(b-a)\lambda + \frac{1}{2}f''(\xi_1)(b-a)^2\lambda^2 -\lambda f'(a)(b-a)-\frac{1}{2}\lambda f''(\xi_2)(b-a)^2\\
  &= \frac{1}{2}(b-a)^2\lambda \Big(f''(\xi_1)\lambda -f''(\xi_2)\Big)
\end{align*}
and
\begin{align*}
  g(\lambda) &= f'(b)(b-a)(\lambda-1) + \frac{1}{2}f''(\xi_2)(b-a)^2(1-\lambda)^2 -(\lambda-1)\big(f(b)-f(a)\big)\\
  &= \frac{1}{2}(b-a)^2(1-\lambda) \Big(f''(\xi_3)(1-\lambda) -f''(\xi_4)\Big).
\end{align*}
Multiplying the first equality by $1-\lambda$, the second by $\lambda$ and adding up yields the desired estimate.
\end{proof}

\section{Exponential decay of profiles}
\begin{lemma}[Exponential decay of $\fraku_1,\frakv_1$]\label{lemma:expconvergenceforall}\ \\
Consider $\fraku_1$ as defined in Lemma \ref{lem:u1} and $\frakv_1$ as in equation \eqref{3eq:def-v1}.
Then $\fraku_1,\frakv_1,\frakv_1'$ decay at an exponential rate to $0$ at $\pm\infty$.
\end{lemma}
\begin{proof}
We observe that $\xi:=(\frakv_0,\frakv_0',\frakv_1,\frakv_1')$ is a solution to
the ODE system
\begin{align*}
	\xi_1' &= \xi_2,\\
	\xi_2' &= \xi_1-f^{-1}(\xi_1),\\
	\xi_3' &= \xi_4,\\
	\xi_4' &= -\xi_2+\frac{\sigma\xi_2}{f'(f^{-1}(\xi_1))^2}+\xi_3
	-\frac{\xi_3}{f'(f^{-1}(\xi_1))}.
\end{align*}
that converges to the stationary states $(\pm 1,0,0,0)$ with $r\to\pm\infty$.
We compute that these stationary points are
hyperbolic and deduce from the stable manifold theory, see for example
\cite{Perko1996}, that the solution approaches the
stationary states exponentially. The exponential convergence of $\fraku_1$
follows from the representation
\begin{align}
	\fraku_1=\frac{1}{f'(\fraku_0)}\big(-\sigma\fraku_0'+\frakv_1\big),
	\label{eq:u1formula}
\end{align}
and the vanishing of $\fraku_0'$ at $\pm\infty$ with an exponential rate,
see the corresponding statements in Theorem \ref{thm:1dminimizer}.
\end{proof}

\section{The solution operator $(-\partial^2+\Id)^{-1}$}
\begin{lemma}\label{lemma:expdecaysolution6}
We consider $J(z):=\frac{1}{2}e^{-|z|}$ for $z\in\R$. The convolution operator
$\frakA_0:w\longmapsto J\ast w$ has the following properties:
\begin{enumerate}
	\item $\frakA_0:L^2(\R)+L^\infty(\R)\longrightarrow L^2(\R)+L^\infty(\R)$
	is a well-defined mapping with $\frakA_0(L^2(\R))\subseteq H^2(\R)\cap L^\infty(\R)$,
	$\frakA_0(L^\infty(\R))\subseteq W^{2,\infty}(\R)$ and $\frakA_0(C^0_b(\R))
	\subseteq C^2_b(\R)$.
	\item If $w\in L^2(\R)$ the function $\frakA_0w$ is the unique solution in $H^2(\R)$ of
	\begin{align*}
		(-\partial^2+\Id)\frakA_0w = w\quad\text{ almost everywhere in }\R.
	\end{align*}
  \item If $w\in L^\infty(\R)$ and $\lim_{z\to\pm\infty}w(z)$ exist we also get
		\begin{align*}
			\lim\limits_{z\to\pm\infty}\frakA_0w(z)=\lim\limits_{z\to\pm\infty}w(z).
		\end{align*}
	\item If $w\in L^\infty(\R\times\omega)$ we define
	\begin{align*}
		\frakA_0w(z,x):=\int_{\R} J(z-\zeta)w(\zeta,x)d\zeta=[J\ast w(\cdot,x)](z).
	\end{align*}
	If $w\in C^{j_1}(\R;C^{j_2}(\omega))$, $j_1,j_2\in\N_0$ then
	$\frakA_0w\in C^{j_1+2}(\R;C^{j_2}(\omega))$.
	\item If $w\in X^{\mu,\Lambda}(\R;\Gamma)$ for some $\Lambda>0$ and $\mu\in(0,1)$
	we get $\frakA_0w,(\frakA_0w)',(\frakA_0w)''\in X^{\mu,\tilde \Lambda}(\R;\Gamma)$ for some $\tilde\Lambda=\tilde\Lambda(\Lambda,\mu)$.
\end{enumerate}
\end{lemma}
\begin{proof}
	The first two items follow from \cite{LiebLoss2001} and properties of the convolution (see for example \cite{Stei70}).

  The third and fourth claim follow from Lebesgue's dominated convergence theorem and standard theory of parameter dependent integrals.

  For the last item we estimate
  \begin{align*}
		|e^{\mu|z|}\frakA_0w(z,x)|&\leq\int_\R J(\zeta)e^{\mu|z|}|w(z-\zeta,y)|d\zeta
		\leq\Lambda\int_\R J(\zeta)e^{\mu|z|}e^{-\mu|z-\zeta|}d\zeta\\
		&\leq\Lambda\int_\R J(\zeta)e^{\mu|\zeta|}d\zeta
		=\frac{\Lambda}{2}\int_\R e^{-(1-\mu)|\zeta|}d\zeta
		=\frac{\Lambda}{1-\mu}.
	\end{align*}
	The estimate for $(\frakA_0w)'$ follows similarly since $J'\in L^\infty(\R)$ also decays exponentially.

  Finally, these properties also yield the decay of $(\frakA_0w)''=(\frakA_0w) - w$.
\end{proof}

\section{Expansion of the solutions operator}
\label{sec:A-expansion}
In this section we give some more details for the proof of Proposition \ref{pro:expansion-v}.

\begin{lemma}
\label{lemma:expdecaysolution1}
Let $\Lambda>0$, $\delta\in(0,1)$ and $\mu\in (0,1)$ be given.
There exists $\eps_0=\eps_0(\delta,\mu,\Gamma)>0$ with the following property:
Let $\eps\in(0,\eps_0)$, $\tilde u_\eps\in X^{\frac{\mu}{\eps},\Lambda}_\delta(\Omega)$ be given
and assume $\tilde v_\eps\in C^2(\Omega)\cap C^1(\overline \Omega)$ solves
\begin{align*}
	(-\eps^2\Delta+\Id)\tilde v_\eps&=\tilde u_\eps\quad\text{in}\quad\Omega\\
	\nabla\tilde v_\eps\cdot\nu_{\partial\Omega}&=0\quad\text{on}\quad\partial\Omega.
\end{align*}
Then we have $\tilde v_\eps\in X^{\frac{\mu}{\eps},\tilde\Lambda}_\delta(\Omega)$ for
$\tilde\Lambda=\frac{2}{1-\mu^2}\Lambda$.
\end{lemma}

\begin{proof}
We obtain that the function $Z:=e^{\frac{\mu}{\eps}|d_\delta|}\tilde v_\eps$ satisfies
\begin{align*}
	-\eps^2\Delta Z+2\mu\eps\nabla |d_\delta|\cdot \nabla Z+(1+\mu\eps
	\Delta |d_\delta|-\mu^2|\nabla |d_\delta||^2)Z
	=e^{\frac{\mu}{\eps}|d_\delta|}\tilde u_\eps\quad&\text{in}\quad\Omega\\
	\partial_{\nu_{\partial\Omega}}Z=0\quad\quad\quad\quad&\text{on}\quad\partial\Omega.
\end{align*}
Choose $\eps_0>0$ sufficiently small such that
\begin{equation*}
  \inf_\Omega \Big(1+\mu\eps \Delta |d_\delta|-\mu^2|\nabla |d_\delta||^2 \Big)\geq \frac{1-\mu^2}{2}>0.
\end{equation*}
Assume that $M:=\max_{\overline \Omega} Z>\tilde\Lambda:=\frac{2}{1-\mu^2}\Lambda$ and observe that in $\{Z>\tilde\Lambda\}$
\begin{equation}
  -\eps^2\Delta (Z-\tilde\Lambda)+2\mu\eps\nabla |d_\delta|\cdot \nabla (Z-\tilde\Lambda)+\frac{1-\mu^2}{2}(Z-\tilde\Lambda)
  \leq e^{\frac{\mu}{\eps}|d_\delta|}\tilde u_\eps -\Lambda
  \leq 0
  \label{eq:Hopf-max}
\end{equation}
since $\tilde u_\eps\in X^{\frac{\mu}{\eps},\Lambda}_\delta(\Omega)$.
If the maximum of $Z$ is attained at a point $x_0\in\partial\Omega$ we choose an open ball
$B\subset\Omega\cap \{Z>\tilde\Lambda\}$ with $\overline{B}\cap\partial\Omega=\{x_0\}$.
We deduce from the Hopf Lemma and $\nabla Z\cdot\nu_{\Omega}=0$ that $Z=M$ in $B$ holds.
This implies that the maximum of $Z$ is always attained in $\Omega$, which yields by \eqref{eq:Hopf-max} that $Z\leq\tilde\Lambda$.
Similarly we obtain $-Z\leq\tilde\Lambda$.
\end{proof}

The next lemma shows a {\em quasi-locality} of the operator $-\eps^2\Delta+\Id$ if applied to functions that are exponentially close to $\pm 1$ away from the interface.
\begin{lemma}\label{lemma:expdecaysolution4}
Let $\delta\in(0,1)$, $\mu,\Lambda>0$ and a cut-off
function $\eta_1$ as in Assumption \ref{ass:modified_distances} be given.
Then there exists $\eps_0=\eps(\delta,\Gamma,\eta_1)$
such that for all $\eps\in(0,\eps_0)$ and any $w_\eps\in C^2(\omega)$ the following
property holds: If for all $x\in\{|d|\geq 3\delta\}$
\begin{align*}
	|\eps\nabla w_\eps(x)|\leq\Lambda e^{-\frac{\mu}{\eps}|d_\delta(x)|}\quad\text{and}\quad
	|w_\eps(x)-\sgn(d(x))|\leq\Lambda e^{-\frac{\mu}{\eps}|d_\delta(x)|},
 \end{align*}
then there exists $R_\eps\in X^{\frac{\mu}{\eps},\Lambda}_\delta(\Omega)$ such that
\begin{align}
	&(-\eps^2\Delta+\Id)(\eta_\delta w_\eps+(1-\eta_\delta)\sgn(d))
  \nonumber\\
	&\qquad\qquad =\eta_\delta\cdot(-\eps^2\Delta+\Id)w_\eps+(1-\eta_\delta)\sgn(d)+\chi_{\{|d|\geq 3\delta\}}R_\eps.
  \label{eq:app:C3}
\end{align}
\end{lemma}
\begin{proof}
We calculate
\begin{align*}
	&(-\eps^2\Delta+\Id)(\eta_\delta w_\eps +(1-\eta_\delta)\sgn(d))\\
	&=\eta_\delta\cdot(-\eps^2\Delta+\Id)w_\eps+(1-\eta_\delta)\sgn(d)
  -2\eps^2\nabla w_\eps\nabla\eta_\delta
	-\eps^2\Delta\eta_\delta\cdot(w_\eps-\sgn(d)).
\end{align*}
For the last two terms we obtain
\begin{align*}
	2\eps^2\big|\nabla w_\eps\nabla\eta_\delta\big|\leq
	\frac{4\eps^2}{\delta}\chi_{\{|d|\geq 3\delta\}}\big|\nabla w_\eps\big|
	\leq\frac{4\eps\Lambda}{\delta}\chi_{\{|d|\geq 3\delta\}}e^{-\frac{\mu}{\eps}|d_\delta|}
\end{align*}
and
\begin{align*}
	|\eps^2\Delta\eta_\delta\cdot(w_\eps-\sgn(d))|\leq
	\frac{\eps^2\Lambda C(\eta_1,\Gamma)}{\delta^2}\chi_{\{|d|\geq 3\delta\}}e^{-\frac{\mu}{\eps}|d_\delta|}.
\end{align*}
Choosing $\eps_0\leq\min\{\frac{1}{8},\big(2C(\eta_1,\Gamma)\big)^{-\frac{1}{2}}\}\delta$ yields the claim.
\end{proof}

We need a corresponding statement for functions that are defined in terms of the inner variables.
\begin{lemma}\label{lemma:expdecaysolution5}
There exists $\eps_0>0$ such that for all $\eps\in(0,\eps_0)$ and any
$w\in C^2(\R\times\omega)$ the following holds:

Assume
\begin{align*}
	w-\sgn,\, \partial_z w,\, \nabla_x w \in X^{\mu,\Lambda}(\R;\Gamma)
\end{align*}
and define $\wei:=w\circ\Psi_\eps^{-1}\in C^2(\omega)$.
Then there exist $R_\eps^w\in X^{\frac{\tilde\mu}{\eps},2\Lambda}_\delta(\Omega)$ such that \eqref{eq:app:C3} holds for $w_\eps=\wei$.
\end{lemma}

\begin{proof}
We observe that $|d_\delta(x)|\leq \frac{2}{3}|d(x)|$ in $\{|d|\geq 3\delta\}$ and deduce in $\{3\delta\leq |d|\leq 5\delta\}$
\begin{align*}
	\eps\big|\nabla \wei\big|(x) &\leq \big|
	\partial_z w(z,x)\big|+ \eps \big|(\nabla_x w)(z,x)\big|
	\leq 2\Lambda e^{-\mu|z|}
  \leq 2\Lambda e^{-\frac{3\mu d_\delta(x)}{2\eps}}
\end{align*}
and
\begin{align*}
	|\wei(x)-\sgn(d(x))|&\leq
	\Lambda e^{-\mu|z|}\leq \Lambda e^{-\frac{3\mu d_\delta(x)}{2\eps}}.
\end{align*}
The claim then follows from Lemma \eqref{lemma:expdecaysolution4}.
\end{proof}

We finally can prove Proposition \ref{pro:expansion-v}.
\begin{proof}[Proof of Proposition \ref{pro:expansion-v}]
Let $v_\eps=\calA_\eps u_\eps$ and consider the profile functions $v_0,v_1,v_2$ as defined in \eqref{eq:veps-0}-\eqref{eq:veps-2}.
In this proof $\Lambda$ may change from line to line but will always be independent of $\eps$.

Let $u_0,u_1,u_2\in X^{\mu,\Lambda}(\R;\Gamma)$, $0<\mu<1$.
We first observe from Lemma \ref{lemma:expdecaysolution6} and $(\frakA_0\sgn)
=\sgn(z)(1-e^{-|z|})$ that $v_0$ inherits the exponential decay to $\pm 1$ from $u_0$, since
\begin{align*}
	v_0-\sgn=\frakA_0(u_0-\sgn)
	+\big(\frakA_0\sgn-\sgn\big)
	\in X^{\mu,\Lambda}(\R;\Gamma).
\end{align*}
From  Lemma \ref{lemma:expdecaysolution6} we also get $\partial_z v_0\in X^{\mu,\Lambda}
(\R;\Gamma)$ and $v_0\in C^2(\R\times\omega)$ is $C^4$-regular with respect to $x$.

The same arguments apply to the next order. This yields
\begin{align*}
	v_1=\frakA_0u_1+H \frakA_0\partial_z v_0\in
	C^2(\R\times\omega)\cap X^{\mu,\Lambda}(\R;\Gamma),
\end{align*}
and $\partial_z v_1\in X^{\mu,\Lambda}(\R;\Gamma)$.

Similarly, we have
\begin{align*}
	v_2=\frakA_0u_2+H \frakA_0\partial_z v_1+\frakA_0(\Delta-z|\FF|^2\partial_z)v_0
	\in X^{\frac{10}{11}\mu,\Lambda}(\R;\Gamma),
\end{align*}
with $v_2\in C^2(\R\times\omega)$, and $\partial_z v_2\in X^{\frac{10}{11}\mu,\Lambda}(\R;\Gamma)$.

We then obtain from \eqref{eq:uei-vei}
\begin{align}
	(-\eps^2\Delta+\Id)\vei&=\uei + \eps^3R_\eps^v,
  \label{eq:pf3.7-1}
\end{align}
with
\begin{align*}
  R_\eps^v\circ \Psi_\eps &= -|z|^2R_\eps^H\partial_zv_0 + H\partial_zv_2+(\Delta-z|\FF|^2\partial_z)v_1 \\
  &\qquad - \eps |z|^2R_\eps^H\partial_zv_1 +\eps(\Delta-z|\FF|^2\partial_z)v_2 - \eps^2|z|^2R_\eps^H\partial_zv_2.
\end{align*}

From the properties of the profile functions and their derivatives we can conclude
$R_\eps^v\circ \Psi_\eps\in X^{\frac{10}{11}\mu,\Lambda}(\R;\Gamma)$.
Since $\phi_1(z)\leq \frac{9}{10}z$ for all $z\geq 0$ by Assumption \ref{ass:modified_distances} we deduce $|d_\delta|\leq \frac{9}{10}|d|$ and therefore obtain $R_\eps^v\in X^{\mu,\Lambda}_\delta(\Omega)$.

The previous observations show that $w:=v_0+\eps v_1 +\eps^2 v_2$ fulfills the assumptions of Lemma \ref{lemma:expdecaysolution5} with $\mu$ replaced by $\frac{10}{11}\mu$.
Applying the Lemma to $\wei=\vei$ and using \eqref{eq:pf3.7-1} we therefore obtain for some $R_\eps^w\in X^{\frac{15\mu}{11\eps},\Lambda}_\delta(\Omega)$
\begin{align*}
	(-\eps^2\Delta+\Id)\big(&\eta_\delta\vei +(1-\eta_\delta)\sgn(d)\big)\\
	&=\eta_\delta\cdot(-\eps^2\Delta+\Id)\vei +(1-\eta_\delta)\sgn(d)
	+ \chi_{\{|d|\geq 3\delta\}}R_\eps^w\\
	&=\eta_\delta\uei+\eps^3\eta_\delta R_\eps^v+(1-\eta_\delta)\sgn(d)
	+ \chi_{\{|d|\geq 3\delta\}}R_\eps^w\\
	&=u_\eps+\eps^3 R_\eps,
\end{align*}
where $R_\eps\in X^{\frac{\mu}{\eps},\Lambda}_\delta(\Omega)$ due to the stronger exponential decay of $R_\eps^w$.

Since $\nabla\big(\eta_\delta\vei +(1-\eta_\delta)\sgn(d)\big)\cdot\nu_\Omega=0$ at $\partial\Omega$ we deduce from Lemma \ref{lemma:expdecaysolution1}
\begin{align*}
	v_\eps=\eta_\delta \vei +(1-\eta_\delta)\sgn(d)
	+\eps^3{\tilde R}_\eps,
\end{align*}
with ${\tilde R}_\eps\in X^{\frac{\mu}{\eps},\tilde\Lambda}_\delta(\Omega)$.
\end{proof}

\hbadness=1800
\bibliographystyle{vancouver}

\end{document}